\documentclass[review]{elsarticle}

\usepackage{hyperref}

\usepackage[english]{babel}

\usepackage{graphicx}
\usepackage{subcaption}

\usepackage{amsmath,amssymb,amsthm,amsfonts}

\usepackage{dsfont}

\usepackage{blindtext}

\journal{Spatial Statistics}

\begin{document}

\begin{frontmatter}

\title{Information and Complexity Analysis of Spatial Data}

\author[1]{José M. Angulo\corref{mycorrespondingauthor}}
\cortext[mycorrespondingauthor]{Corresponding author}
\ead{jmangulo@ugr.es}
\author[1]{Francisco J. Esquivel}
\ead{jesquivel@ugr.es}
\author[1]{Ana E. Madrid}
\ead{anaesther@ugr.es}
\author[1]{Francisco J. Alonso}
\ead{falonso@ugr.es}
\address[1]{University of Granada, Department of Statistics and Operations Research, Campus Fuente Nueva s/n, 18071 Granada, Spain}

\begin{abstract}
Information Theory provides a fundamental basis for analysis, and for a variety of subsequent methodological approaches, in relation to uncertainty quantification. The transversal character of concepts and derived results justifies its omnipresence in scientific research, in almost every area of knowledge, particularly in Physics, Communications, Geosciences, Life Sciences, etc. Information-theoretic aspects underlie modern developments on complexity and risk.   
A proper use and exploitation of structural characteristics inherent to spatial data motivates, according to the purpose, special considerations in this context.       
      
In this paper, some relevant approaches introduced regarding the informational analysis of spatial data, related aspects concerning complexity analysis, and, in particular, implications in relation to the  structural assessment of  multifractal point patterns, are reviewed under a conceptually connective evolutionary perspective.

\end{abstract}


\end{frontmatter}


\section{Introduction}
\label{section:introduction}

Since the appearance of the seminal paper by Shannon (1948), there has been an immense literature devoted to the foundations and developments of Information Theory, a vast discipline with transversal implications in every field of knowledge. Besides many other proposals, Shannon entropy still remains as a benchmark concept from which, in particular, a number of generalizations have been introduced, with Rényi entropy (Rényi 1961) constituting the most representative one under preservation of extensivity for independent systems. 

In the last decades, `complexity' has become a multiconceptual aspect of increasing interest for research. In particular, from the point of view of quantification regarding structural characteristics of a random system, Information Theory provides a meaningful and solid support for objective formalization of complexity measures. In the context of Geography, Batty, Morphet, Masucci and Stanilov (2014), extending Batty's (1974) original approach for `spatial entropy' from a partitioning-based formal relation between the discrete and continuous Shannon entropy versions, proposed a measure of `complexity' under a restricted notion of departure from equilibrium (equiprobability). In parallel, in the ambit of Physics, a different conception of `complexity' was introduced (see, for instance, Huberman and Hogg 1986), as departure from both equilibrium and singularity (i.e. degeneracy to a single state with probability 1). Under this approach, the discrete product-type formulation by López-Ruiz, Mancini and Calbet (1986) for a balanced assessment of information and disequilibrium factors, and its exponential-entropy version for continuous probability distributions proposed by Catalán, Garay and López-Ruiz (2002), lead to the formulation of a generalized complexity measure by López-Ruiz, Nagy, Romera and Sañudo (2009), a two-parameter family which, in the end, quantifies in an exponential scale the variations of Rényi entropy with respect to changes in the deformation parameter, for a given continuous probability distribution. In fact, based on the notion by Campbell (1966) of exponential (Shannon and Rényi) entropy as a measure of extent (`diversity index') of a probability distribution, either continuous or discrete, we justify the significance of the two-parameter generalized complexity measure also for the discrete case, giving then a proper complementary interpretation as a diversity ratio for both scenarios. In particular, we show that Batty et al.'s approach to complexity can then be embedded as a specific, partially degenerated case, into the generalized complexity measure family, for some fixed parameter values.   

A step further in the perspective presented in this paper is concerned with systems characterized by multifractal measures (see, for instance, the accessible presentation by Harte 2001). In some general sense, multifractality has been generally identified and referred to in the literature as a form of complexity.  As a reference object of the multifractal formalism, the generalized Rényi dimensions (and the singularity spectrum, its pair Legendre transform), defined in terms of the scaling behaviour of Rényi entropy, constitute a key instrument for practical analysis in many fields of application, and particularly in Geophysics. In this respect, a direct limiting connection between the above mentioned two-parameter generalized complexity measure and the increments of the generalized dimension curve was formally established in Angulo and Esquivel (2014), thus highlighting a precise interpretation of the variational properties of the latter for quantitative complexity assessment in the multifractal domain. In Esquivel, Alonso and Angulo (2017), the complementary  informativeness of the maps of relative increments and, in particular, of the curve of derivatives of generalized dimensions was further justified. 

Beyond and complementarily with the informational `global' comparison of two possible distributions on a given system directly from their marginal entropies, divergence measures are formulated with the aim of quantifying their probabilistic structural coherence by a `local' (state-by-state) examination. Historically, as before,  Kullback-Leibler (1951) divergence (as the counterpart for Shannon entropy)  and its generalization by Rényi (1961) divergence (correspondingly for Rényi entropy) constitute, among many other available definitions, reference concepts in the construction and multidisciplinary applications of Information Theory. In a parallel scheme to the discussion outlined above, we analyse their implications in the context of complexity and multifractality. In this direction, Romera, Sen and Nagy (2011), adopting a similar product-type approach, now based on Rényi divergence, proposed a two-parameter generalized relative complexity measure, for `local' assessment of complexity coherence between two continuous probability distribution. From a formal and conceptual statement of a divergence-based `relative diversity index', we also justify the significance of the corresponding discrete version and, as before, give a complementary related interpretation for both cases. Finally, based on a natural formulation of generalized relative dimensions regarding the scaling behaviour of Rényi divergence in a multifractal context, we also establish the direct limiting relation between the corresponding increments and the generalized relative complexity family, with a subsequent interpretation and derivation of related variational tools for practical (multifractal) relative complexity assessment. 

In summary, the aim of this paper is to present, in a schematic synthesis, a perspective properly connecting, formally and conceptually,  all the above mentioned aspects, from information (section \ref{section:information-entropy-divergence}) to complexity (section \ref{section:complexity}) and, in particular, to multifractality (section \ref{section:multifractality}), with a parallel  reference to measures for `global'  and  (comparative) `local' structural assessment. Through the study of a real seismic series (section \ref{section:example}), focused on the temporal distribution and magnitudes of registered events during different phases in the circumstance of a volcanic episode, we demonstrate, in particular, the complementary usefulness of generalized dimensions and generalized relative dimensions for analysis and interpretation. In final section \ref{section:further-remarks-conclusion}, a concise conclusion is given, with an added reference to some further aspects and alternative approaches.

\section{Information Entropy and Divergence}
\label{section:information-entropy-divergence}

In this section, formal and conceptual aspects focused on the definition of Shannon and Rényi entropies, as well as Kullback-Leibler and Rényi divergences, are introduced. These constitute the basis for the construction of product-type generalized complexity and generalized relative complexity measures, respectively, as discussed in section \ref{section:complexity}, as well as for the definition of generalized dimensions  and generalized relative dimensions in the multifractal domain, see section \ref{section:multifractality}. (Related notions based on Tsallis entropy and divergence, as an alternative form of a deformation-parameter generalization, are referred in section \ref{section:further-remarks-conclusion}.) In analogy to Campbell diversity index, a divergence-based relative diversity index is proposed, both of which provide a specific interpretation for product complexity and relative complexity measures addressed in section \ref{section:complexity}. Batty's approach to spatial entropy, which connects with the discussion of discrete vs. continuous entropy versions, is also introduced here, as a former step to a related notion of complexity associated with the information difference, referred too in section \ref{section:complexity}.

\subsection{Entropy measures (Hartley, Shannon, Rényi)}
\label{subsection:entropy}

\subsubsection{Discrete probability distributions}
\label{subsubsection:entropy-discrete}

\begin{paragraph}{`Information content' (Hartley entropy) -- Hartley (1928)}
\medskip

For a random selection of one element within a finite set $A$ of cardinality $n$, an appropriate measure of information under monotonicity, additivity and normalization requirements is given by
$$ \log (n). $$
An extended interpretation, from this germinal notion introduced and justified by Hartley (1928) in the context of telecommunications, leads to the following general formulation: For an event with probability $p$ of occurrence, the quantity   
$$ - \log (p) .$$
represents an appropriate measure of information, as before, under monotonicity, additivity and normalization requirements. This can be viewed as the `information content' provided by the occurrence of that particular event.  

Since the measure is intended for quantitative comparison purposes, and hence it can be more generally defined except for any specified multiplicative constant, an arbitrary logarithmic base can be adopted (Hartley 1928, p.540). In what follows, the natural logarithm `$\log_{e}$', denoted as `$\ln$', is used in all the formulations.
\end{paragraph}

\begin{paragraph}{Shannon entropy -- Shannon (1948)}
\medskip
For a discrete probability distribution $\bar{p} = (p_1,\dots,p_n)$,  the (Shannon) entropy (or `information entropy')  is defined as 
$$ H(\bar{p}) := - \sum_{i=1}^n p_i \ln (p_i) = - E[\ln (\bar{p})] .$$
Thus, Shannon (1948) entropy is interpreted as the expected information content provided by the system realization (outcome generation) based on $\bar{p}$, i.e. the $\bar{p}$-mean Hartley (extended) entropy.

Among others, two well-known basic properties of Shannon entropy are: 
\begin{itemize}
\item Minimum and maximum values of $H$:
$$  H_{min}= 0 \qquad H_{max}= \ln(n), $$
with $H_{min}$ being related to degenerate systems concentrating the probability mass in only one of the possible states, e.g. $\bar{p} = (1, 0, \dots, 0)$, and $H_{max}$ being reached only in the case of equiprobability, i.e. $\bar{p} = (\frac{1}{n}, \dots, \frac{1}{n}) =: \left[\frac{1}{n}\right]$ (from now on, although the latter should be properly written as $H_{max}(n):= H\left(\left[\frac{1}{n}\right]\right)$  in reference to the specific number $n$ of possible states considered, this argument is omitted whenever it is implicitly understood).   
  \item `Extensivity' (or additivity property): The (joint) Shannon entropy of a system composed by two independent systems, which are assumed to be characterized by two independent respective discrete probability distributions, is equal to the sum of their individual (marginal) Shannon entropies. 
\end{itemize}

For illustration, Figure \ref{figure:Shannon_entropy} represents in a simplex plot the Shannon entropy values based on a ternary system, characterized by a discrete probability distribution with $n=3$ possible elementary events.

\begin{figure}
            \includegraphics[height=8.0cm]{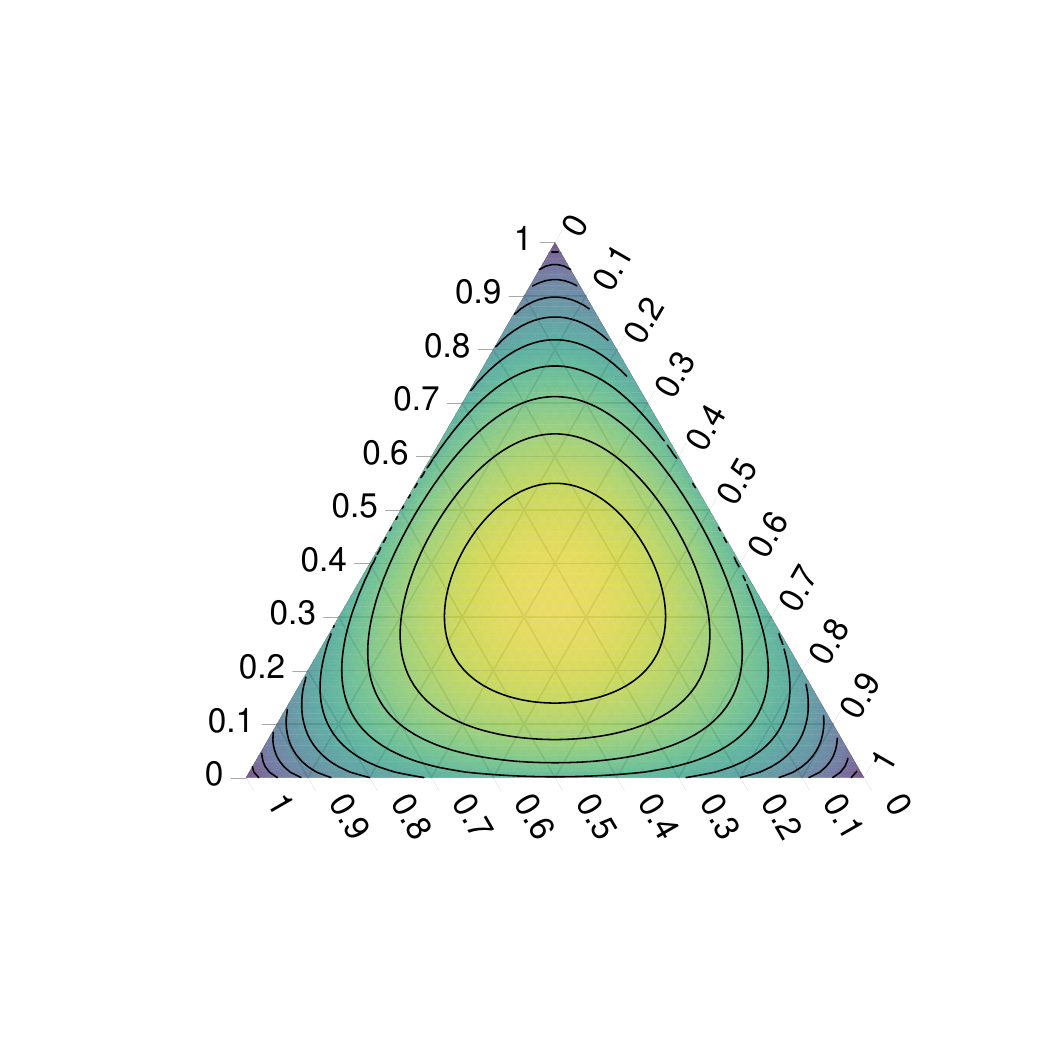}\hspace*{1mm}
            \includegraphics[trim=0mm -19mm 133mm -5mm, clip, height=8.0cm]{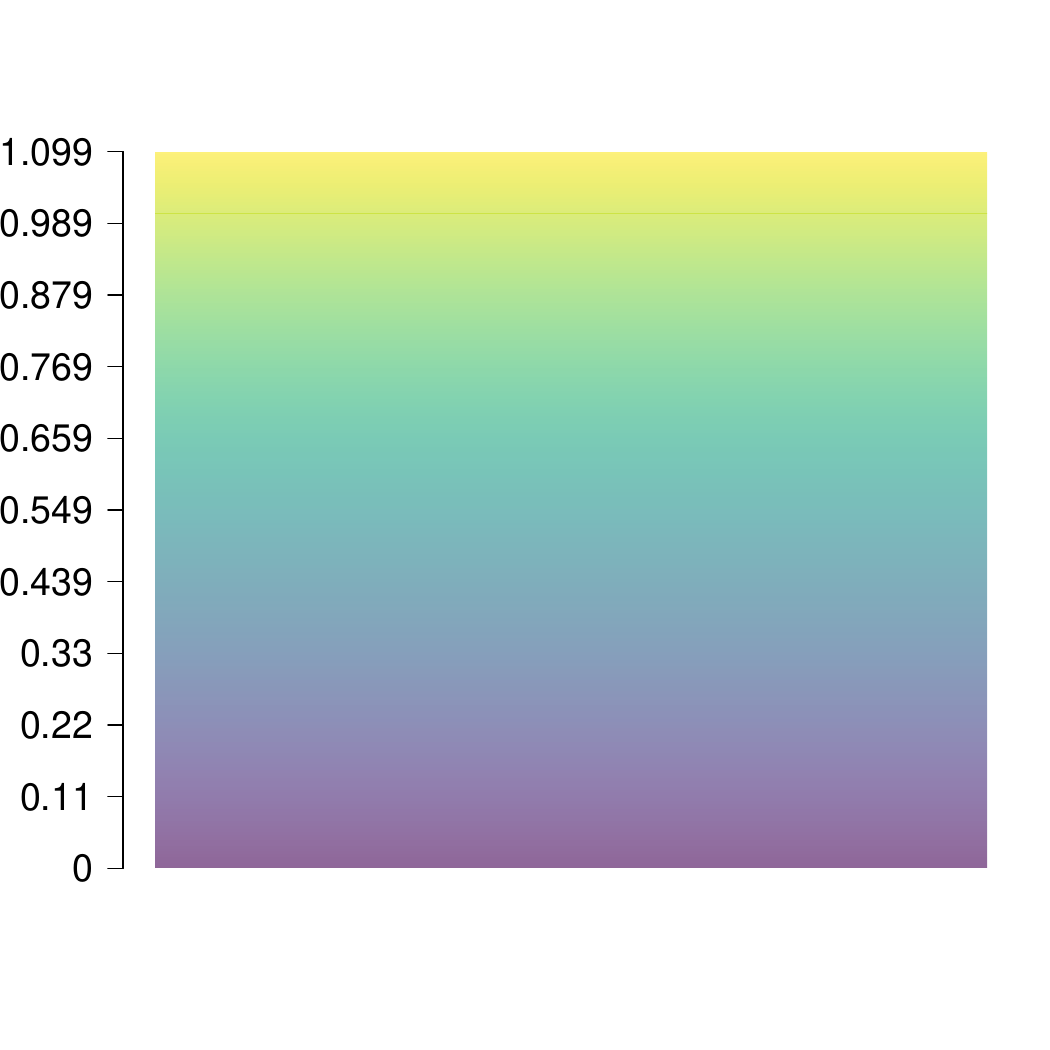}

\caption{Shannon entropy  for $\bar{p}= (p_1, p_2, p_3)$, with values varying from 0, the minimum uncertainty at the degenerate distributions corresponding to the vertices, to $\ln(3)$, the maximum uncertainty associated with the equiprobability central point.} 
\label{figure:Shannon_entropy}
\end{figure}
           
\end{paragraph}

\begin{paragraph}{R\'{e}nyi entropy of order $q$ -- Rényi (1961)}

\medskip
For a discrete probability distribution $\bar{p} = (p_1,\dots,p_n)$, the (Rényi) entropy of order $q$ is defined as 
$$ H_q(\bar{p}) :=  \frac{1}{1-q} \ln \left(\sum_{i=1}^n p_i^q\right) = \frac{1}{1-q}\ln\left( E[\bar{p}^{q-1}]\right) \qquad (q\neq1). $$
Historically, Rényi (1961) entropy perhaps constitutes the best-known and most-used generalization of Shannon entropy, being formulated in terms of a `deformation parameter' $q$ whose meaning is related to the $q$-power distortion effect derived on the reference probability distribution $\bar{p}$:
$$
\bar{p}= (p_1, \dots, p_n) \qquad  \longrightarrow  \qquad  \bar{p}^{q,\ast}= (p_1^{q,\ast}, \dots, p_n^{q, \ast}), \qquad p_i^{q,\ast} = \frac{p_i^{q}}{\sum_j p_j^{q}}.
$$
Figure \ref{figure:power_distortion} displays the paths followed by the ternary probability distribution $\bar{p}^{q,\ast}= (p_1^{q,\ast}, p_2^{q,\ast}, p_3^{q, \ast})$ as $q$ tends from 1 to $\infty$ (left plot) or to 0 (right plot), assuming different starting distributions. In particular, it can be observed that the center of the triangle (three-dimensional simplex), the central points in the three edge segments (degenerate marginal two-dimensional sub-simplexes), and the three vertices (degenerate marginal one-dimensional sub-simplexes) represent stationary and/or absorbing points depending on the case.

\begin{figure}
\includegraphics[trim=50mm 100mm 50mm 100mm, clip, height=5.0cm]{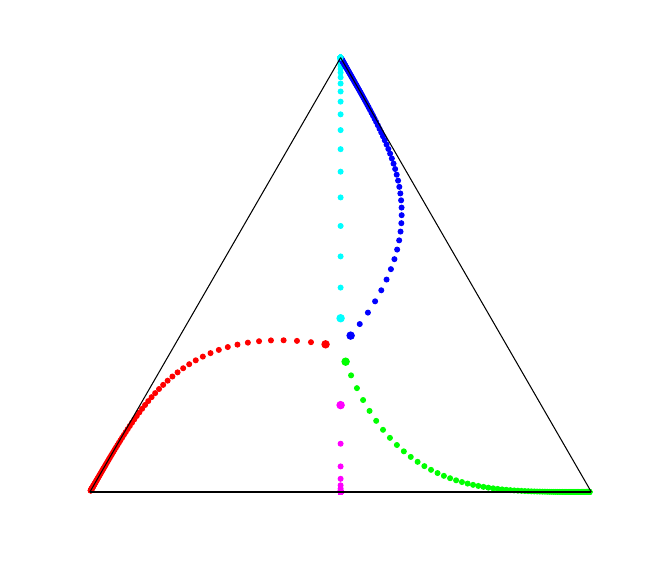}
\includegraphics[trim=50mm 100mm 50mm 100mm, clip, height=5.0cm]{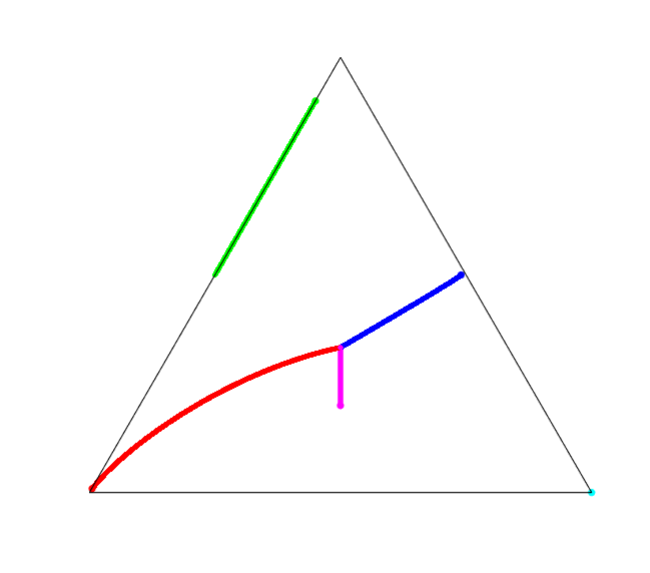}
\caption{Paths followed by $\bar{p}^{q,\ast}= (p_1^{q,\ast}, p_2^{q,\ast}, p_3^{q, \ast})$  as $q$ tends from 1 to $\infty$ (left plot) or to 0 (right plot), from different starting distributions. In the first case, power distortion makes any non-equiprobable distribution to move towards the edges and/or vertices; in particular, middle edge points attract those distributions having two dominant equiprobable states. In the second case, conversely, any distribution without probability degeneracy for any of the states moves by power distortion towards the simplex center representing equiprobability; partial degeneracy with null probability for only one state makes the distribution tend towards the partial equiprobability represented by the middle point of the corresponding edge; complete degeneracy in one state, represented by the vertices, results in invariance.} 
\label{figure:power_distortion}
\end{figure}

Among others, as before, some well-known basic properties of Rényi entropy are: 
 
\begin{itemize}
\medskip
\item Shannon entropy is the limiting case of Rényi entropy as $q\rightarrow 1$.
\medskip
  \item Minimum and maximum values of $H_q$:
$$ {H_q}_{min}= 0 \qquad {H_q}_{max}= \ln(n). $$

\item Rényi entropy satisfies `extensivity' for independent systems. 
\end{itemize}

Plots in Figure \ref{figure:Renyi_entropy_0.5-2-5-100}  show the values of Rényi entropy for a ternary  probability distribution under different specifications for the deformation parameter,  namely $q=0.5$, 2,  5 and 100.

\begin{figure} \hspace*{-11mm}
 \includegraphics[height=5.6cm]{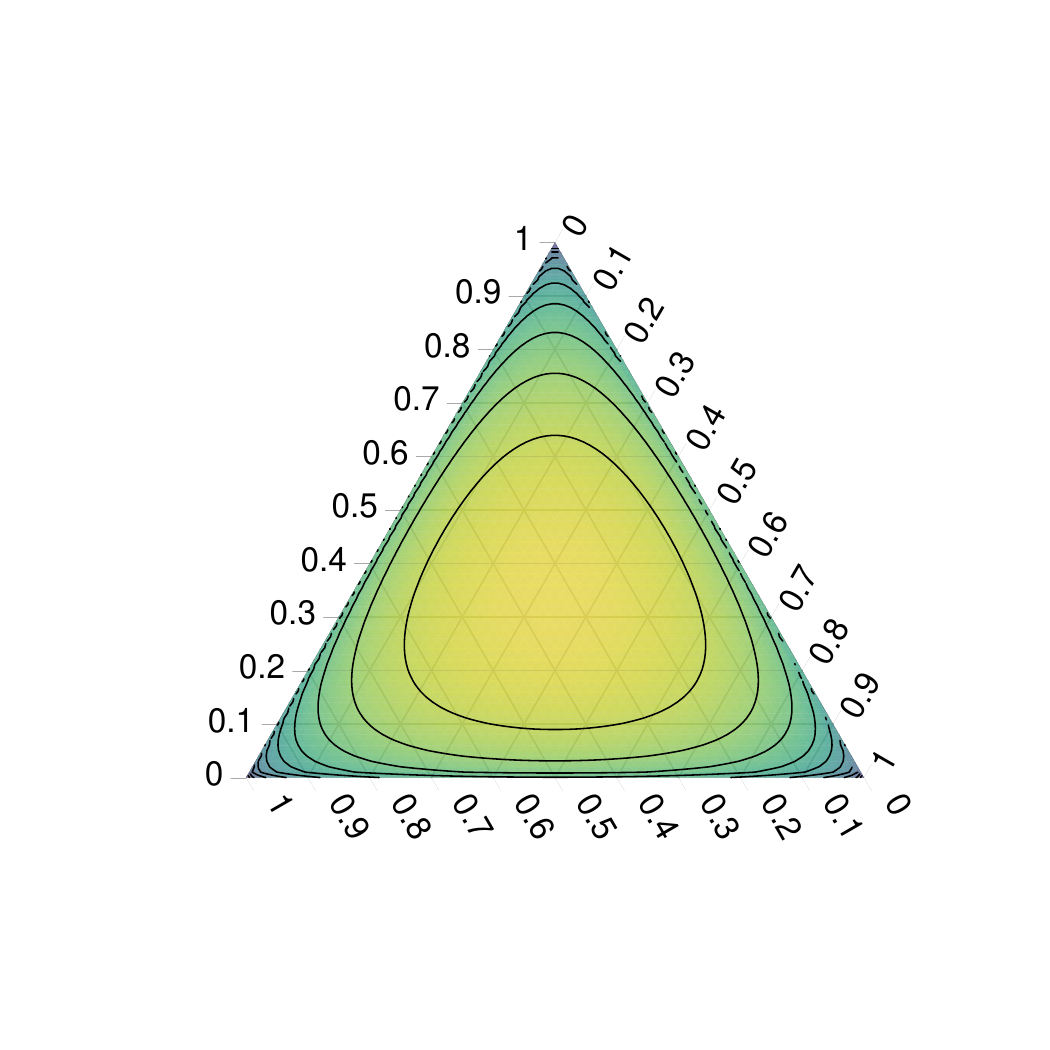}\hspace*{-5mm}
            \includegraphics[trim=0mm -19mm 133mm -5mm, clip, height=5.6cm]{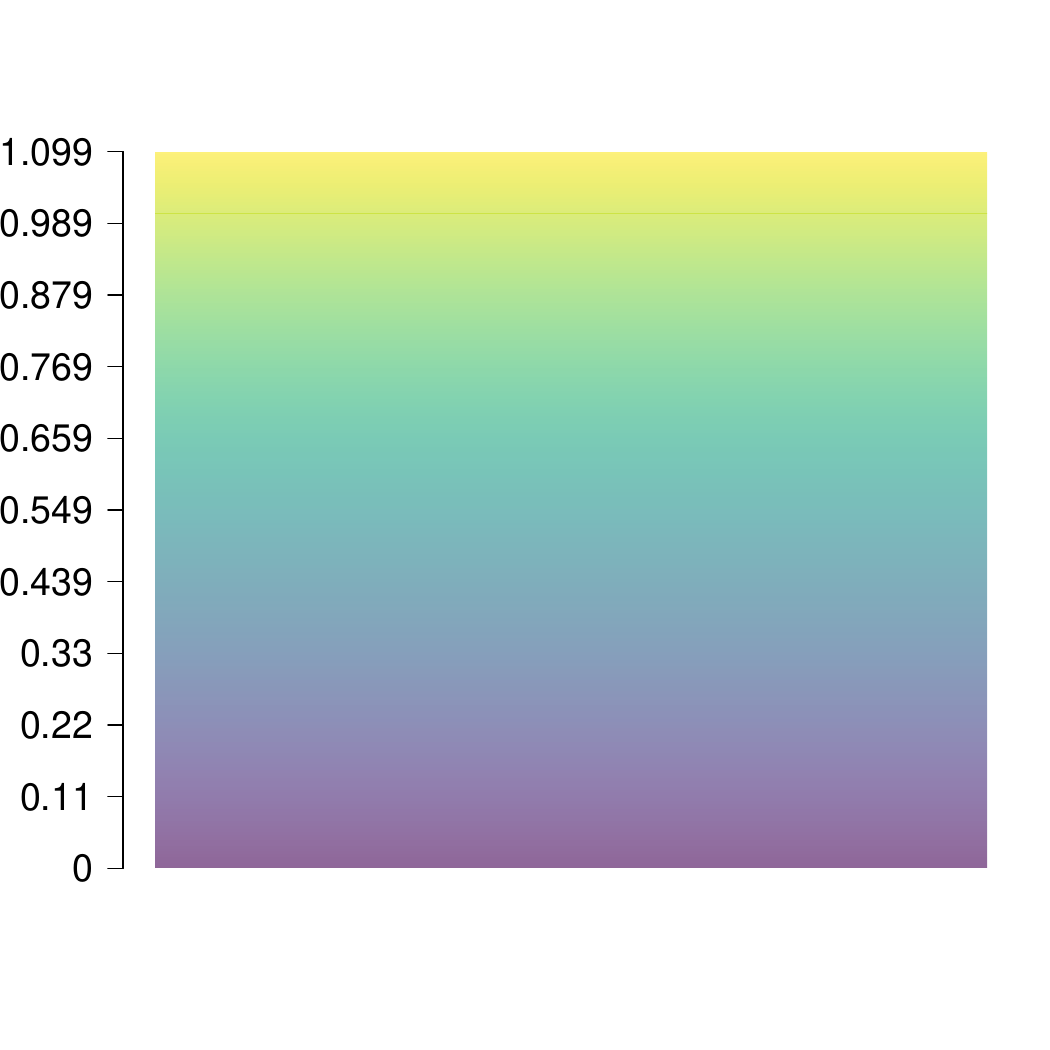}
\hspace*{-1mm}
 \includegraphics[height=5.6cm]{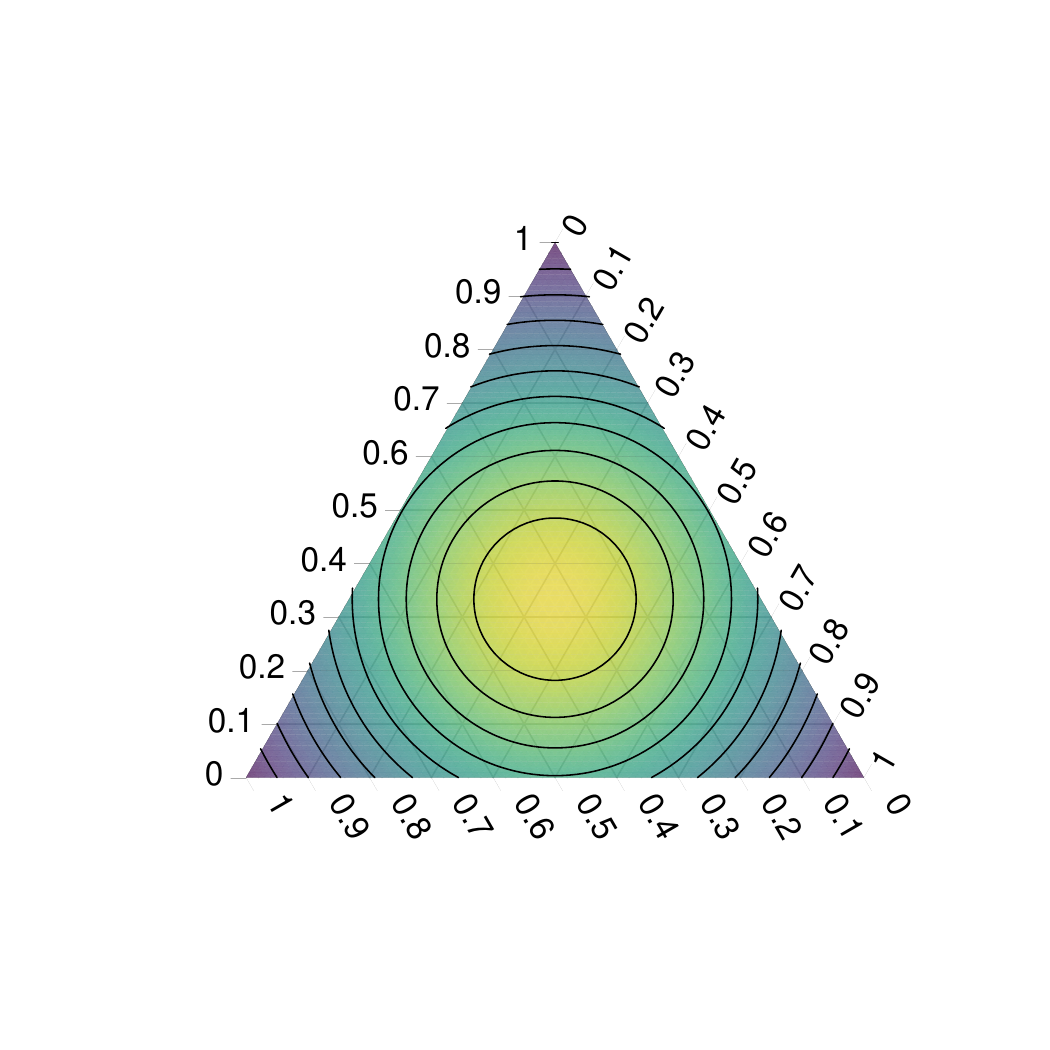}\hspace*{-5mm}
            \includegraphics[trim=0mm -19mm 133mm -5mm, clip, height=5.6cm]{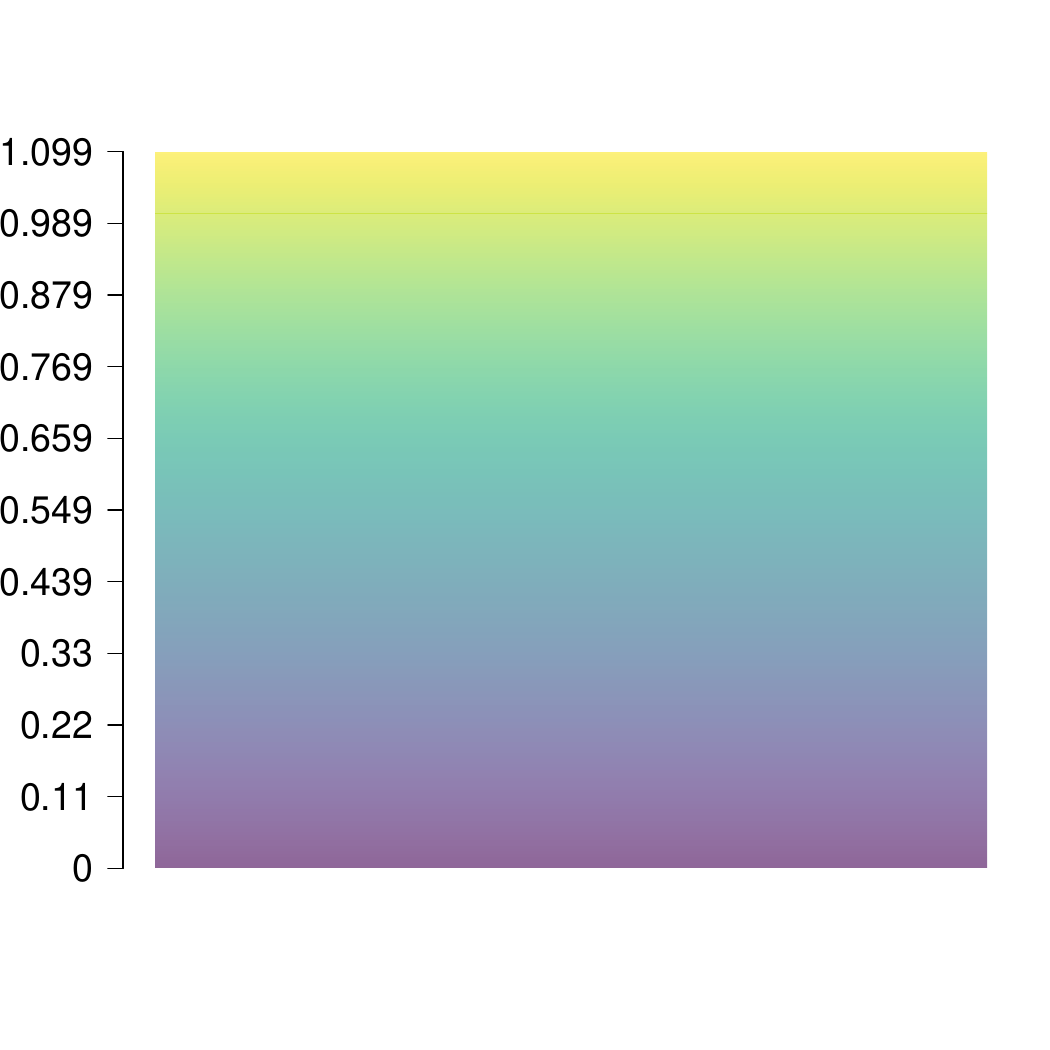}
\vspace*{-12mm}  \newline   \hspace*{-11mm}       
 \includegraphics[height=5.6cm]{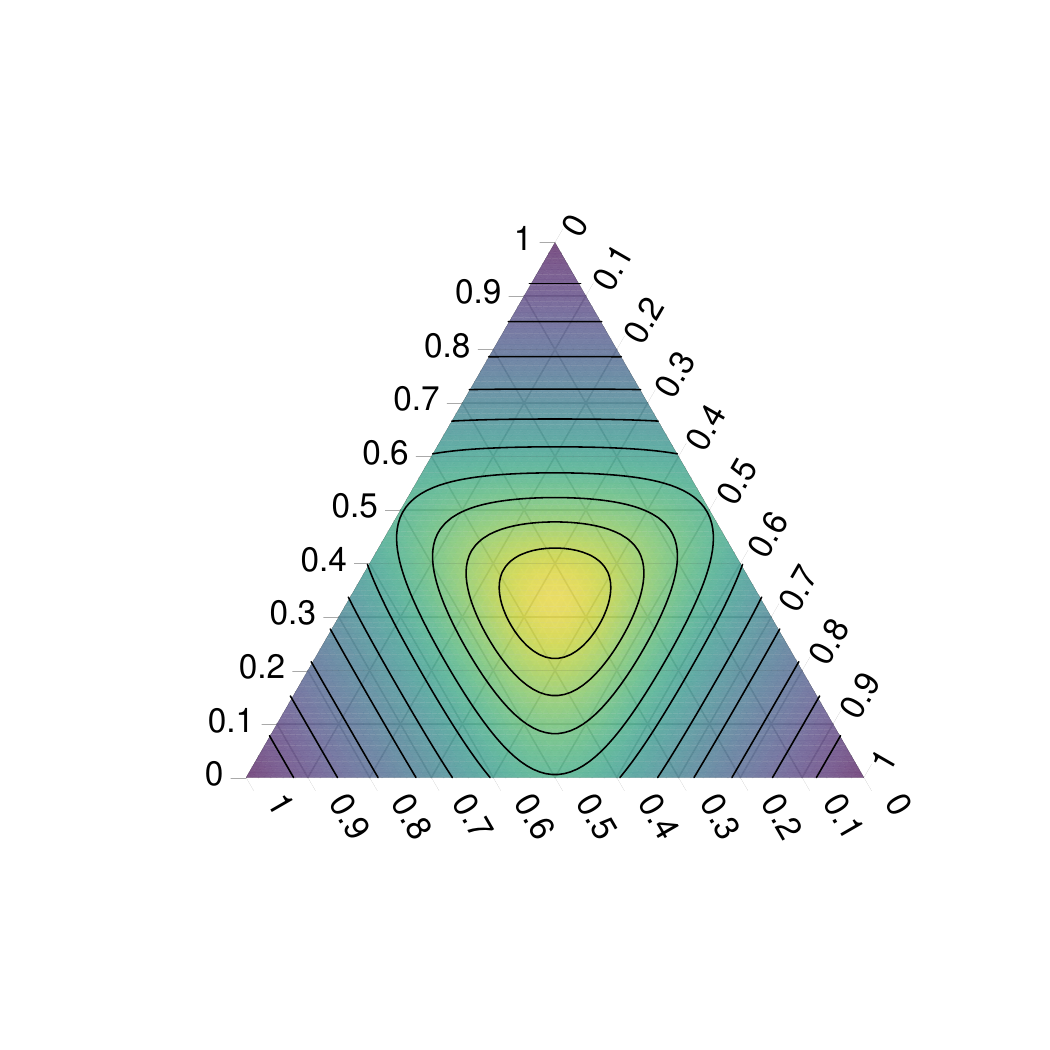}\hspace*{-5mm}
            \includegraphics[trim=0mm -19mm 133mm -5mm, clip, height=5.6cm]{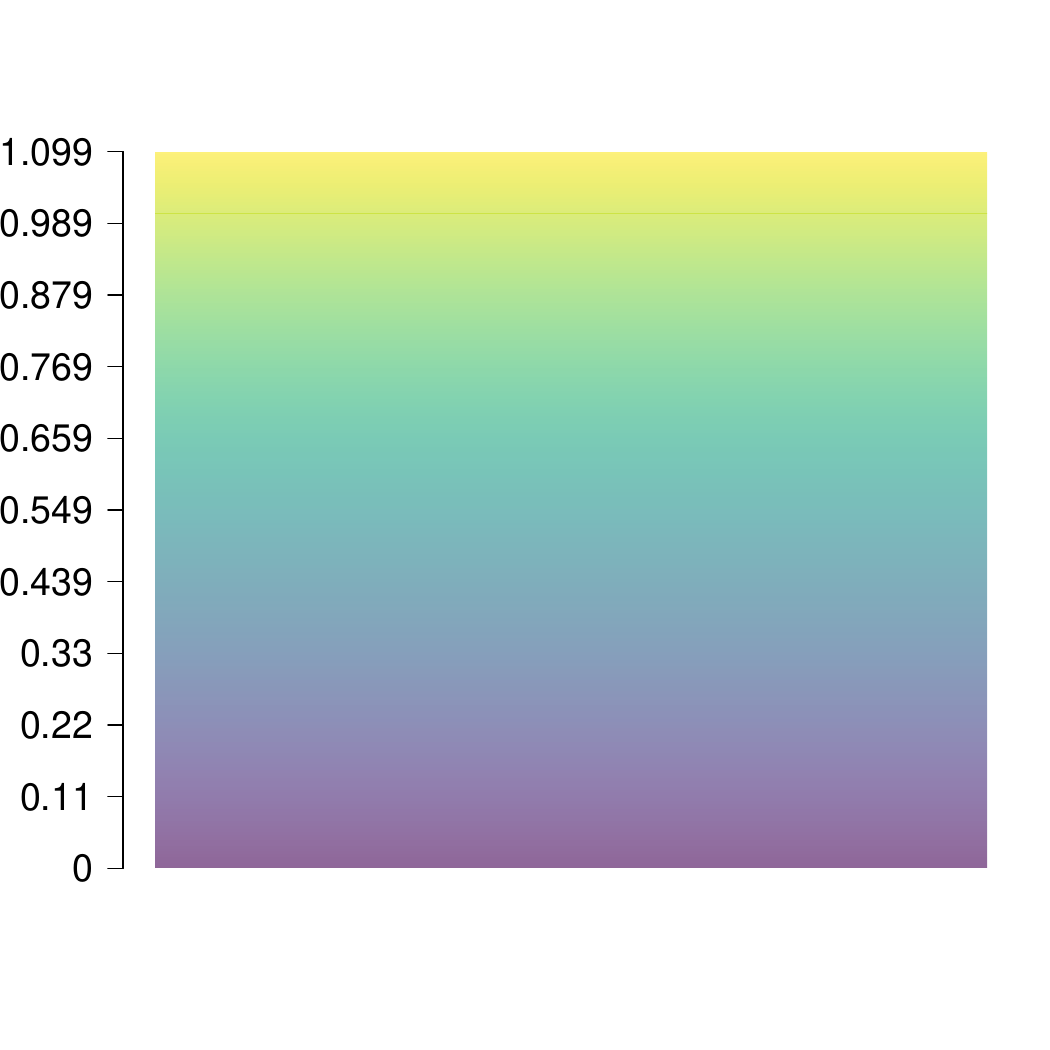}
       \hspace*{-1mm}     
 \includegraphics[height=5.6cm]{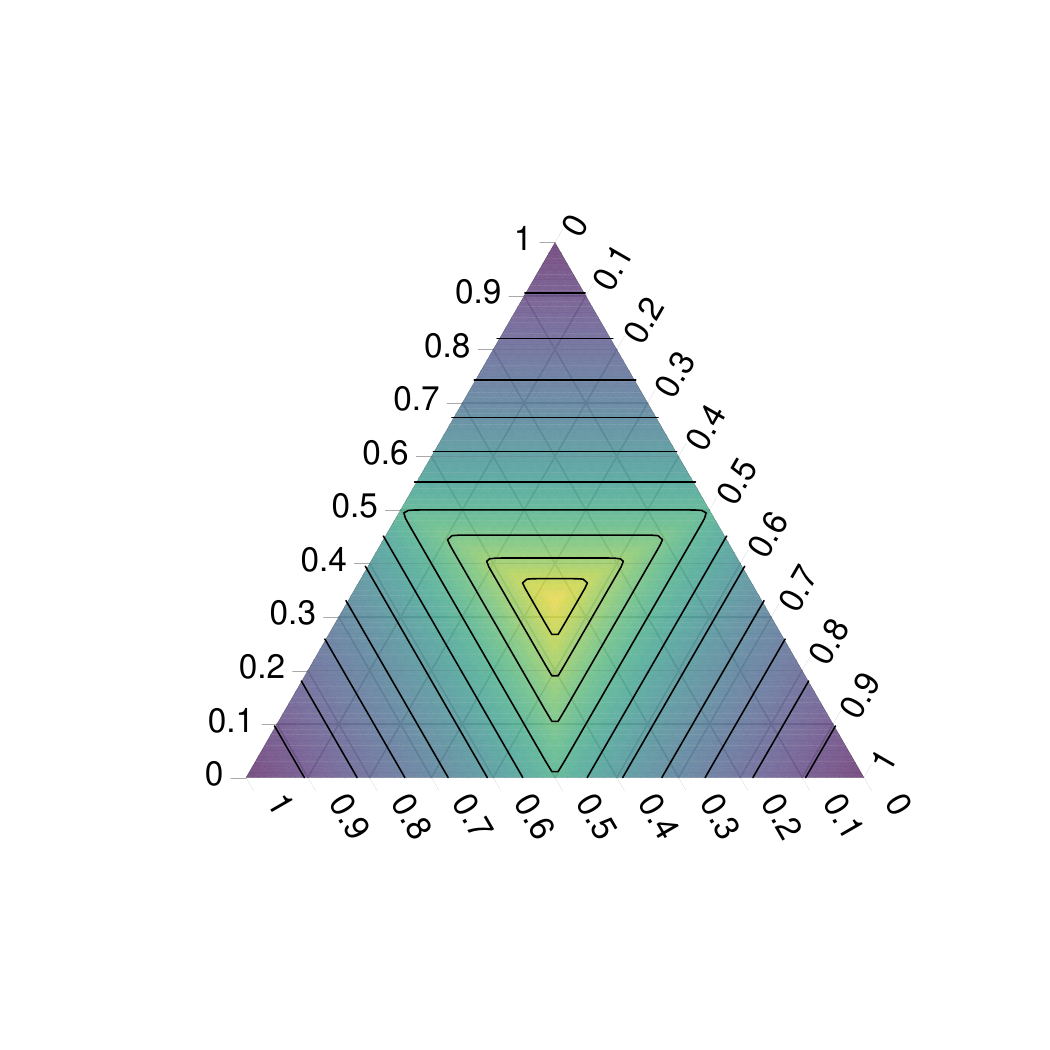}\hspace*{-5mm}
            \includegraphics[trim=0mm -19mm 133mm -5mm, clip, height=5.6cm]{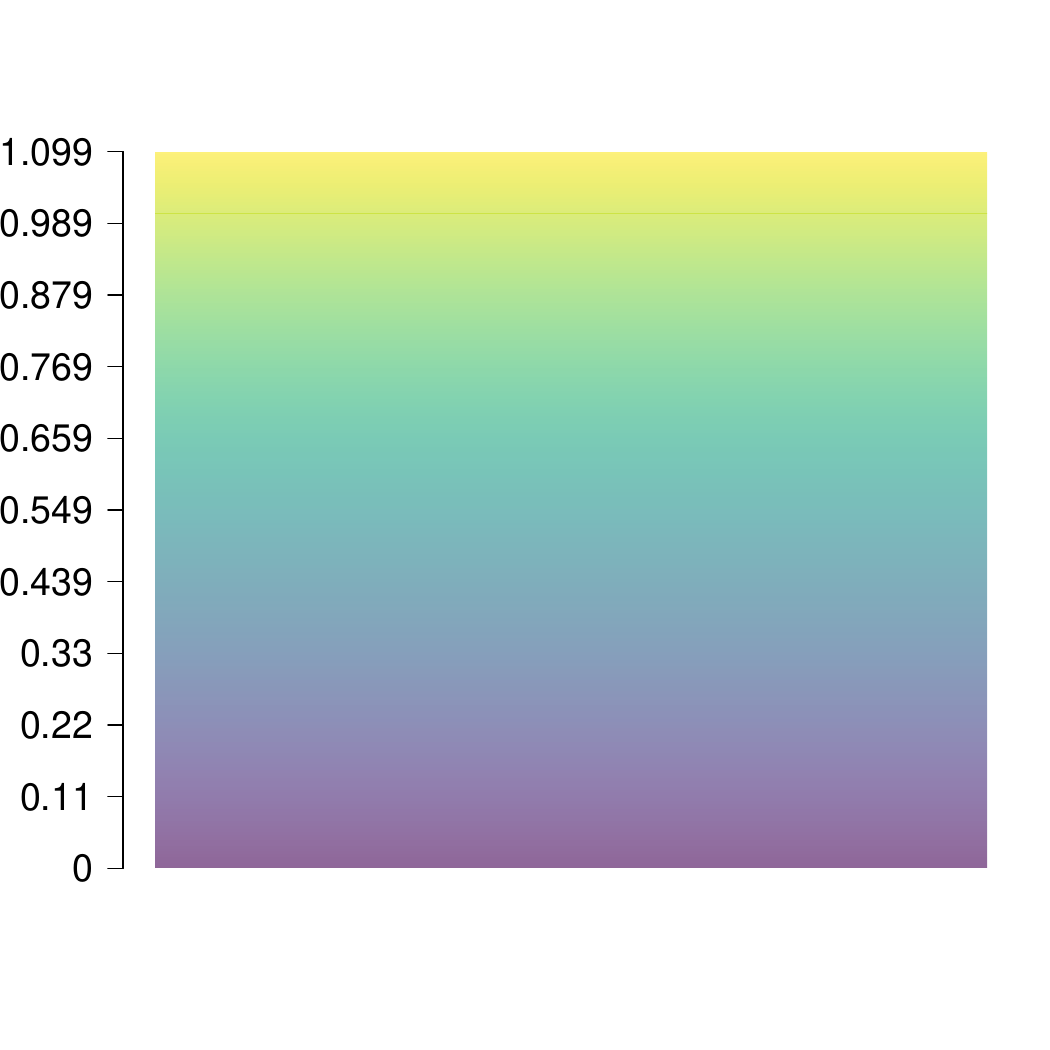}
            \caption{Rényi entropy of orders  $q=0.5$ (top left), 2 (top right), 5 (bottom left) and 100 (bottom right) for $\bar{p}= (p_1, p_2, p_3)$. In all cases, values vary from 0, the minimum uncertainty at the degenerate distributions corresponding to the vertices, to $\ln(3)$, the maximum uncertainty associated with the equiprobability central point, with dissimilar transition patterns reflecting  different distortion effects derived from the deformation parameter specifications. In particular, the monotonic non-increasing character of Rényi entropy with respect to $q$ is visualized.}
\label{figure:Renyi_entropy_0.5-2-5-100}
\end{figure}

\end{paragraph}

\begin{paragraph}{Campbell diversity index -- Campbell (1966)} 

\medskip
Given a discrete probability distribution $\bar{p} = (p_1,\dots,p_n)$, Campbell (1966) justifies measuring the intrinsic number of states of the system, according to $\bar{p}$, by the exponential of Shannon entropy and, as an extension, of Rényi entropy; i.e. the numbers  
$$ DI(\bar{p}) = e^{H(\bar{p})}  \qquad  DI_{q}(\bar{p}) = e^{H_q(\bar{p})} $$
(write $DI_1(\bar{p})$ for $DI(\bar{p})$).

In particular:
\begin{itemize}
\item For a degenerate probability distribution, say $\bar{p} = (1, 0,...,0)$, we have $DI_q = 1$.

\item In the case of equiprobability, $DI_q = n$.
\end{itemize}

Figure \ref{figure:Campbell_DI_1-2} represents the values of Campbell diversity indices $DI_1(\bar{p})$ and $DI_2(\bar{p})$ for a ternary probability distribution, $\bar{p} = (p_1,p_2,p_3)$.

\begin{figure} \hspace*{-11mm}
            \includegraphics[height=5.6cm]{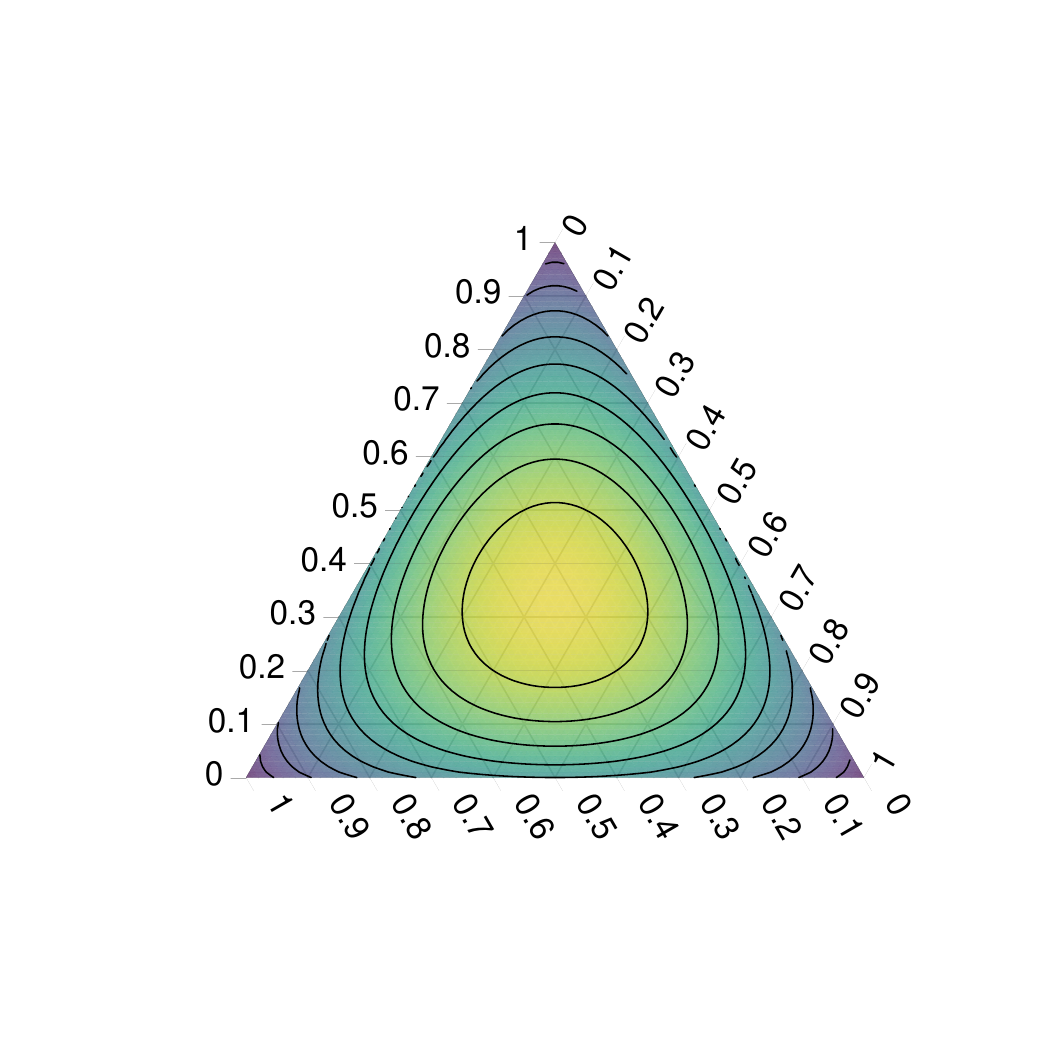}\hspace*{-5mm}
            \includegraphics[trim=0mm -19mm 133mm -5mm, clip, height=5.6cm]{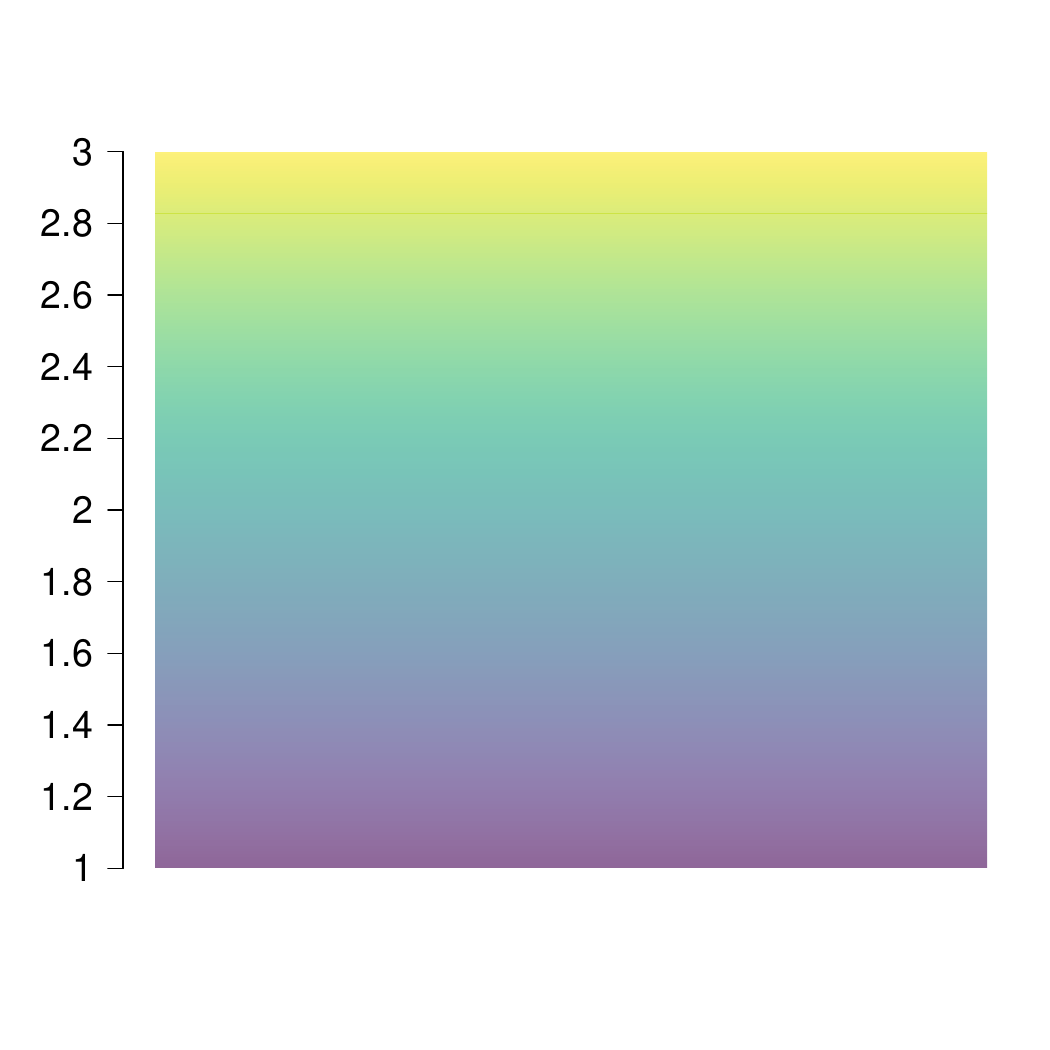}
     \hspace*{-1mm}       
\includegraphics[height=5.6cm]{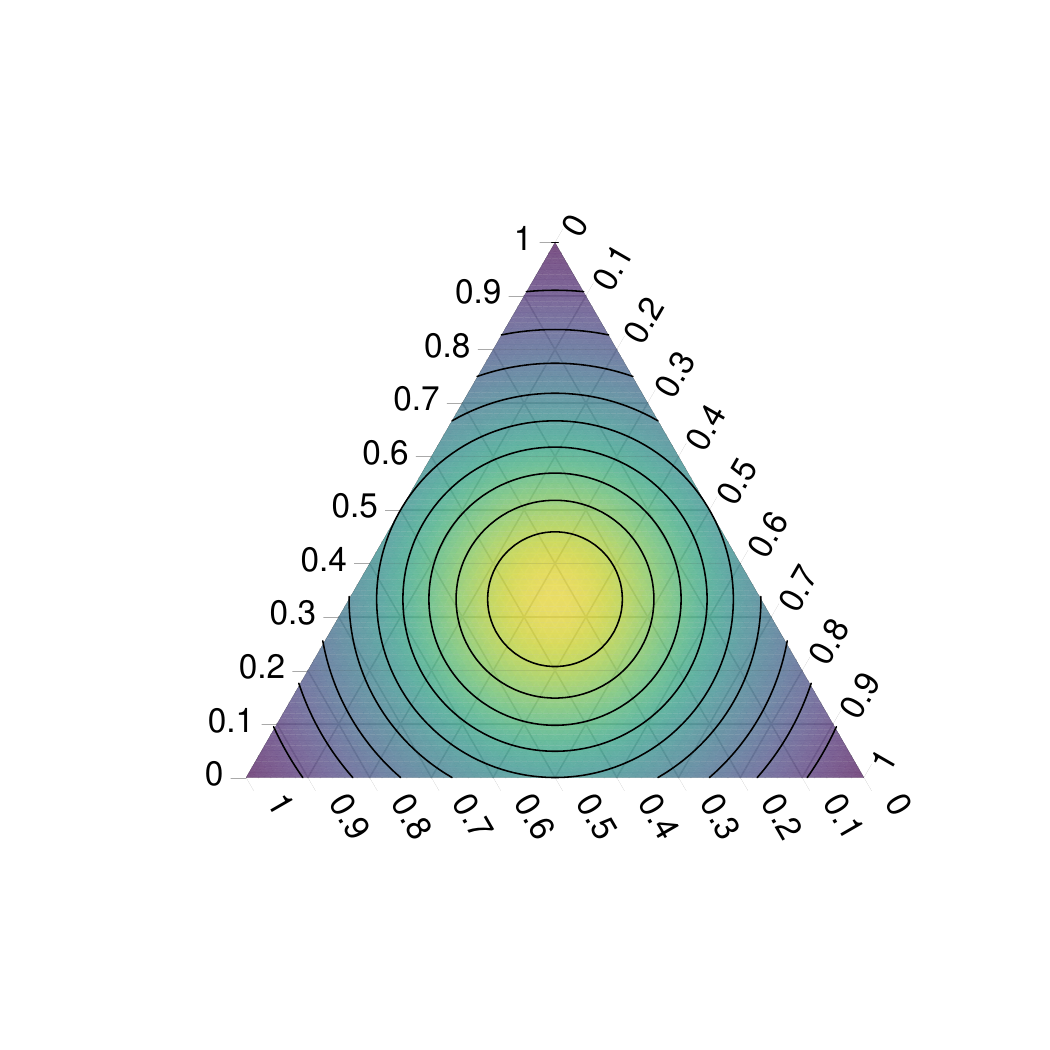}\hspace*{-5mm}
            \includegraphics[trim=0mm -19mm 133mm -5mm, clip, height=5.6cm]{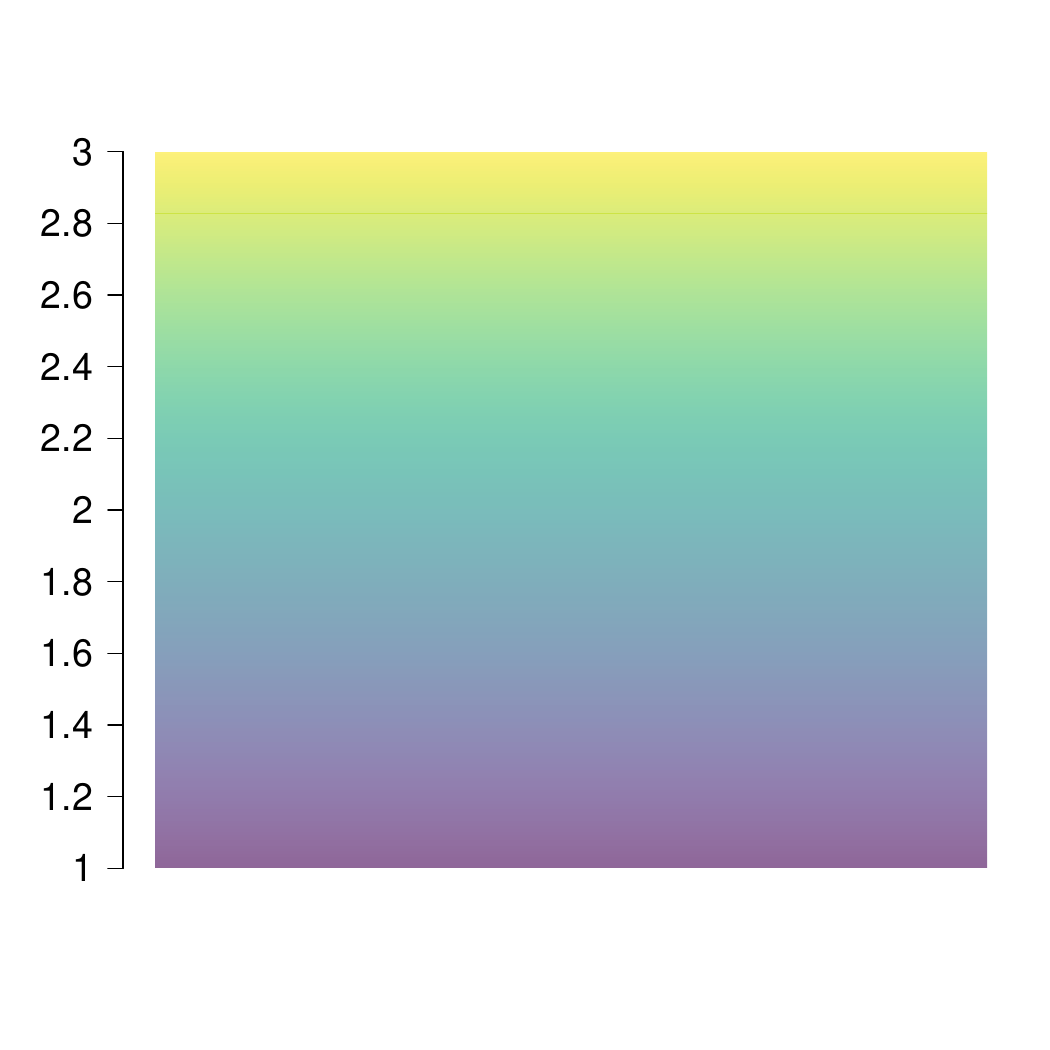}
             \caption{Campbell diversity indices $DI_1(\bar{p})$ (left) and $DI_2(\bar{p})$ (right) for $\bar{p}= (p_1, p_2, p_3)$. The values continuously vary from  1, the minimum number of intrinsic states corresponding to the degenerate distributions at the  vertices, to 3, the maximum associated with the equiprobability distribution at the central point. Transitions differ according to the specified order, i.e. the value of the deformation parameter.} 
\label{figure:Campbell_DI_1-2}
\end{figure}

As commented in section \ref{section:complexity}, Campbell diversity index provides a meaningful direct interpretation of generalized complexity measures based on Rényi entropy.
  
\end{paragraph}

\subsubsection{Continuous probability distributions}
\label{subsubsection:entropy-continuous}

The definitions given in \ref{subsubsection:entropy-discrete} of Shannon and Rényi entropies, as well as of Campbell diversity index,  are conveniently reformulated for continuous probability distributions. For a proper correspondence between the discrete and continuous versions, the latter must be specifically interpreted in terms of uncertainty at unit scale resolution, as discussed below. This constitutes the basis for Batty's formalization of a concept of \emph{spatial entropy} (Batty 1974), later reinterpreted, in a certain sense, as a form of \emph{complexity} (see discussion in section \ref{section:complexity}).

\begin{paragraph}{Shannon entropy, continuous version -- Shannon (1948)}
\medskip
For a continuous distribution with probability density function $\{f(\mathbf{x}): \mathbf{x} \in \mathds{R}^d\}$, the (Shannon) entropy is defined by

$$ H(f) := - \int_{\mathds{R}^d}\ln (f(\mathbf{x})) f(\mathbf{x}) d\mathbf{x} .$$

Well-known basic properties of continuous Shannon entropy are: 
\begin{itemize}
  \item Infimum and supremum of $H$:
  
$$ H_{inf}= -\infty  \qquad     H_{sup}= + \infty .$$

  \item Continuous Shannon entropy satisfies `extensivity'. 
  \medskip
\end{itemize}

\end{paragraph}

\begin{paragraph}{R\'{e}nyi entropy of order $q$, continuous version -- Rényi (1961)}

For a continuous distribution with probability density function $\{f(\mathbf{x}): \mathbf{x} \in \mathds{R}^d\}$, the (Rényi) entropy of order $q$ is defined by
$$ H_q(f) :=  \frac{1}{1-q} \ln \left(\int_{\mathds{R}^d} f^q(\mathbf{x}) d\mathbf{x} \right) = \frac{1}{1-q}\ln\left( E[f^{q-1}]\right) \qquad (q\neq1). $$

The following basic properties hold:
\begin{itemize}
\item Shannon entropy $H(f)$ is the limiting case of Rényi entropy $H_q(f)$ as $q\rightarrow 1$.

\item Infimum and supremum of $H_q$:
$$ {H_q}_{inf}= -\infty   \qquad  {H_q}_{sup}= +\infty . $$

\item Continuous Rényi entropy satisfies `extensivity' for independent systems.
\end{itemize}
\end{paragraph}

\begin{paragraph}{Campbell diversity index, continuous version -- Campbell (1966)}

The definition of Campbell diversity index can be formally extended to systems characterized by  continuous probability distributions. For a continuous distribution with probability density function $\{f(\mathbf{x}): \mathbf{x} \in \mathds{R}^d\}$,
$$ DI(f) = e^{H(f)}  \qquad  DI_{q}(f) = e^{H_q(f)} $$
(as before, $DI_1(f)$ can be written for $DI(f)$). In this case, $DI_{q}(f)$ varies on $(0,+\infty)$, with an appropriate interpretation as a measure of `extent' (or `spread') of the distribution represented by $f$. 
\end{paragraph}

\begin{paragraph}{Continuous vs. discrete entropy versions}

Despite their formal analogy (in fact, the same notation is commonly used, except for the type of argument under consideration, implicit or explicitly), there is not a direct correspondence between the definitions of discrete and continuous Shannon and Rényi entropy versions. However, they are in fact related in terms of a proper interpretation of uncertainty as a concept intrinsically linked to scale in the continuous probability distribution scenario. For simplicity, this is justified below in reference to the particular case of equiprobability in dimension $d=1$.

For a continuous Uniform distribution $U(I_L)$ on the interval $I_L = [0,L]$ \/ ($L \in \mathds{R}^{+}$),  we have
$$
f_{U(I_L)}(\cdot) \equiv \frac{1}{L} \qquad \longrightarrow \qquad H_q(f_{U(I_L)}) = \ln (L) \quad (\forall q). 
$$  
Dividing the interval $I_L$ into $N$ subintervals of equal length, $\Delta = \frac{L}{N}$, we obtain the corresponding discrete distribution $\bar{p}^{\Delta} = (p_1^{\Delta}, \dots, p_N^{\Delta})$ with 
$$
p_i^{\Delta} =    \frac{1}{L}\cdot \Delta = \frac{1}{N} \quad (\forall i) \qquad \longrightarrow \qquad H_q(\bar{p}^{\Delta}) = \ln (N) \quad (\forall q).
$$
In general, we then have
$$
H_q(f_{U(I_L)}) = \ln (L) = \ln (N\cdot \Delta) = H_q(\bar{p}^{\Delta}) + \ln (\Delta) .
$$
In particular, for $L \in \mathds{N}$ and $N=L$, i.e. $\Delta = 1$, the continuous and discrete entropy values coincide. This equivalence supports the interpretation that the continuous entropy versions implicitly assume that uncertainty is quantified as a degree of determination of the system state up to unit scale resolution.  In fact, for $L < 1$,  we have $H_q(f_{U(I_L)}) < 0$, consistently with 
the fact that in this case a degree of determination strictly more precise than the unit scale is a priori implied by the specified range for the distribution support. 

For  general continuous probability distributions with support in $\mathds{R}$, negative entropy values can thus be reinterpreted as `a positive degree of certainty' about the system state at each realization, with respect to the reference zero-entropy case represented by the standard $U([0,1])$ distribution. (This argument can be extended to the $d$-dimensional case under appropriate considerations.)    

\end{paragraph}

\begin{paragraph}{Batty's approach to `spatial entropy' -- Batty (1974), etc} 

Based on a Riemann-type approximation of the integral defining the continuous Shannon entropy, Batty (1974) suggested to write
$$
H(f) = \lim_{\Delta \mathbf{x}_i \rightarrow 0} \left[- \sum_i p_i \ln \left(\frac{p_i}{\Delta \mathbf{x}_i} \right)   \right],
$$
for a partition of space into cells of size $\Delta \mathbf{x}_i$, respectively, with $p_i$ being the probability mass concentrated on cell $i$. The term within brackets is interpreted as the \emph{spatial entropy} corresponding to that particular partition of space.

Batty, Morphet, Masucci and Stanilov (2014) proposed rewriting this approximation in terms of the discrete Shannon entropy, as given by the decomposition  
 
$$ 
H(\bar{p}) = - \sum_i p_i \ln \left(\frac{p_i}{\Delta \mathbf{x}_i} \right)  - \sum_i p_i \ln \left(\Delta \mathbf{x}_i \right) = S + Z,
$$
where, as mentioned, $S$ is the approximation to the continuous entropy based on the particular partition of space adopted, and $Z$ is explicitly interpreted by the authors  as ``the approximation to the information associated with the sizes of the events comprising the distribution which enable densities to be measured'', pointing out that the interesting aspect to be examined on $H(\bar{p})$ is the ``numerical co-variation of its elements $S$ and $Z$''. 

In fact, from our previous interpretation in terms of scale, the term $\ln \left(\Delta \mathbf{x}_i \right)$ can be regarded as the `spatial information content' (in analogy to Hartley's approach) of cell $i$ intrinsic to its size, a quantity which becomes negative for $\Delta \mathbf{x}_i < 1$; hence, in this sense, the value
$$
-Z = \sum_i p_i \ln \left(\Delta \mathbf{x}_i \right) 
$$
represents the $\bar{p}$-mean spatial information content for this particular partition, noting that it tends to $-\infty$ as the partition gets indefinitely finer.    

\end{paragraph}

\subsection{Divergence measures (Kullback-Leibler, Rényi)}  
\label{subsection:divergence}

Measures of entropy can be used for a first comparison, in `global' terms, of the uncertainty associated with different probability distributions on a given system. However, since very different probability distributions may lead to the same entropy values, a `local' comparison based on the specific assignment of probabilities for each possible state provides the relevant information for structural assessment. This is the general purpose of divergence measures. Here, the definitions of the Kullback-Leibler and Rényi divergence measures are introduced. A notion of a `relative diversity index', in analogy to Campbell's diversity index approach, is stated.

\begin{paragraph}{Kullback-Leibler (1951) and R\'enyi (1961) divergence measures}

For two discrete probability distributions $\bar{p}_1=(p_{11},p_{12},...,p_{1n})$ and $\bar{p}_2=(p_{21},p_{22},...,p_{2n})$ (with $\bar{p}_1$ being absolutely continuous with respect to $\bar{p}_2$), on a given set of $n$ states,  Kullback-Leibler (directed) divergence and Rényi (directed) divergence of order $q$ of $\bar{p}_1$ from $\bar{p}_2$ are respectively defined as  
\begin{align} 
KL(\bar{p}_1\|\bar{p}_2) & := \sum_{i} p_{1i} \ln \left(\frac{p_{1i}}{p_{2i}}\right)  = E_{\bar{p}_1}\left[\ln\left( \frac{\bar{p}_1}{\bar{p}_2}\right) \right],\nonumber \\
 & \nonumber \\
H_q(\bar{p}_1\|\bar{p}_2) & :=  \frac{1}{q-1} \ln \left(\sum_{i} p_{1i} \left( \frac{p_{1i}}{p_{2i}} \right)^{q-1} \right) \nonumber \\
   & = \frac{1}{q-1} \ln \left( E_{\bar{p}_1}\left[\left( \frac{\bar{p}_{1}}{\bar{p}_{2}} \right)^{q-1}\right]\right) \qquad (q\neq1).\nonumber
\end{align}
For $q \rightarrow 1$, $H_q(\bar{p}_1\|\bar{p}_2)$ tends to $KL(\bar{p}_1\|\bar{p}_2)$ (hence with the latter being also denoted as $H_1(\bar{p}_1\|\bar{p}_2)$).

Similarly, for two continuous distributions with respective probability density functions $\left\{f(\mathbf{x}): \mathbf{x} \in \mathds{R}^d \right\}$ and $\left\{g(\mathbf{x}): \mathbf{x} \in \mathds{R}^d \right\}$ (with $f$ being absolutely continuous with respect to $g$), the (directed) divergences of $f$ from $g$ are correspondingly defined as  
\begin{align} 
KL(f\|g) & := \int_{\mathds{R}^d} f(\mathbf{x}) \ln \left(\frac{f(\mathbf{x})}{g(\mathbf{x})} \right) d\mathbf{x} = E_{f}\left[\ln\left( \frac{f}{g}\right) \right],\nonumber \\
 & \nonumber \\
H_q(f\|g) & :=  \frac{1}{q-1} \ln \left(\int_{\mathds{R}^d} f(\mathbf{x}) \left( \frac{f(\mathbf{x})}{g(\mathbf{x})} \right)^{q-1} \right) \nonumber \\
   & = \frac{1}{q-1} \ln \left( E_{f}\left[\left( \frac{f}{g} \right)^{q-1}\right]\right) \qquad (q\neq1).\nonumber
\end{align}
As before, for $q \rightarrow 1$, $H_q(f\|g)$ tends to $KL(f\|g)$ (also denoted as $H_1(f\|g)$).

\end{paragraph}

\begin{paragraph}{A `pivotal' quantity: The `information difference'}

There is one particular instance for which Kullback-Leibler and Rényi divergences are equal to the  differences of the corresponding Shannon and Rényi entropies, respectively. This is related to the concept of `information difference', which is reviewed later in section \ref{section:complexity}       
in the context of complexity.

For a discrete probability distribution $\bar{p}=(p_1,p_2,...,p_n)$, the \emph{information difference} is defined as the entropy defect with respect to the case of equiprobability,
$$
ID(\bar{p}) := H_{max} - H(\bar{p}) = \ln(n) - H(\bar{p}) = KL\left(\bar{p} \| \left[\frac{1}{n}\right]\right).
$$
More generally, the \emph{information difference of order $q$} is defined as 
$$
ID_q(\bar{p}) := {H_q}_{max} - H_q(\bar{p}) = \ln(n) - H_q(\bar{p}) = H_q\left(\bar{p} \| \left[\frac{1}{n}\right]\right),
$$
which includes the previous case for $q=1$ under the aforementioned limiting equivalence. The divergence expression reflects the fact that the `locality' aspect, in comparing probabilities for each specific state, becomes irrelevant in this particular case due to uniformity in the reference distribution.

In the continuous case, since in general Shannon and Rényi entropies are unbounded, such a concept of information difference can be defined only under restrictions for which a certain distribution with maximum entropy may exist. In particular, and for simplicity, let us consider the subclass, denoted here as $\mathcal{P}_I$, of continuous probability distributions with support contained in a given compact subset $I \subset \mathds{R}$ with positive Lebesgue measure, $\lambda(I) >0$. The distribution within  $\mathcal{P}_I$ with the maximum (Shannon, Rényi) entropy , $H_{q_{max}}^{\mathcal{P}(I)}$, is the continuous Uniform distribution on $I$, denoted $U(I)$, with probability density function 
$$
f_{U(I)} \equiv \frac{1}{\lambda(I)} 
$$
(except possibly for a subset $I_0 \subset I$ with $\lambda(I_0) = 0$), noting that
$$
H_{q_{max}}^{\mathcal{P}(I)} = H_q(f_{U(I)}) = \ln(\lambda(I)) , 
$$
for all $q$. For any probability distribution in $\mathcal{P}_I$, say with probability density function $f$, the corresponding information difference and information difference of order $q$, both relative to the subclass $\mathcal{P}_I$, can be defined, respectively, as 
$$
ID(f) := H_{max}^{\mathcal{P}(I)} - H(f) = \ln(\lambda(I)) - H(f) = KL\left(f \| f_{U(I)}\right),
$$
and
$$
ID_q(f) := H_{q_{max}}^{\mathcal{P}(I)} - H_q(f) = \ln(\lambda(I)) - H_q(f) = H_q\left(f \| f_{U(I)}\right).
$$

\end{paragraph}

\begin{paragraph}{A `relative diversity index'}

Following Campbell's `diversity index' approach, a related notion of a `relative diversity index' is formulated here is terms of divergence.

For two discrete probability distributions, denoted $\bar{p}_1=(p_{11},p_{12},...,p_{1n})$ and $\bar{p}_2=(p_{21},p_{22},...,p_{2n})$, on a given set of $n$ states, the `relative diversity index' and `relative diversity index of order $q$' of $\bar{p}_1$ with respect to $\bar{p}_2$ are respectively defined as  
$$
DI(\bar{p}_1\|\bar{p}_2) = e^{KL(\bar{p}_1\|\bar{p}_2)},  \qquad  DI_q(\bar{p}_1\|\bar{p}_2) = e^{H_q(\bar{p}_1\|\bar{p}_2)} \qquad (q\neq 1)
$$
(as before, $DI_1(\bar{p}_1\|\bar{p}_2)$ can be written for $DI(\bar{p}_1\|\bar{p}_2)$). 

In particular, for $\bar{p}_2 \equiv \left[\frac{1}{n}\right]$, 
$$
DI_q(\bar{p}_1\|\left[\frac{1}{n}\right]) = \frac{n}{DI_q(\bar{p}_1)}.
$$
Also, for $\bar{p}_1 \equiv \bar{p}_2$, we  have that $DI_q(\bar{p}_1\|\bar{p}_2) = 1$. In general, the relative diversity index (of order $q$) can be interpreted as a measure of the structural departure of $\bar{p}_1$ from $\bar{p}_2$ in terms on the state-by-state probability contribution to diversity. 

The concept of relative diversity index (of order $q$) can be similarly formulated and interpreted for continuous probability distributions, by exponentiation, under the corresponding definitions of Shannon and Renyi divergences given above (details are omitted). 
 
\end{paragraph}

\section{Complexity and Relative Complexity Measures}
\label{section:complexity}

In this section, we address some formulations of complexity and relative complexity measures, based on the concepts of entropy and divergence introduced in section \ref{section:information-entropy-divergence}, which have received special attention in the literature. As a preliminary form, proposed in the context of Geography, we refer to the notion of complexity by Batty, Morphet, Masucci and Stanilov (2014), specifically understood as departure from equilibrium. A more general statement, emerged from Physics originally under the notion of complexity as departure from both equilibrium and degeneracy, was developed as a  product-type two-parameter measure defined by the exponential of the difference of Rényi entropies, for complexity, and of Rényi divergences, for relative complexity, of different orders. The former mentioned preliminary form  can be interpreted  as a very special case of the generalized complexity measure for certain fixed parameter values. Both generalized complexity and relative complexity measures are properly interpreted from their straightforward expression in terms of the diversity and relative diversity indices introduced in section \ref{section:information-entropy-divergence}. The usefulness of maps based on relative increments and, in particular, of the curves of derivatives of Rényi entropy and divergence with respect to the deformation parameter is justified for complexity and relative complexity assessment (see related elements in section \ref{section:multifractality} in the context of multifractal complexity and relative complexity assessment).

\subsection{Entropy and Complexity}
\label{subsection:entropy-complexity}

\begin{paragraph}{Complexity uniquely as departure from equilibrium}

Shannon (1948) defined the \emph{relative entropy} of a discrete probability distribution $\bar{p} = (p_1,\dots,p_n)$ as the [entropy] ratio
$$
\frac{H(\bar{p})}{H_{max}},
$$
that is, the normalization of entropy with respect to the maximum entropy corresponding to equiprobability. He then called  [relative] \emph{redundancy} its complement,
$$
R(\bar{p}) = 1 - \frac{H(\bar{p})}{H_{max}}. 
$$  
Rewritten in terms of the \emph{information difference}   $ID(\bar{p}) = H_{max} - H(\bar{p})$, as  
$$ 
R(\bar{p}) = \frac{H_{max} - H(\bar{p})}{H_{max}} = \frac{ID(\bar{p})}{H_{max}} = \frac{KL(\bar{p}\|\left[\frac{1}{n}\right])}{H_{max}},
$$
the redundancy of $\bar{p}$ can be interpreted as a normalized measure of departure from equilibrium, hence representing, in a certain sense, the relative loss of freedom in the occurrence of the system states implied by its probabilistic structure.

Batty (1974) incorporated these concepts in his work, adopting from Theil (1967) (see also Theil 1972) the name `information gain' for $ID(\bar{p})$ (in this sense, `[relative] redundancy' can be identified as `relative information gain'). Further, Batty, Morphet, Masucci and Stanilov (2014) explicitly call $ID(\bar{p})$ the `complexity difference', and  $R(\bar{p})$ the `complexity ratio', of $\bar{p}$. Thus, the concept of `complexity' is interpreted by the authors in the specific sense of divergence from equiprobability.
  
The left plot of Figure \ref{figure:complexity-difference_vs_LMC-complexity} shows the values of information difference (`complexity difference'), $ID(\bar{p})$, for a ternary probability distribution, $\bar{p} = (p_1,p_2,p_3)$.

\begin{figure} \hspace*{-11mm}
            \includegraphics[height=5.6cm]{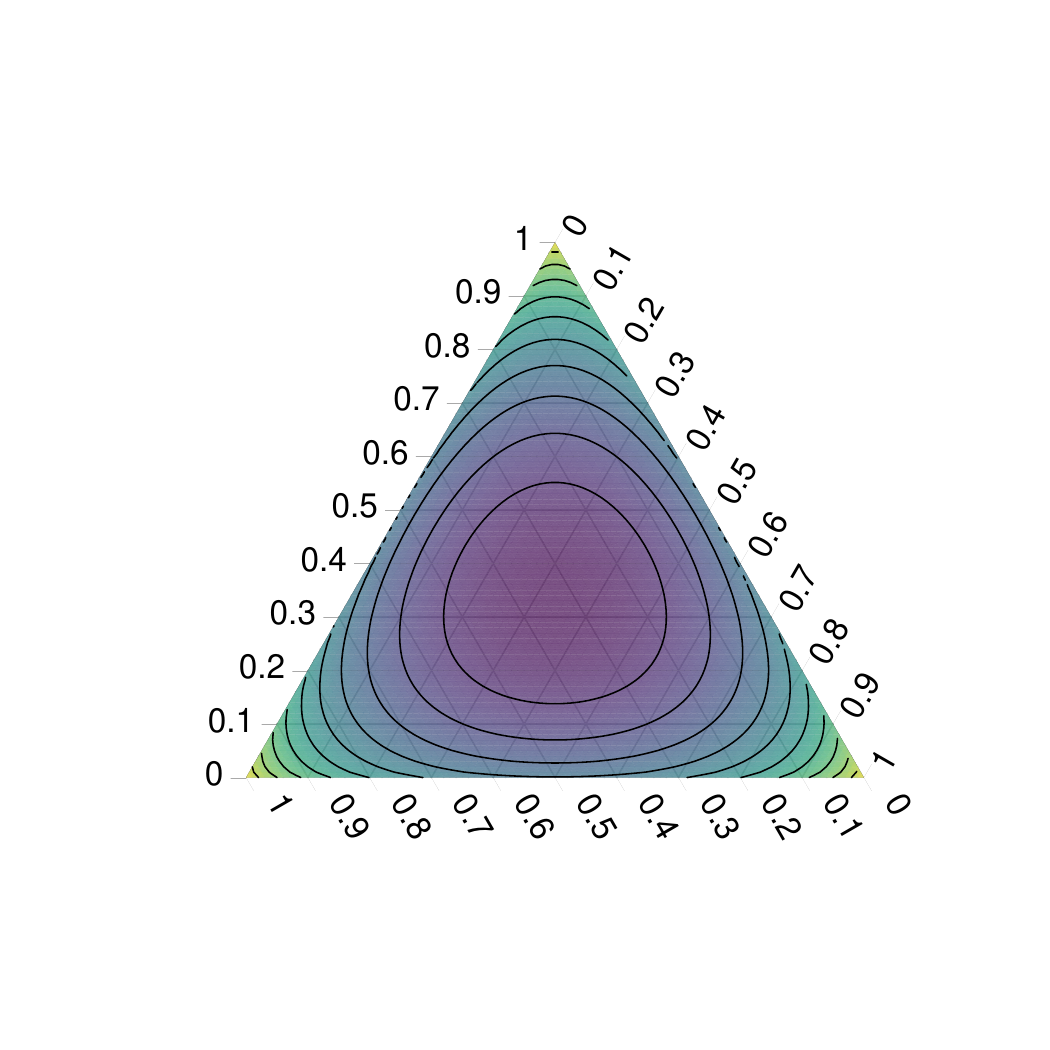}\hspace*{-5mm}
            \includegraphics[trim=0mm -19mm 133mm -5mm, clip, height=5.6cm]{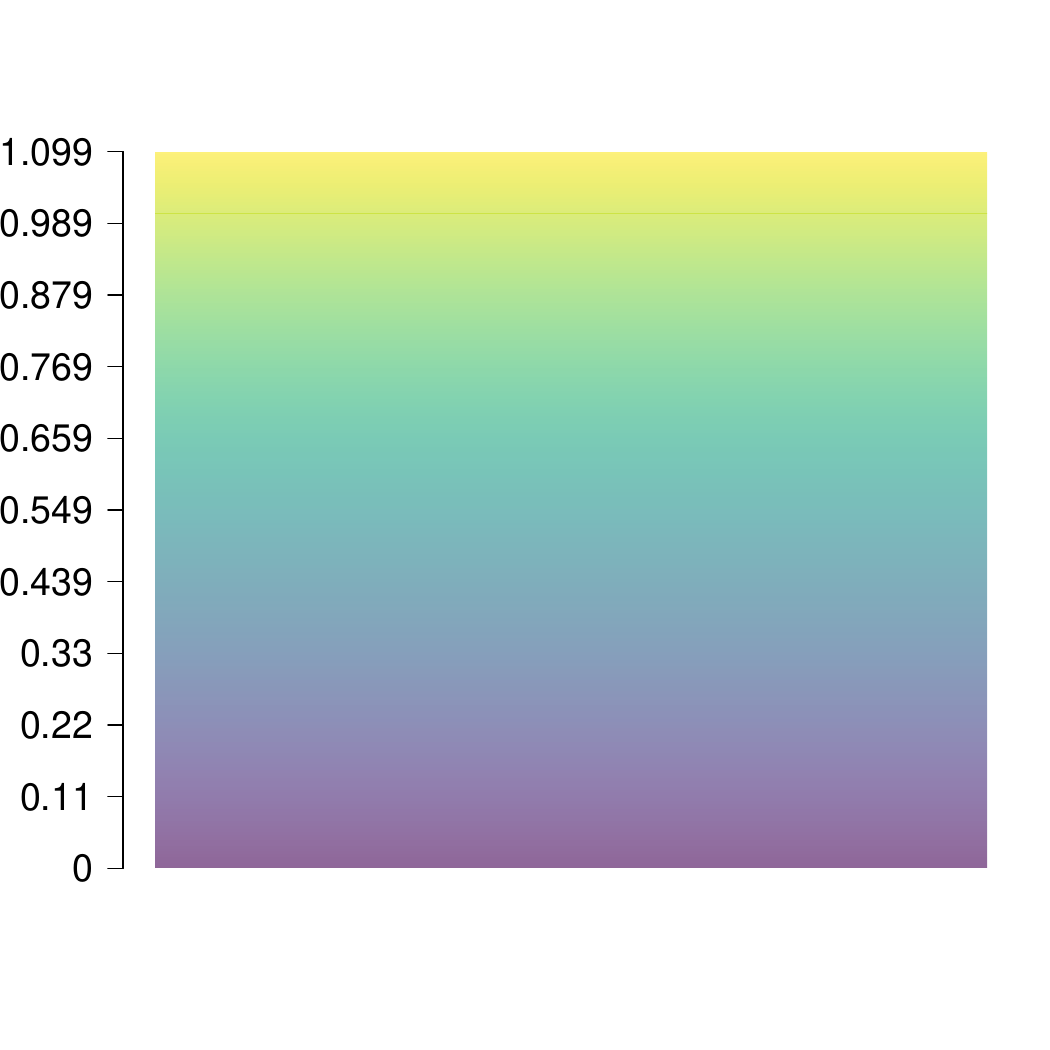}
            \hspace*{-1mm}
            \includegraphics[height=5.6cm]{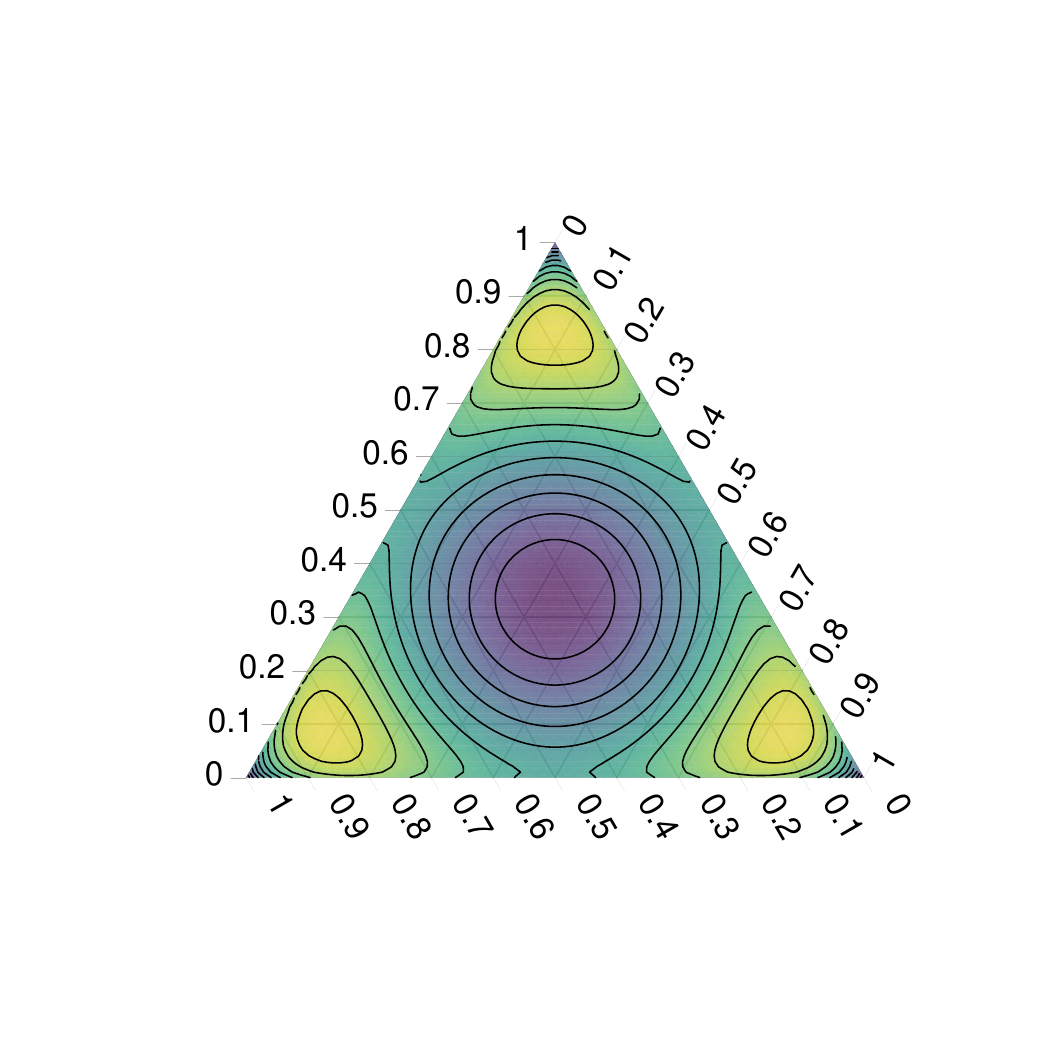}\hspace*{-5mm}
            \includegraphics[trim=0mm -19mm 133mm -5mm, clip, height=5.6cm]{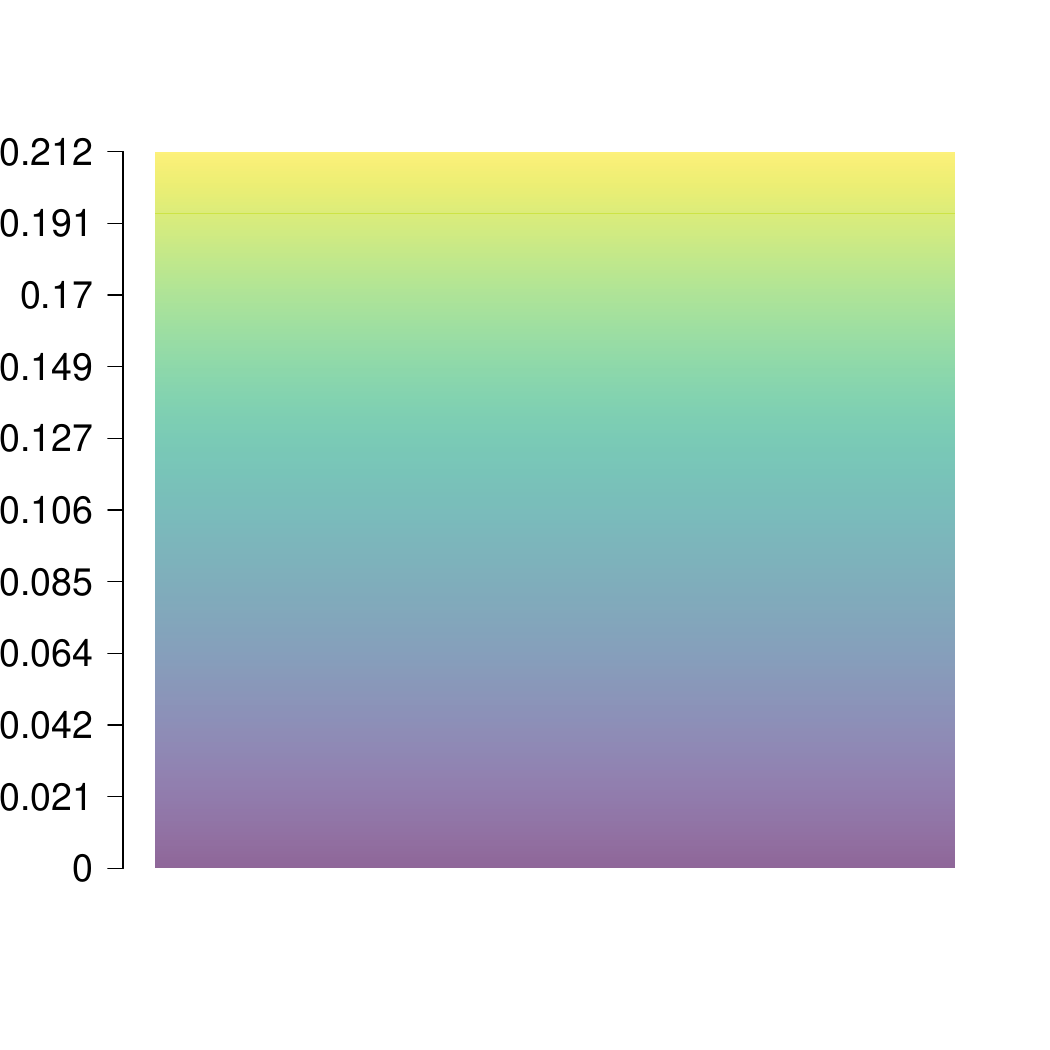}
\caption{Information difference (or `complexity difference') $ID(\bar{p})$ (left), and LMC complexity measure  $C_{LMC}(\bar{p})$ (right), for $\bar{p}= (p_1, p_2, p_3)$. In these two plots, the main conceptual difference for `complexity', understood as a measure of uniquely departure from equilibrium in the first case, and as a balance measure of departure from both equilibrium and degeneracy in the second case, is clearly visualized comparing the respective structural patterns. As a result, in particular, these two measures have an opposite behaviour for distributions in proximity to the vertices.} 
\label{figure:complexity-difference_vs_LMC-complexity}                        
\end{figure}

It is clear that similar concepts, such as `relative entropy of order $q$' and `[relative] redundancy of order $q$', can be defined in terms of Rényi entropy. Furthermore, as commented before, related versions for continuous probability distributions are meaningful under restriction to subclasses for which a maximum (Shannon, Rényi) entropy distribution may exist.    

\end{paragraph}

\begin{paragraph}{`Product' complexity measures}

In Physics, a different concept of complexity arose (e.g., Huberman and Hogg 1986) under the idea that perfect order and complete disorder both reflect, although in an opposite sense, lack of structural richness in a system, and hence minimal complexity scenarios.

In particular, this notion was adopted by L\'opez-Ruiz, Mancini and Calbet (1995) in their formulation of a product-type complexity measure defined as 
$$
 C_{LMC}(\bar{p}) =H(\bar{p}) \cdot D(\bar{p}) = \left(-\sum_{i=1}^{n}p_i\ln(p_i) \right)\left(\sum_{i=1}^{n}\left(p_i-\frac{1}{n}\right)^2\right),
$$
where the factors $H(\bar{p})$  and $D(\bar{p})$ quantify, respectively, the  `information' and   `disequilibrium' of the system characterized by the probability distribution $\bar{p}$. Thus, the balance of these two `inner' and `outer' structural aspects of the distribution is evaluated for complexity assessment.  

The right plot of Figure \ref{figure:complexity-difference_vs_LMC-complexity} represents the LMC complexity values for a ternary system.

For a proper adaptation to the continuous case, Catal\'an, Garay and L\'opez-Ruiz (2002) proposed the modified exponential version defined, for a given probability density function $f$ on $\mathds{R}$, as 
\begin{align}
 C_{LMC}^{exp}(f) &:= e^{H(f)}\cdot D(f)=\left(\exp\left\{-\int_{\mathds{R}}f(x) \ln(f(x)) dx\right\} \right)  \left( -\int_{\mathds{R}} f^2(x) dx \right)  \nonumber \\
 & = e^{H_1(f)-H_2(f)} \nonumber
\end{align}
(formal extension to distributions with support in $\mathds{R}^d$ is immediate). 

But in fact, following Campbell's (1966) approach to measuring diversity, this formulation can also be adopted as properly meaningful for the discrete case, 
$ C_{LMC}^{exp}(\bar{p}) := e^{H_1(\bar{p}) - H_2(\bar{p})}.$
Hence, for both cases, the $C_{LMC}^{exp}$ complexity measure can be interpreted as the ratio of Campbell  diversity indices of orders 1 and 2,
$$
C_{LMC}^{exp}(\bar{p}) = \frac{DI_1(\bar{p})}{DI_2(\bar{p})} \qquad C_{LMC}^{exp}(f) = \frac{DI_1(f)}{DI_2(f)}. 
$$

A further step, for a `natural' generalization of R\'enyi-entropy-based product complexity measures,  was given, in the continuous case, by L\'opez-Ruiz, Nagy, Romera and Sa\~nudo (2009) under the two-parameter formulation 
$$
C_{\alpha, \beta}(f) := e^{H_{\alpha}(f)-H_{\beta}(f)}, 
$$
for $0 < \alpha, \beta < \infty$. (As is assumed in section \ref{subsection:entropy-complexity-multifractality}, this definition can be formally extended also covering the negative range for parameters $\alpha$ and $\beta$, with proper specifications regarding the properties satisfied by the measure.)  

As before, this measure can be meaningfully adopted for the discrete case,
$$ 
C_{\alpha, \beta}(\bar{p}) := e^{H_{\alpha}(\bar{p})-H_{\beta}(\bar{p})} .
$$
Again, in terms of Campbell diversity indices of orders $\alpha$ and $\beta$,
$$
C_{\alpha, \beta}(\bar{p}) = \frac{DI_{\alpha}(\bar{p})}{DI_{\beta}(\bar{p})}, \qquad C_{\alpha, \beta}(f) = \frac{DI_{\alpha}(f)}{DI_{\beta}(f)} .
$$
In particular, we then have 
$$
C_{LMC}^{exp}(\cdot) \equiv C_{1,2}(\cdot).  
$$
Furthermore, noting that, in the discrete case, 
$$
H_{\alpha}\left(\left[\frac{1}{n}\right]\right) = \ln (n) = H_{0}(\bar{p}), 
$$
for all $\alpha \geq 0$, and for all $\bar{p}=(p_1,\dots,p_n)$ with $p_i > 0$, $i=1,\dots,n$, the `complexity difference' $ID(\bar{p})$ and the `complexity ratio'  $R(\bar{p})$ (as named by Batty, Morphet, Masucci and Stanilov 2014) can also be seen as special limit cases of the generalized two-parameter complexity measure for $\alpha \rightarrow 0$ and $\beta = 1$, according to the relations
$$
ID(\bar{p}) =  \ln \left( C_{0,1}(\bar{p}) \right), \qquad R(\bar{p}) = \frac{\ln \left( C_{0,1}(\bar{p}) \right) }{H_{max}}  = \frac{\ln \left( C_{0,1}(\bar{p}) \right)}{\ln (n)} = \ln \left[ \left( C_{0,1}(\bar{p}) \right)^{\frac{1}{\ln(n)}} \right].
$$

Plots in Figure \ref{figure:generalized-complexity(1,2)-(2,10)-(0.5,1)-(0.5,10)} show the values of $C_{\alpha, \beta}(\bar{p})$ for a ternary distribution, $\bar{p}= (p_1, p_2, p_3)$, based on the specifications  $(\alpha, \beta) = (1,2)$, $(2,10)$, $(0.5,1)$ and $(0.5,10)$, respectively.

\begin{figure} \hspace*{-11mm}
            \includegraphics[height=5.6cm]{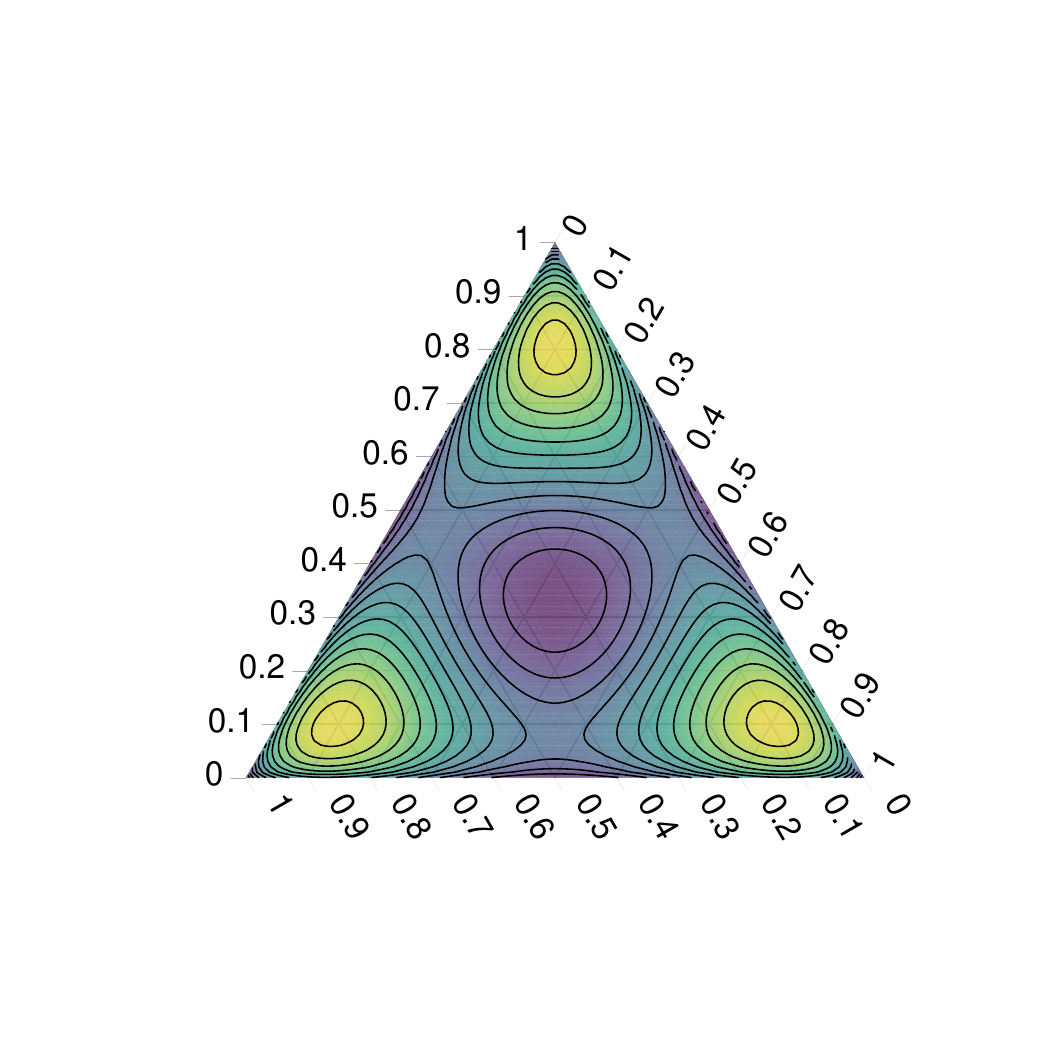}\hspace*{-5mm}
            \includegraphics[trim=0mm -19mm 133mm -5mm, clip, height=5.6cm]{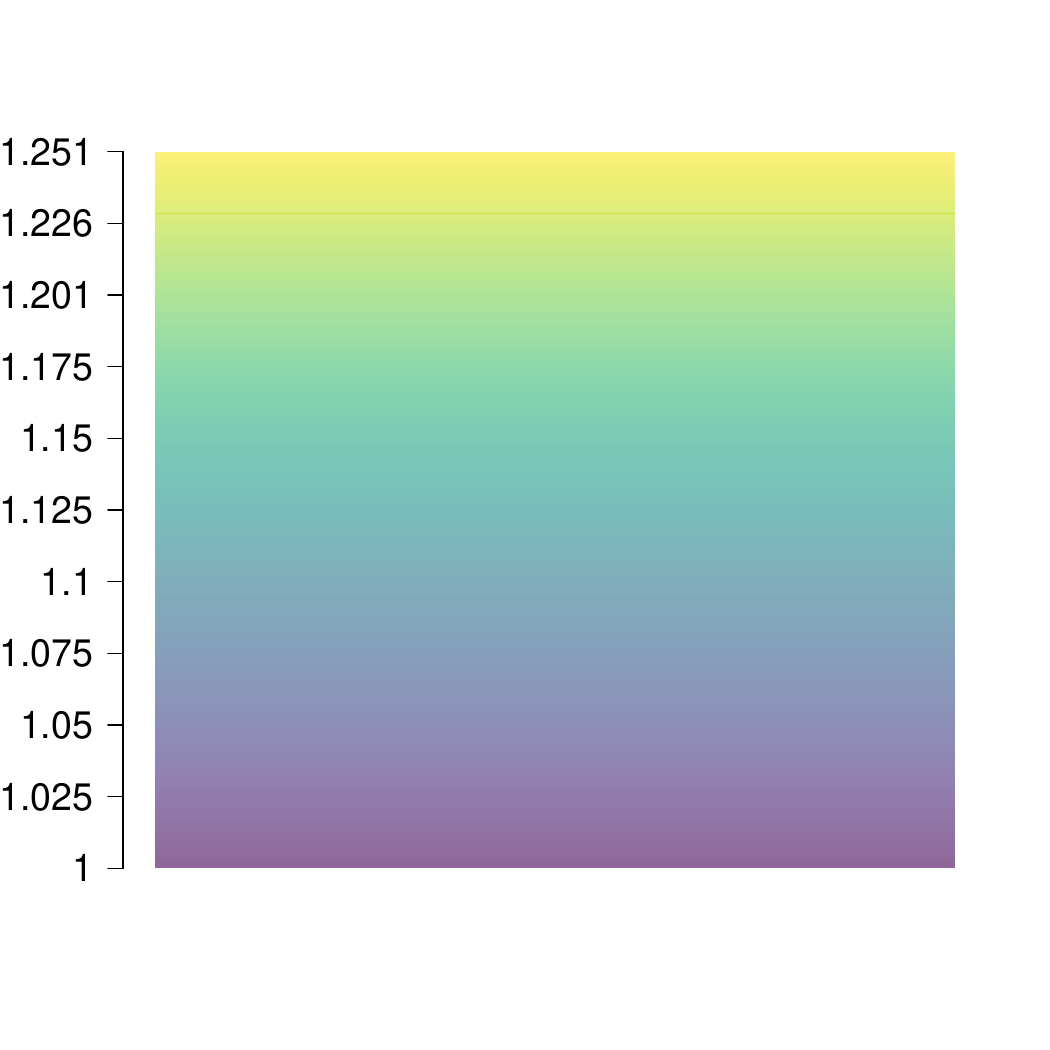}
            \hspace*{-1mm}
            \includegraphics[height=5.6cm]{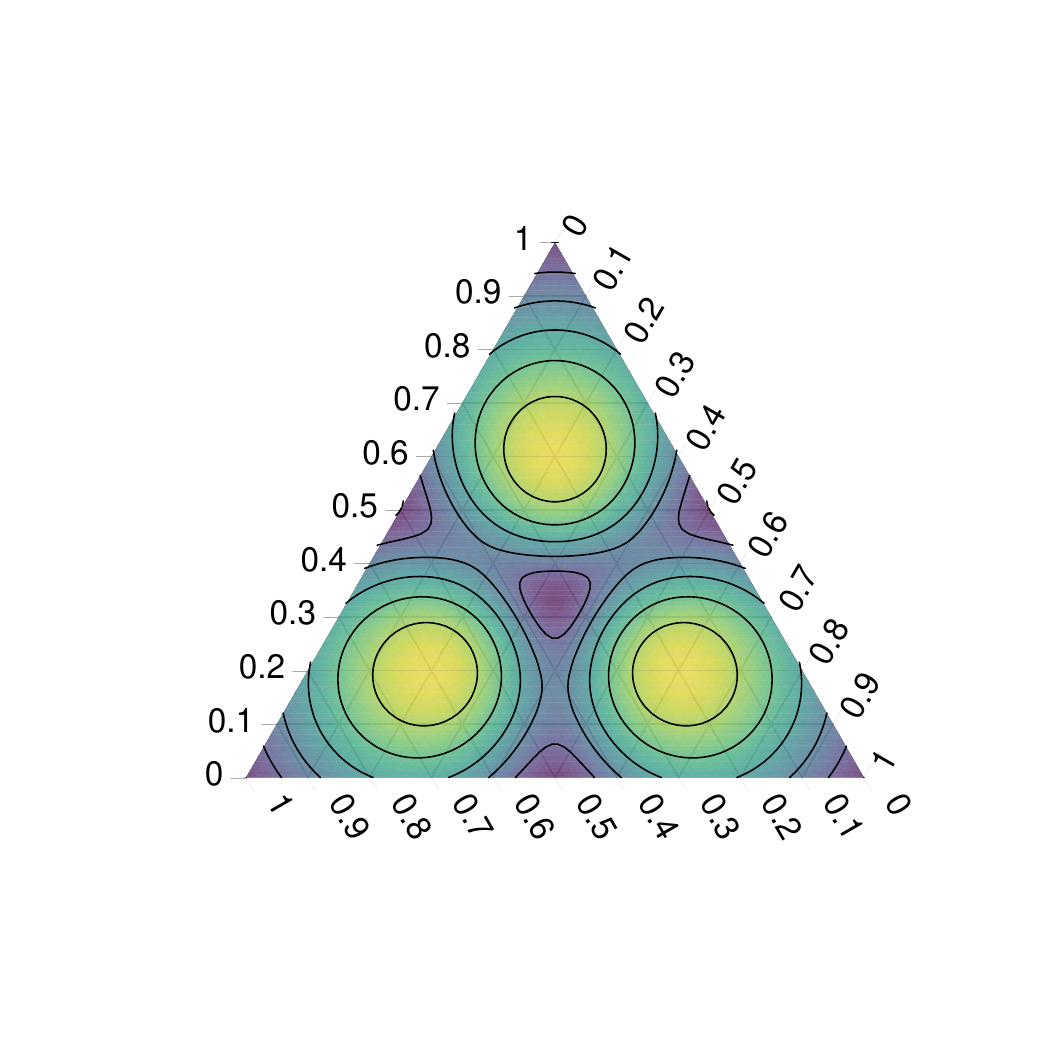}\hspace*{-5mm}
            \includegraphics[trim=0mm -19mm 133mm -5mm, clip, height=5.6cm]{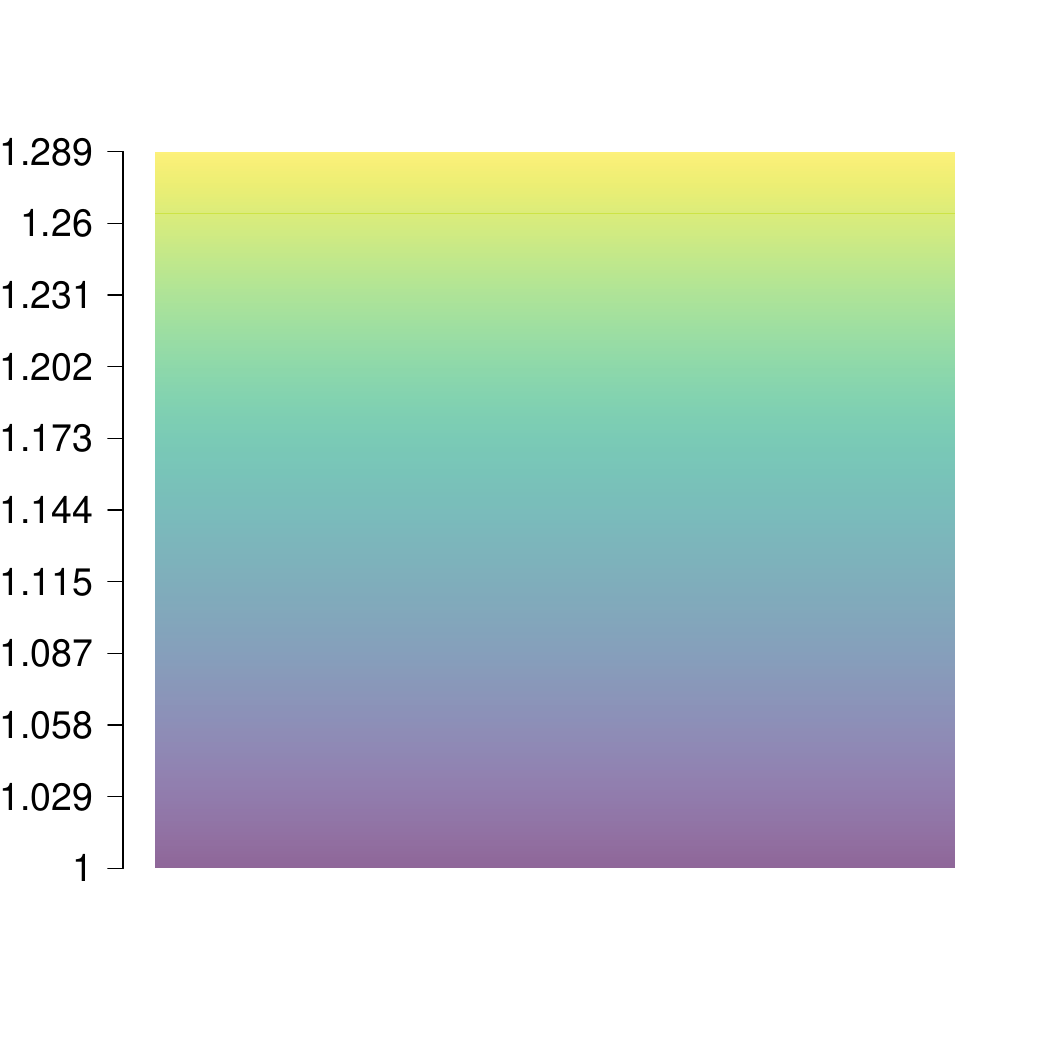}
\vspace*{-12mm} \newline \hspace*{-11mm}
            \includegraphics[height=5.6cm]{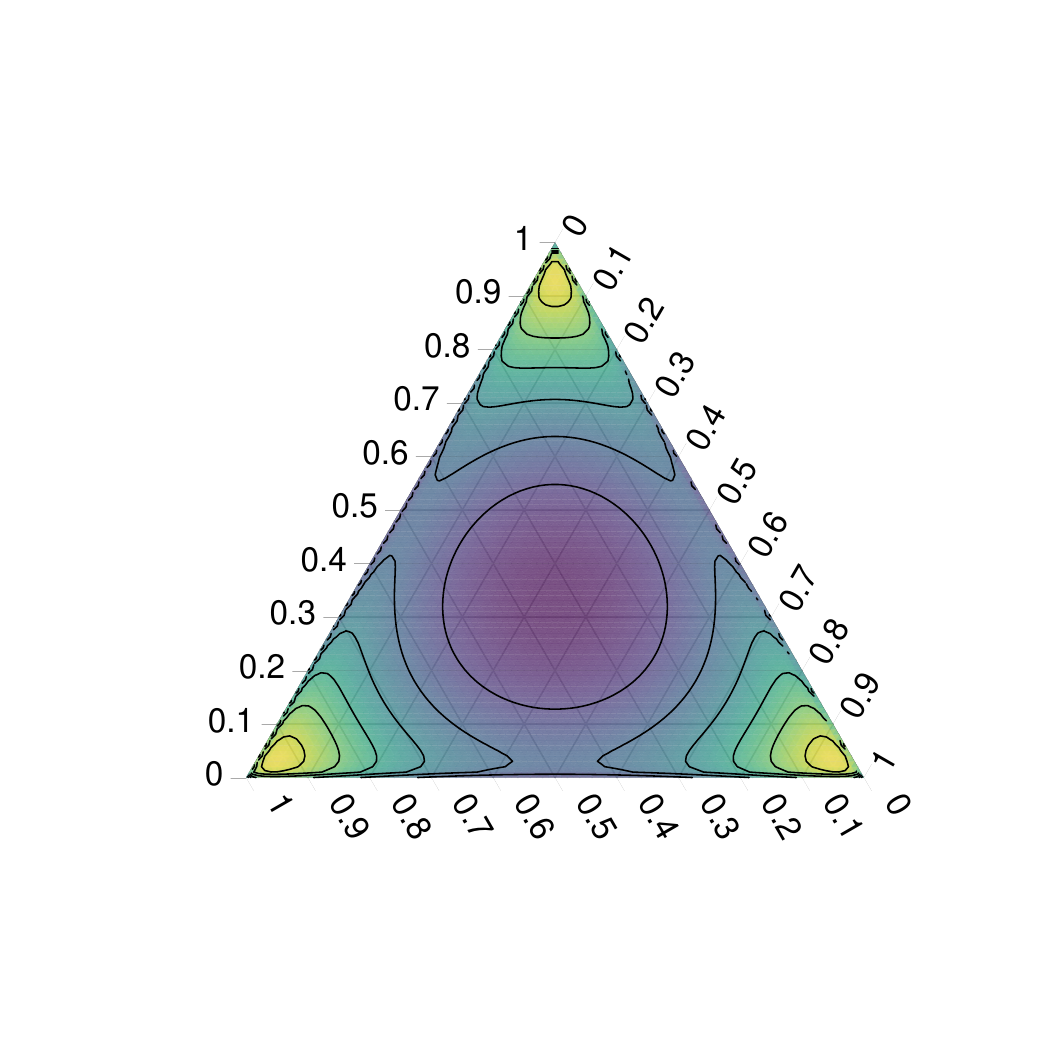}\hspace*{-5mm}
            \includegraphics[trim=0mm -19mm 133mm -5mm, clip, height=5.6cm]{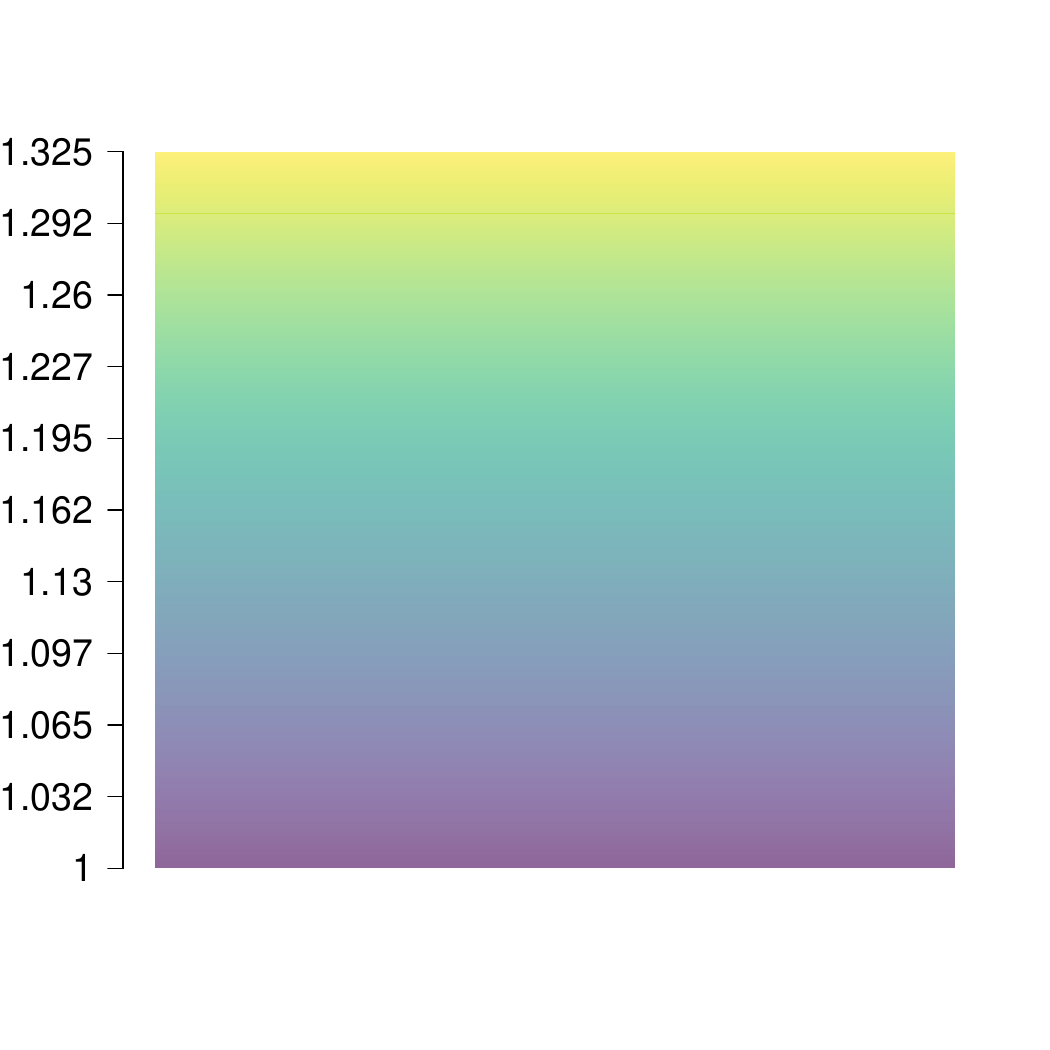}
            \hspace*{-1mm}
            \includegraphics[height=5.6cm]{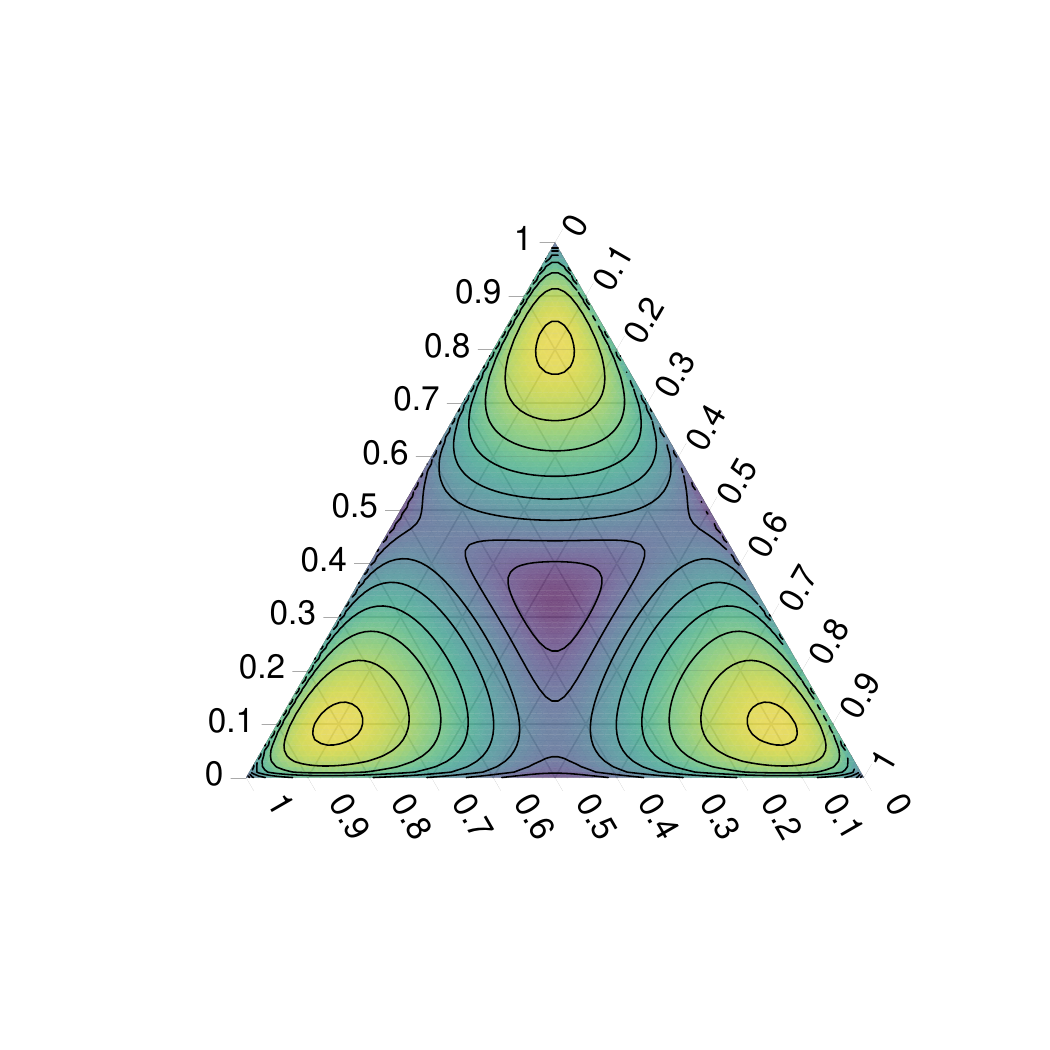}\hspace*{-5mm}
            \includegraphics[trim=0mm -19mm 133mm -5mm, clip, height=5.6cm]{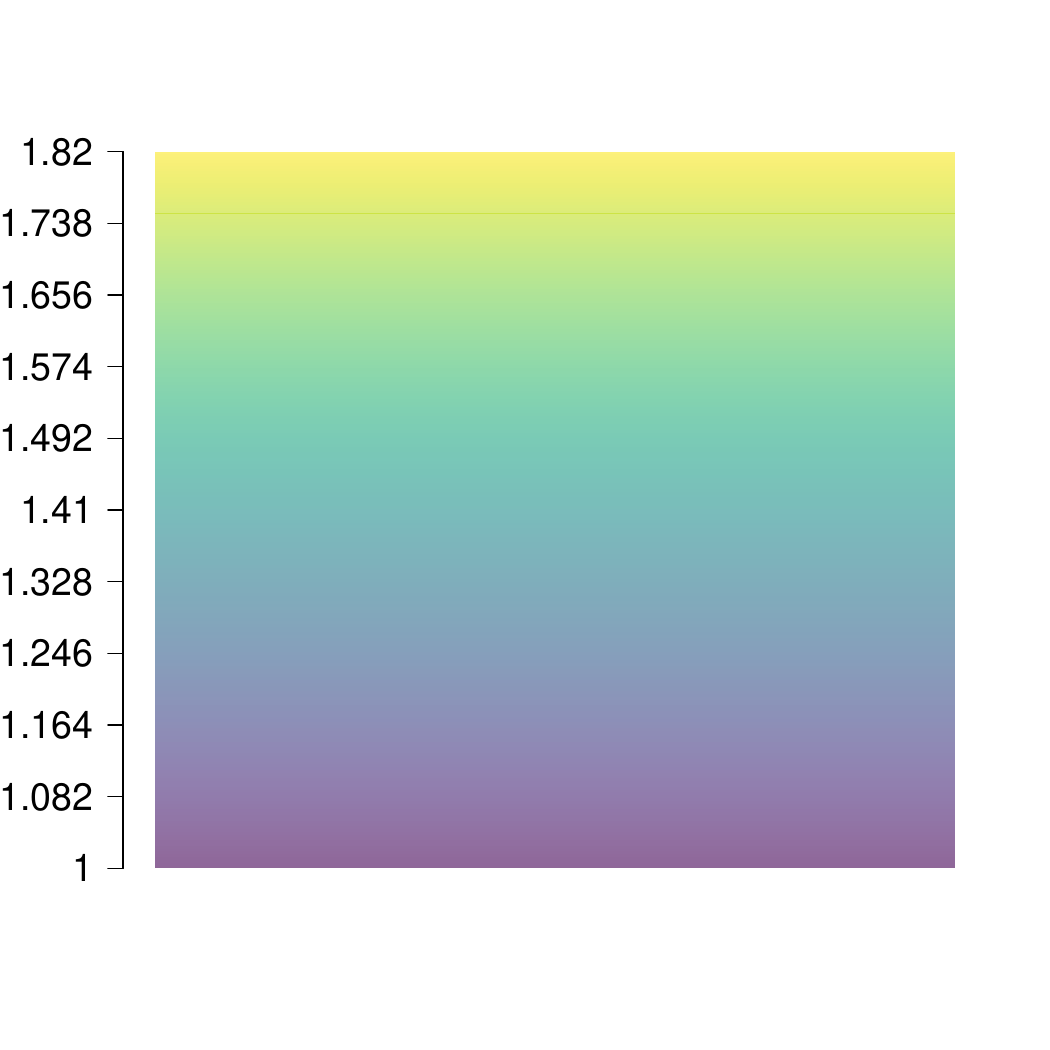}
 \caption{Generalized complexity measure  $C_{\alpha,\beta}(\bar{p})$ for $\bar{p}= (p_1, p_2, p_3)$, with $(\alpha,\beta) = (1,2)$ (top left), (2,10) (top right), (0.5,1) ((bottom left) and (0.5,10) (bottom right). Noting that in all the cases $\alpha$ has been taken lower than $\beta$ (complexity values would be their corresponding multiplicative inverses by interchanging $\alpha$ and $\beta$), the minimum value 1 is always reached for the singular points represented by the center, the three edge middle points, and the three vertices, which in fact represent the invariant ternary distributions under power distortion. The transition patterns, including the locations of the maximal points and their associated complexity values, vary depending on the two-parameter specifications. }  
\label{figure:generalized-complexity(1,2)-(2,10)-(0.5,1)-(0.5,10)}
\end{figure}

In practice, the maps 
$$
\left\{ C_{\alpha, \beta}(\bar{p}): \enspace 0 \leq \alpha, \beta < \infty \right\}, 
$$  
possibly restricted to some subregion within the $(\alpha, \beta)$-parameter space, provide a useful representation of complexity patterns under this approach (and similarly for the continuous case; related formal details and comments are omitted in the remainder of this section). In addition, marginal representations (with a possible rescaling) of 
$$
C_{\alpha, \beta}(\bar{p}) \quad \mbox{vs.} \quad e^{H_{\alpha}(\bar{p})} , \qquad \mbox{and} \qquad
C_{\alpha, \beta}(\bar{p}) \quad \mbox{vs.} \quad e^{H_{\beta}(\bar{p})},
$$ 
for fixed values of $\beta$ and $\alpha$, respectively, are of interest to assess the relative contribution of the factors to the product complexity.  

For $\alpha = \beta$ (which would correspond to the map main diagonal) we invariably have the constant value $C_{\alpha, \beta}(\cdot) \equiv 1$, which is not informative. Observing that the generalized complexity measure $C_{\alpha, \beta}(\bar{p})$ in effect quantifies, in an appropriate scale, the incremental (Rényi) entropic response of the distribution $\bar{p}$ with respect to changes of the deformation parameter, we can justify the complementary use, for complexity assessment, of the relative increment function
$$
\ln \left[ \left(C_{\alpha, \beta}(\bar{p}) \right)^{\frac{1}{\alpha - \beta}} \right] = \frac{H_{\alpha}(\bar{p})-H_{\beta}(\bar{p})}{\alpha  - \beta}, \qquad \mbox{(or} \quad \left(C_{\alpha, \beta}(\bar{p}) \right)^{\frac{1}{\alpha - \beta}}  = e^{\frac{H_{\alpha}(\bar{p})-H_{\beta}(\bar{p})}{\alpha  - \beta}}),
$$
and, in particular, of its limit as $\beta \rightarrow \alpha$ for fixed $\alpha$, that is, the  (curve of) derivatives of $H_{\alpha}(\bar{p})$ with respect to $\alpha$, 
$$
H^{\prime}_{\alpha}(\bar{p}), \qquad \mbox{(or, respectively,} \quad e^{H^{\prime}_{\alpha}(\bar{p})})
$$
(optionally, in both cases, displayed with a change of sign for positiveness, according to the decreasing monotonicity of Rényi entropy with respect to the deformation parameter). This quantity is indeed informative about the entropic sensitivity of $\bar{p}$ to local infinitesimal changes of the deformation parameter. (This idea is introduced in Esquivel, Alonso and Angulo 2017, in relation to multifractal systems and generalized R\'enyi dimensions; see section \ref{section:multifractality} below.)

Figure \ref{figure:Renyi-entropy-derivative} shows the values of the derivative of Rényi entropy, $H^{\prime}_{\alpha}(\bar{p})$, at $\alpha = 2$, for a ternary discrete probability distribution $\bar{p}= (p_1, p_2, p_3)$.

\begin{figure}
            \includegraphics[height=8.0cm]{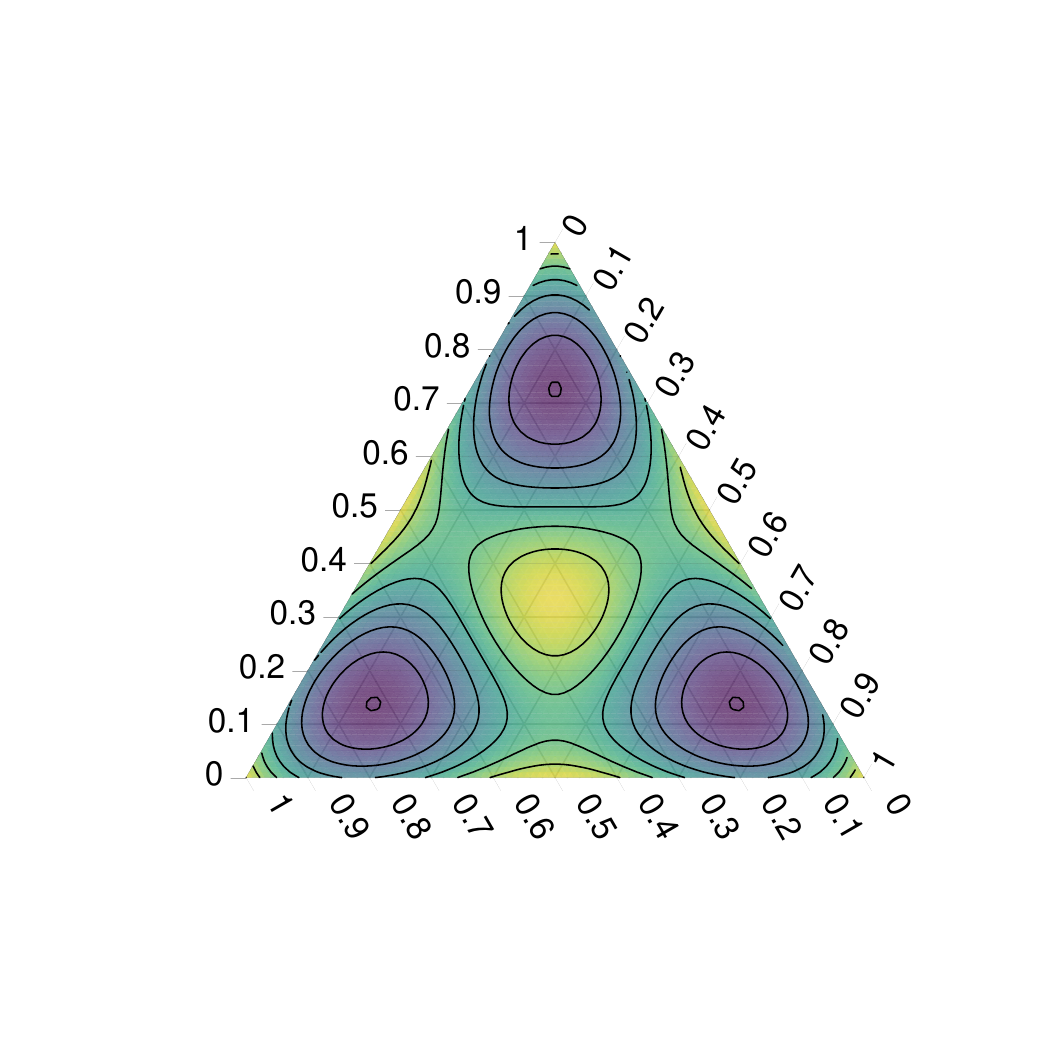}\hspace*{1mm}
            \includegraphics[trim=0mm -19mm 133mm -5mm, clip, height=8.0cm]{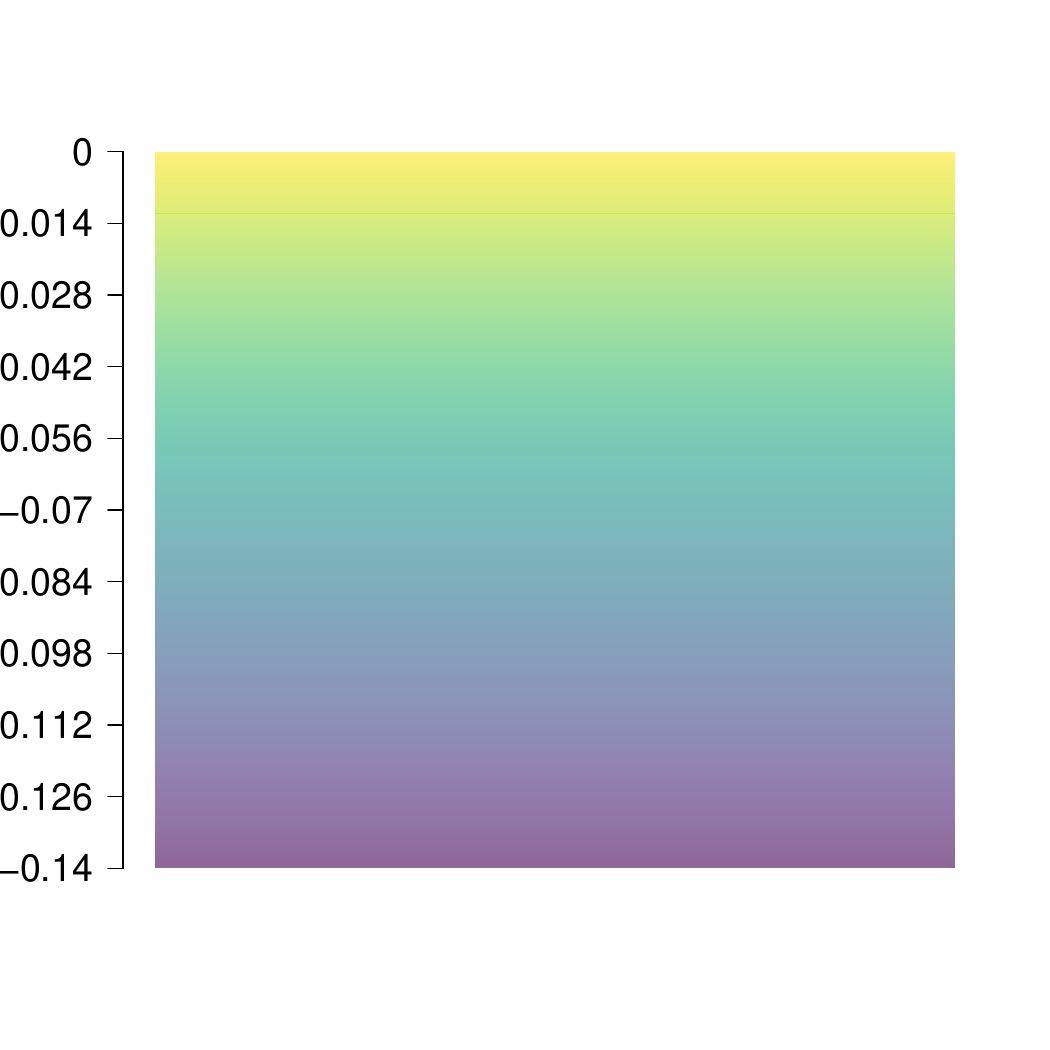}
\caption{Derivative of Rényi entropy, $H_{\alpha}^{\prime}$,  at $\alpha = 2$, for $\bar{p}= (p_1, p_2, p_3)$. The structure of this plot shows the intrinsic value of $-H_{\alpha}^{\prime}$  as a one-parameter generalized complexity measure, under the notion of `complexity' as departure from both equilibrium and degeneracy.}
\label{figure:Renyi-entropy-derivative}           
\end{figure}

\end{paragraph}

\subsection{Relative Entropy and Relative Complexity}
\label{subsection:relative-entropy-complexity}

\begin{paragraph}{Generalized relative complexity -- Romera, Sen and Nagy (2011)}

In analogy with the product-type formulation of the two-parameter generalized complexity measure referred in subsection \ref{subsection:entropy-complexity}, Romera, Sen and Nagy (2011) proposed (as before, in the context of Physics) the formulation of a measure of (directed) relative complexity between two given probability densities,  $\left\{f(\mathbf{x}): \mathbf{x} \in \mathds{R}^d \right\}$ and $\left\{g(\mathbf{x}): \mathbf{x} \in \mathds{R}^d \right\}$ (with $f$ being absolutely continuous with respect to $g$), as  
$$
C_{\alpha, \beta}(f\|g) := e^{H_{\alpha}(f\|g)-H_{\beta}(f\|g)}, 
$$
for $0 < \alpha, \beta < \infty$ (again, the definition can be extended for any real value of the parameters, as considered in section \ref{subsection:relative-entropy-complexity-multifractality}). This measure is then based on the `local' (state-by-state) comparison of the distributions involved, specifically quantifying the sensitivity of Rényi divergence with respect to changes in the deformation parameter.  

As before, this measure can be meaningfully adopted for the discrete case: For two discrete probability distributions $\bar{p}_1=(p_{11},p_{12},...,p_{1n})$ and $\bar{p}_2=(p_{21},p_{22},...,p_{2n})$ on a given set of $n$ states,
$$ 
C_{\alpha, \beta}(\bar{p}_1\|\bar{p}_2) := e^{H_{\alpha}(\bar{p}_1\|\bar{p}_2)-H_{\beta}(\bar{p}_1\|\bar{p}_2)}. 
$$
%

In terms of the relative diversity index introduced in section \ref{section:information-entropy-divergence}, we can therefore also rewrite
$$
C_{\alpha, \beta}(\bar{p}_1\|\bar{p}_2) = \frac{DI_{\alpha}(\bar{p}_1\|\bar{p}_2)}{DI_{\beta}(\bar{p}_1\|\bar{p}_2)}, \qquad C_{\alpha, \beta}(f\|g) = \frac{DI_{\alpha}(f\|g)}{DI_{\beta}(f\|g)},
$$ 
with the correspondingly added interpretation.

Plots in Figure \ref{figure:generalized-relative-complexity(1,2)-(2,10)-(0.5,1)-(0.5,10)_020305} display the ternary simplex representations of the generalized relative complexity values with varying $\bar{p}= (p_1, p_2, p_3)$ and fixed $\bar{q}= (0.2, 0.3, 0.5)$, for different parameter values, namely $(\alpha,\beta) = (1,2)$, $(2,10)$, $(0.5,1)$, $(0.5, 10)$,  respectively. 
The non-symmetric character of the measure with respect to the roles of $\bar{p}$ and $\bar{q}$, as is intrinsic to the definition of Rényi divergence, can be observed comparing the two plots in Figure \ref{figure:generalized-relative-complexity(0.5,1)_q020305-(0.5,1)_p020305}, where in one case the simplex coordinates correspond to varying $\bar{p}= (p_1, p_2, p_3)$, with $\bar{q}= (0.2, 0.3, 0.5)$ as a fixed distribution, and conversely, in the other case the simplex coordinates correspond to varying $\bar{q}= (q_1, q_2, q_3)$, with $\bar{p}= (0.2, 0.3, 0.5)$ as the fixed distribution, both scenarios under the same parameter specifications, $(\alpha,\beta) = (0.5,1)$.

\begin{figure} \hspace*{-11mm}
            \includegraphics[height=5.6cm]{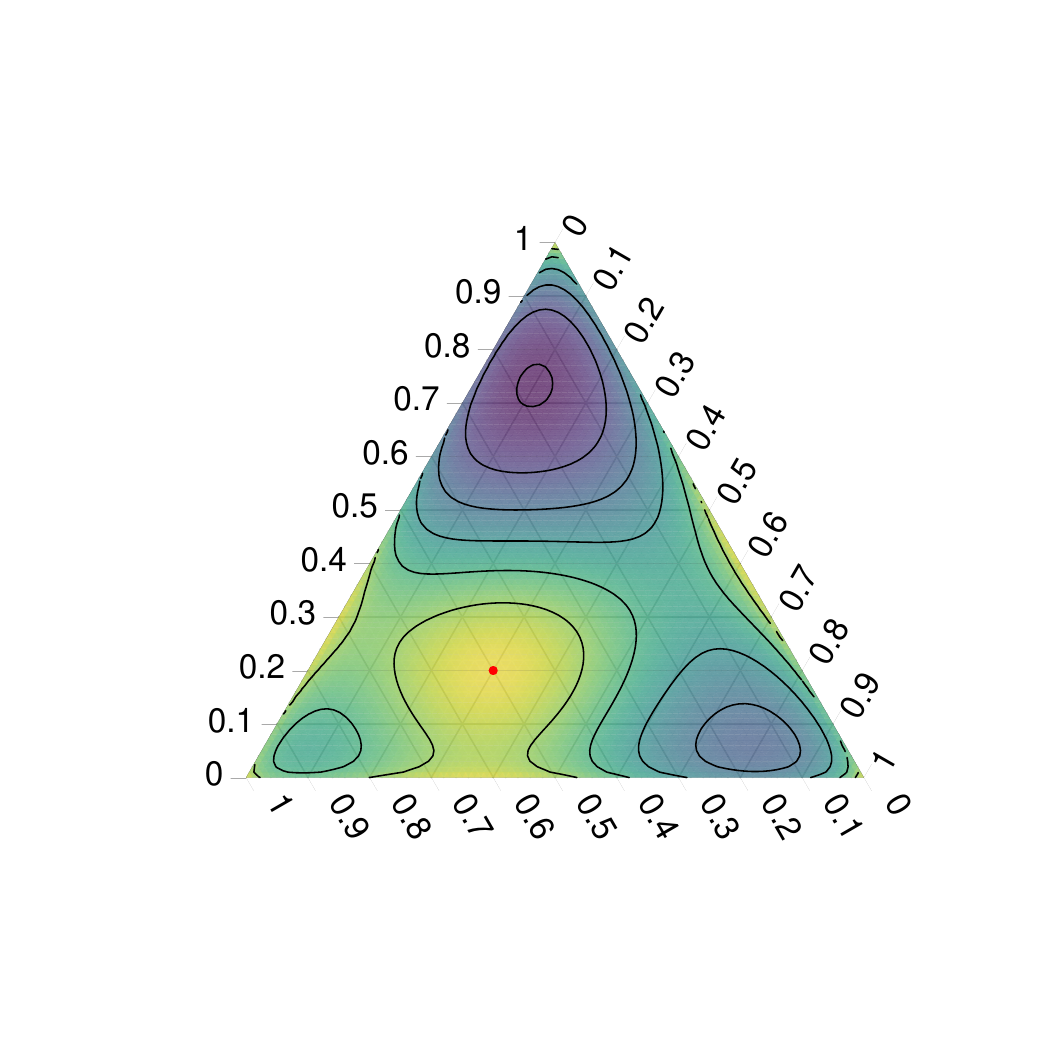}\hspace*{-5mm}
            \includegraphics[trim=0mm -19mm 133mm -5mm, clip, height=5.6cm]{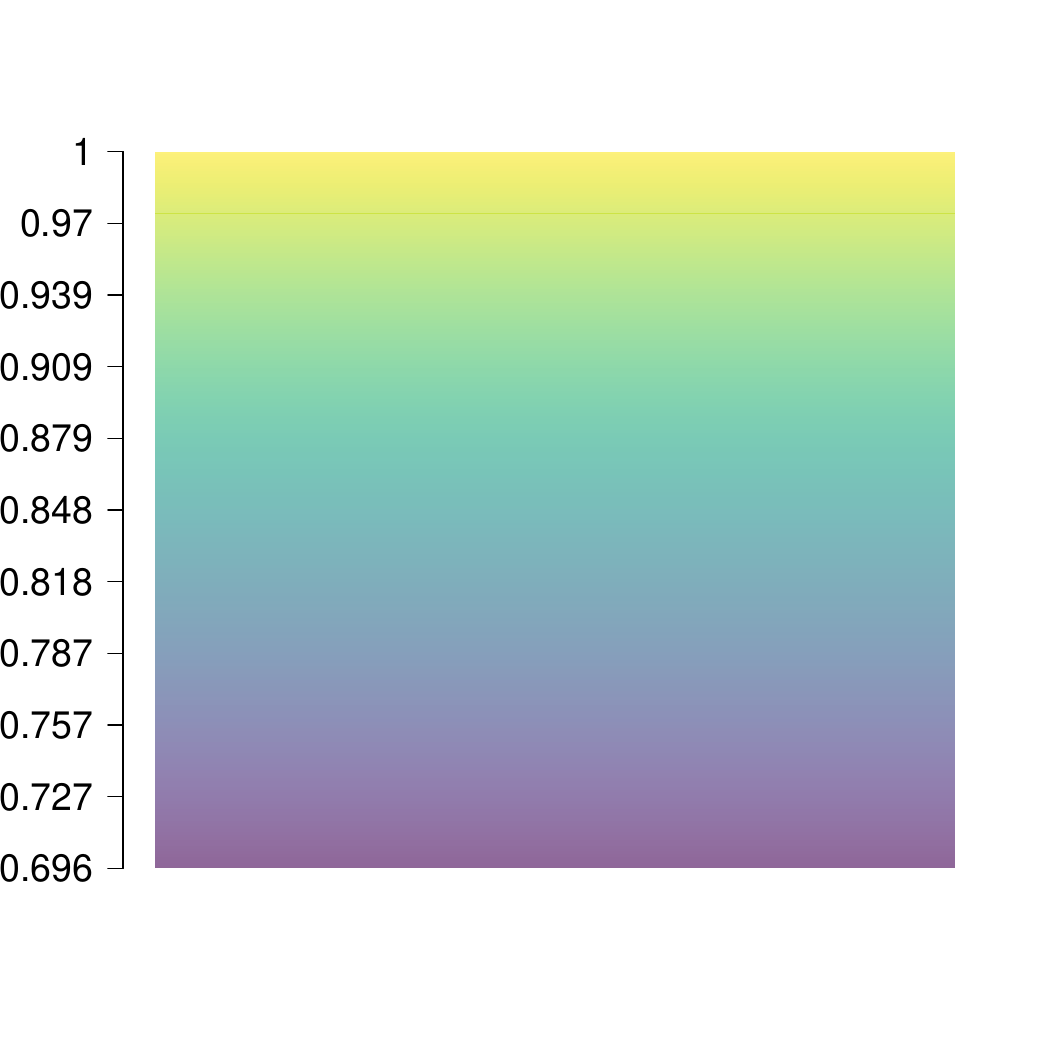}
\hspace*{-1mm}
            \includegraphics[height=5.6cm]{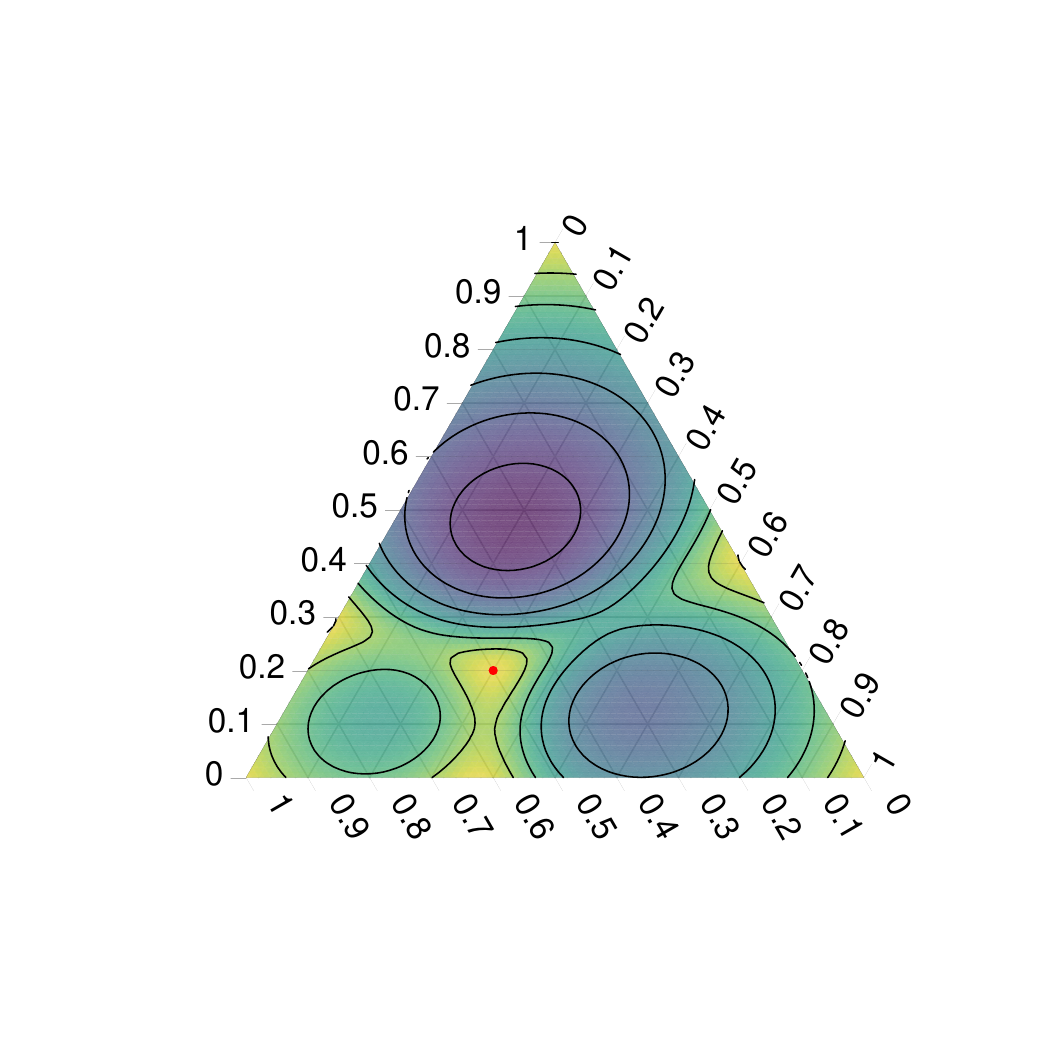}\hspace*{-5mm}
            \includegraphics[trim=0mm -19mm 133mm -5mm, clip, height=5.6cm]{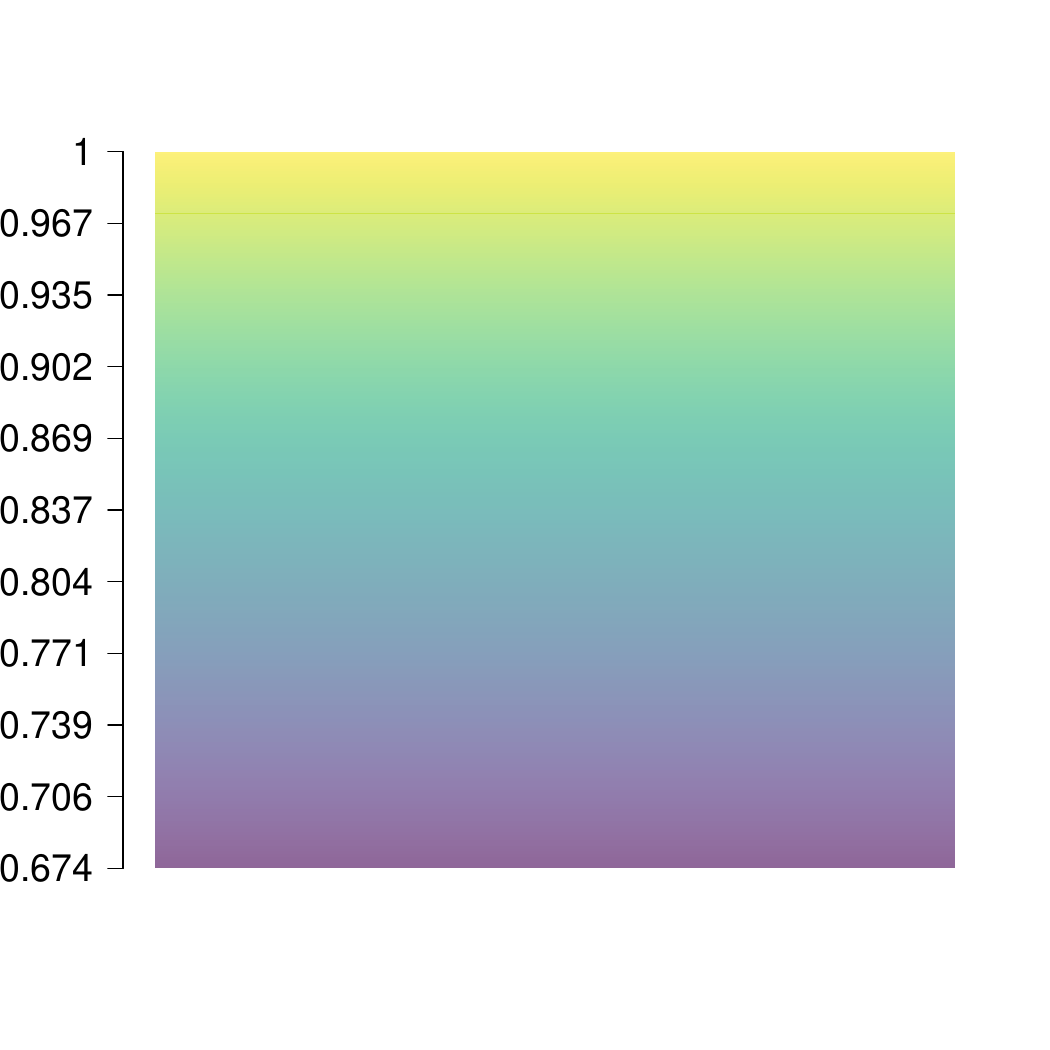}
\vspace*{-12mm} \newline \hspace*{-11mm}            
            \includegraphics[height=5.6cm]{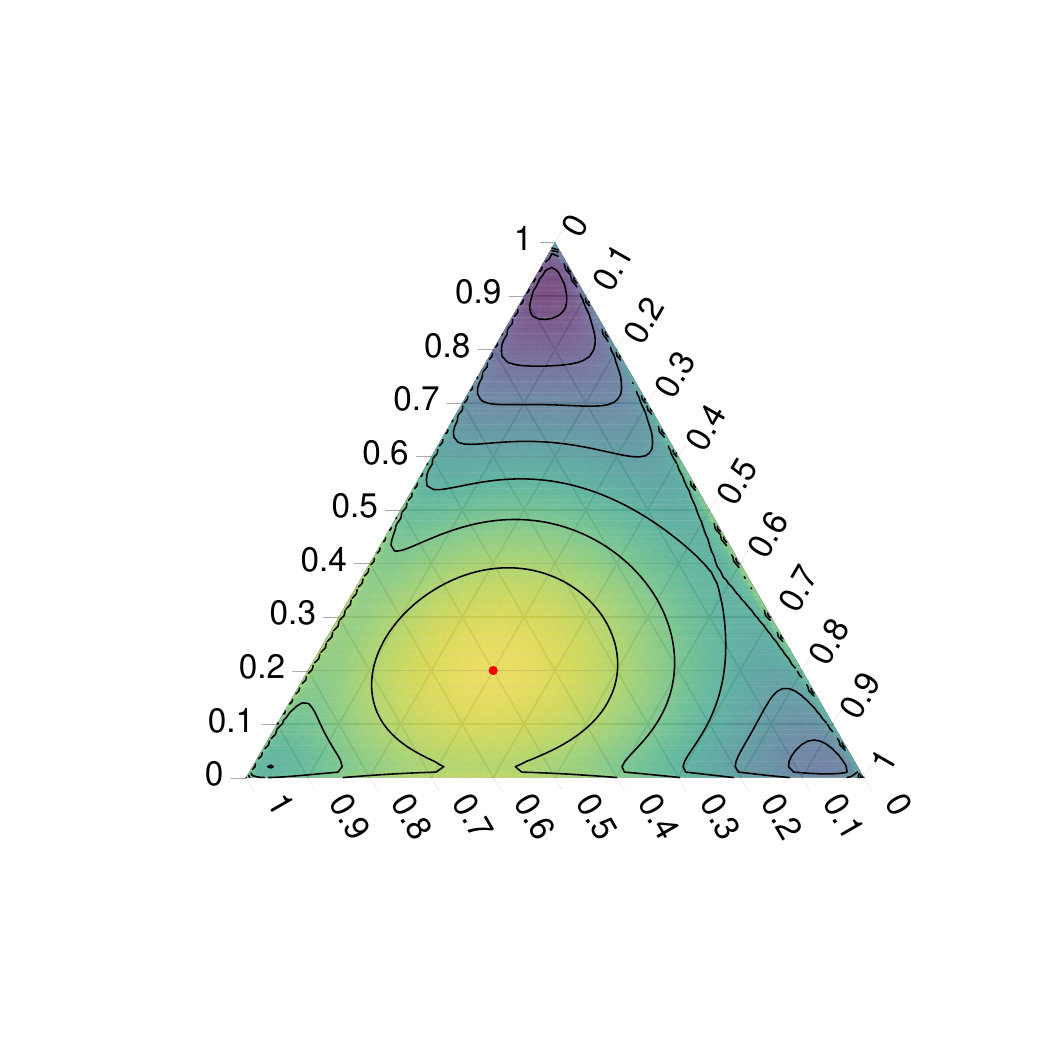}\hspace*{-5mm}
            \includegraphics[trim=0mm -19mm 133mm -5mm, clip, height=5.6cm]{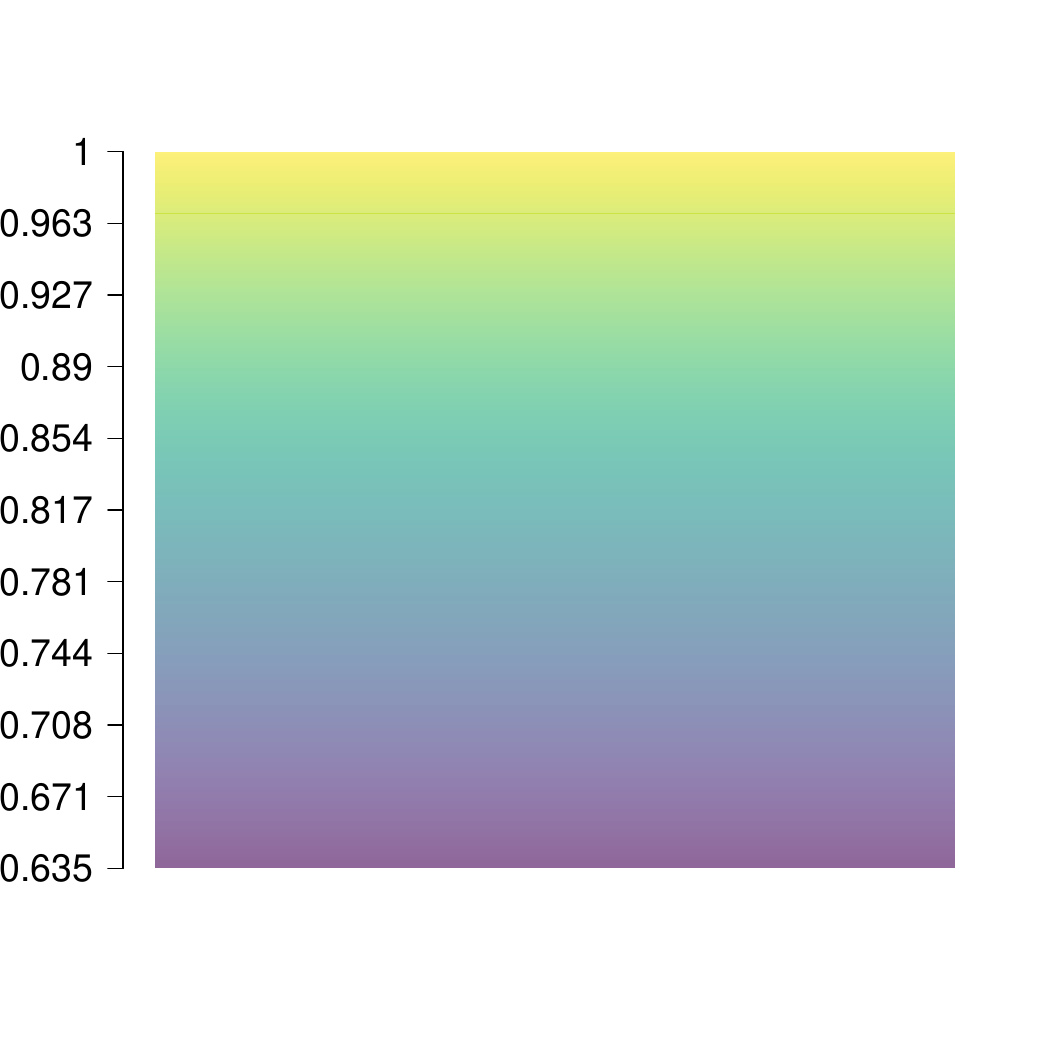}
\hspace*{-1mm}
            \includegraphics[height=5.6cm]{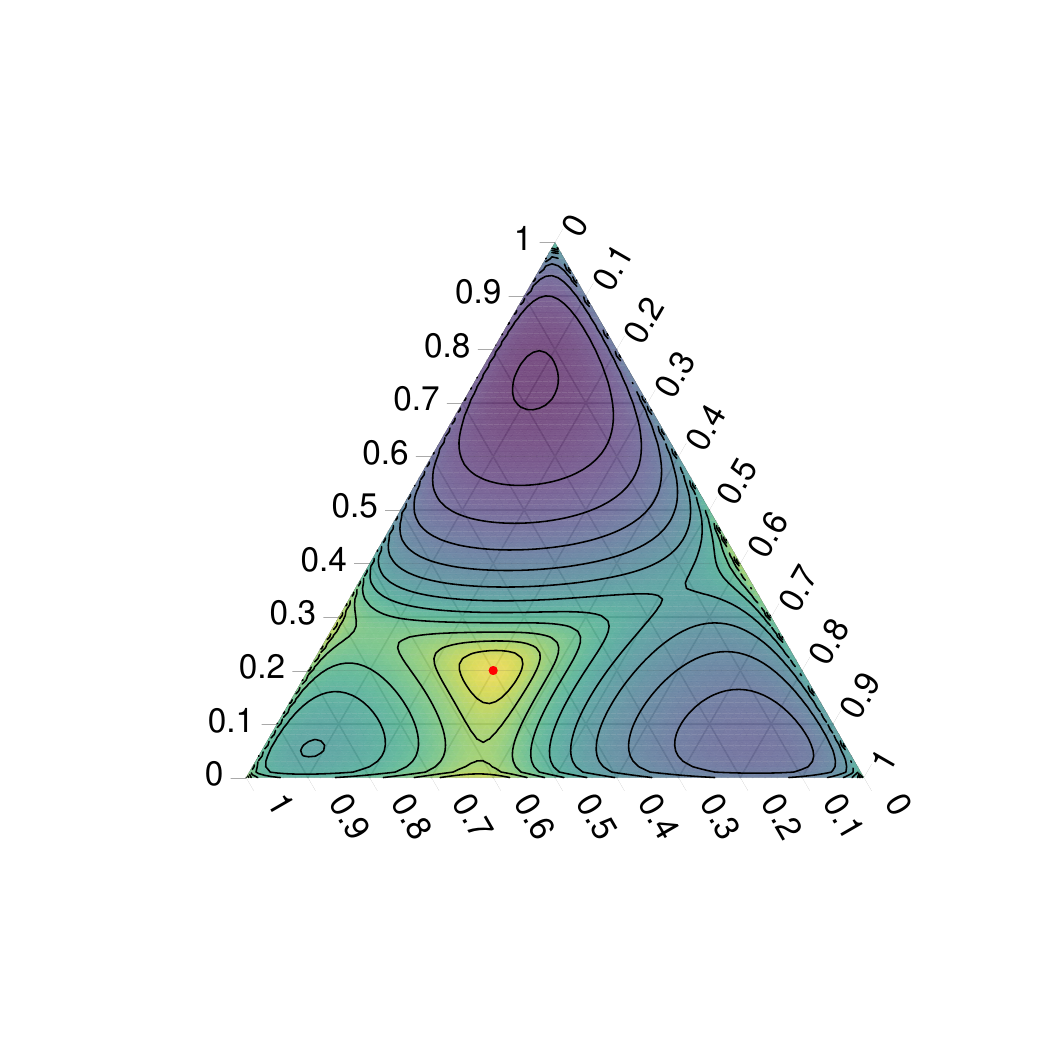}\hspace*{-5mm}
            \includegraphics[trim=0mm -19mm 133mm -5mm, clip, height=5.6cm]{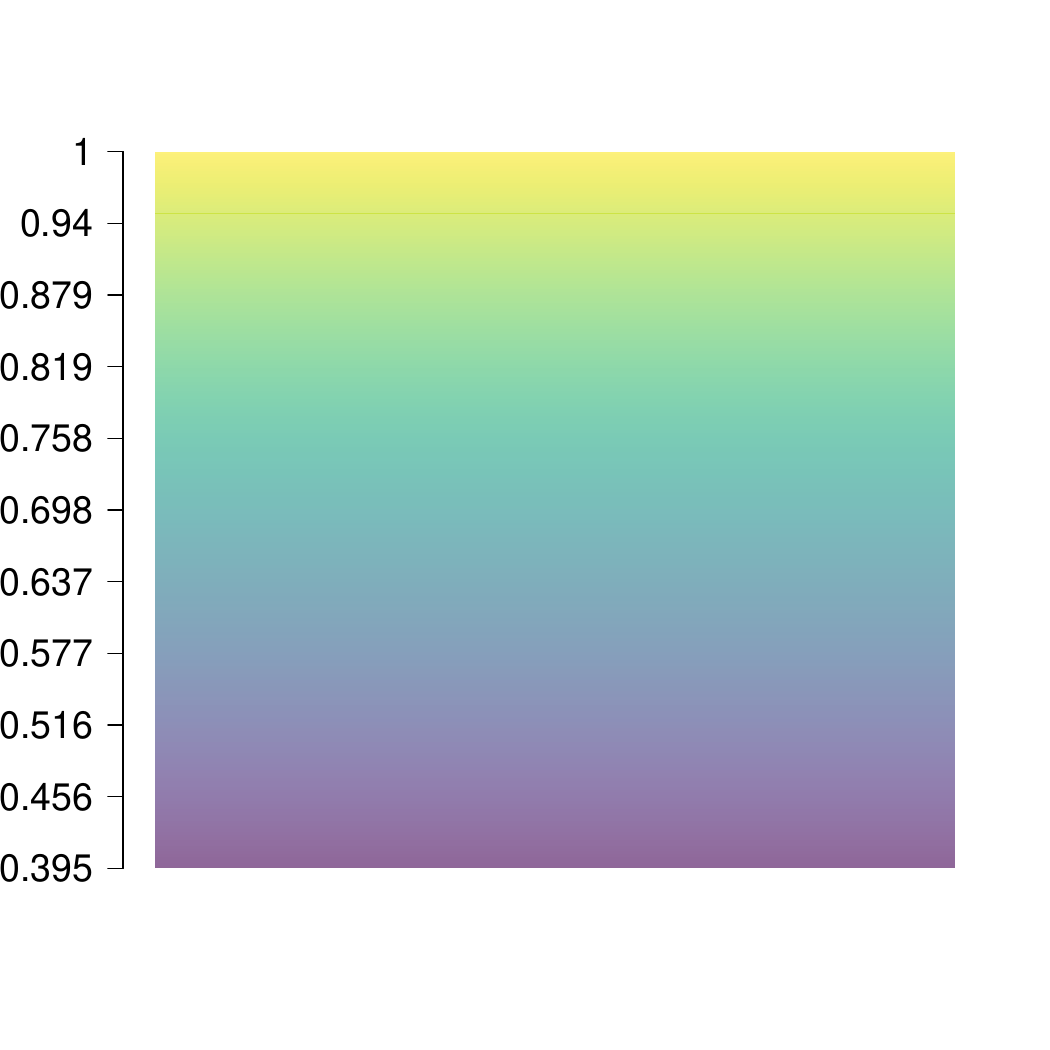}
\caption{Generalized relative complexity measure  $C_{\alpha,\beta}(\bar{p}\|\bar{q})$ for $\bar{p}= (p_1, p_2, p_3)$, $\bar{q}= (0.2, 0.3, 0.5)$, with $(\alpha,\beta) = (1,2)$ (top left), (2,10) (top right), (0.5,1) (bottom left) and (0.5,10) (bottom righ). In all the cases, the maximum value 1 is reached for $\bar{p}$ equal to the fixed reference distribution $\bar{q}$. Structural patterns vary depending on the two-parameter specifications, with lack of triangular symmetry due to departure of $\bar{q}$ from the central equiprobability distribution.}  
\label{figure:generalized-relative-complexity(1,2)-(2,10)-(0.5,1)-(0.5,10)_020305}           
\end{figure}

\begin{figure} \hspace*{-11mm}            
            \includegraphics[height=5.6cm]{Figures-Measures/C_05_1_020305.pdf}\hspace*{-5mm}
            \includegraphics[trim=0mm -19mm 133mm -5mm, clip, height=5.6cm]{Figures-Measures/BarraColor_C_05_1_020305.pdf}
\hspace*{-1mm}
            \includegraphics[height=5.6cm]{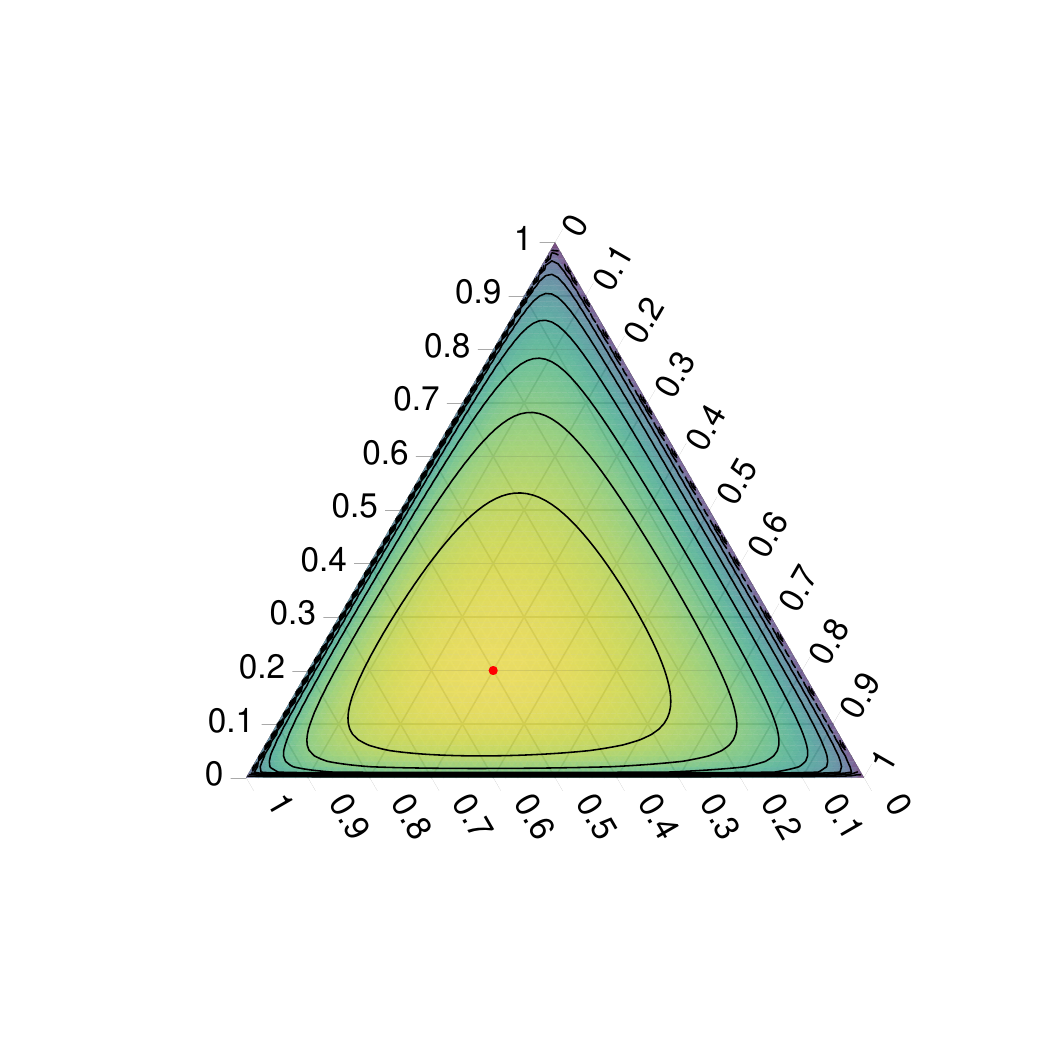}\hspace*{-5mm}
            \includegraphics[trim=0mm -19mm 133mm -5mm, clip, height=5.6cm]{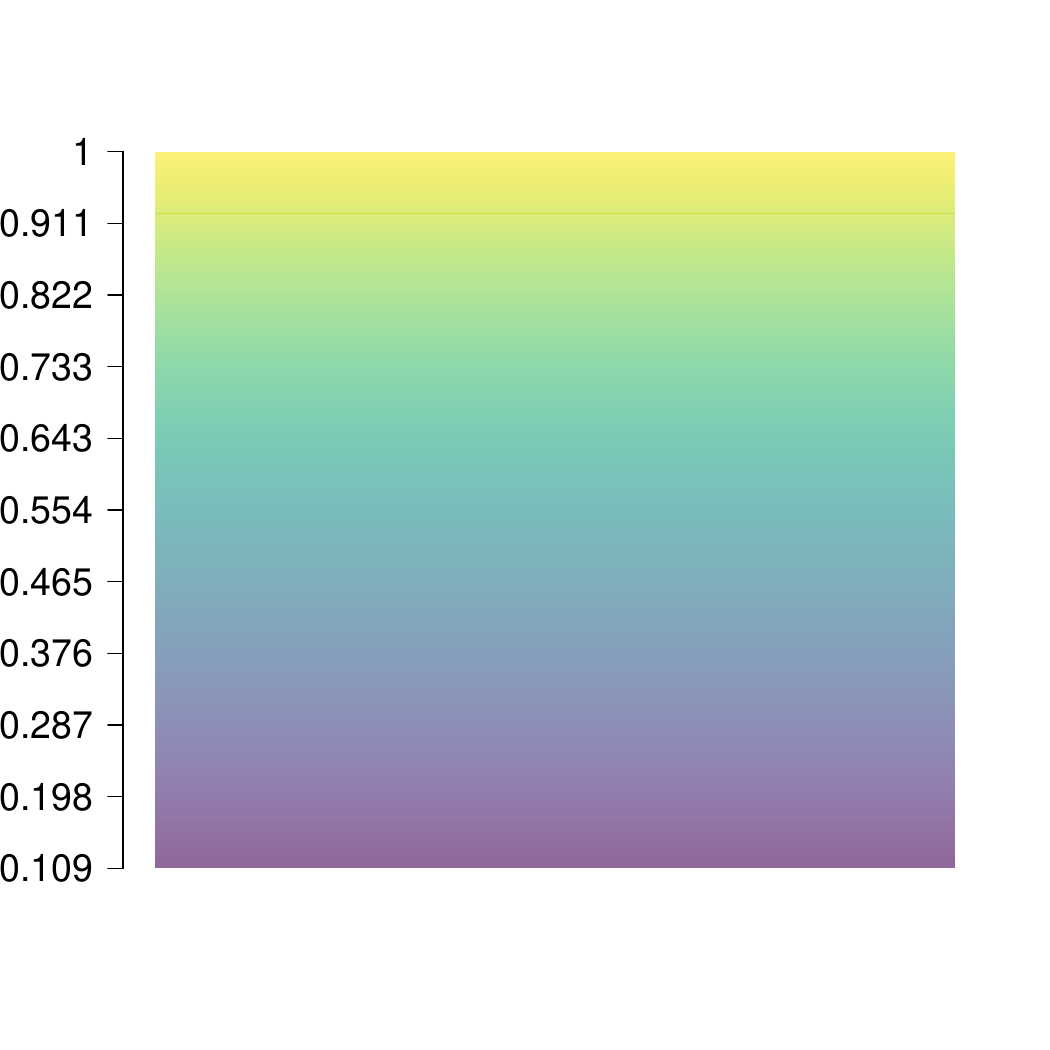}
\caption{Generalized relative complexity measure  $C_{\alpha,\beta}(\bar{p}\|\bar{q})$  for varying $\bar{p}= (p_1, p_2, p_3)$ and fixed $\bar{q}= (0.2, 0.3, 0.5)$ (left), and for varying $\bar{q}= (q_1, q_2, q_3)$ and fixed $\bar{p}= (0.2, 0.3, 0.5)$ (right), in both cases with $(\alpha,\beta) = (0.5,1)$. The different patterns in these plots reflect the non-symmetric roles of the argument distributions $\bar{p}$ and $\bar{q}$, in agreement with the directed nature of Rényi divergence on which the generalized relative complexity measure is based.}  
\label{figure:generalized-relative-complexity(0.5,1)_q020305-(0.5,1)_p020305}  
\end{figure}

In the particular case where the reference distribution is uniform, i.e. for $\bar{p}_2 \equiv \left[\frac{1}{n}\right]$, the generalized relative complexity measure is equivalent to the generalized complexity measure in the following reciprocal sense: 
$$
C_{\alpha, \beta}(\bar{p}_1\|\left[\frac{1}{n}\right]) = \frac{1}{C_{\alpha, \beta}(\bar{p}_1)} = C_{\beta, \alpha}(\bar{p}_1).
$$
   
Along the same lines as described in subsection \ref{subsection:entropy-complexity}, maps based on (relative) increments and, in particular, the curve of Rényi divergence derivatives with respect to the deformation parameter are useful for relative complexity assessment (further details are omitted; see subsection \ref{subsection:relative-entropy-complexity-multifractality} in relation to generalized relative dimensions).

\end{paragraph}

\section{Multifractality: Generalized Dimensions and Complexity}
\label{section:multifractality}

In this section, in the context of multifractality, we first recall the meaningful limiting functional connection between the increments of generalized Rényi dimensions and the two-parameter generalized complexity measure described in subsection \ref{subsection:entropy-complexity}, as established in Angulo and Esquivel (2014) and Esquivel, Alonso and Angulo (2017). Thereafter, from a natural divergence-based formulation of generalized Rényi relative dimensions, regarding the assessment of the local relation between two given multifractal measures, a parallel functional connection is determined, in this case, with the two-parameter generalized relative complexity measure described in \ref{subsection:relative-entropy-complexity}. The instrumental significance of related  maps of increments and, in particular, curves of derivatives, complementarily from both --ordinary and relative-- forms of generalized dimensions, is illustrated in section \ref{section:example}.

\subsection{Entropy and Complexity -- Multifractality}
\label{subsection:entropy-complexity-multifractality}

\begin{paragraph}{Generalized  R\'enyi dimensions -- Hentschel and Procaccia (1983), Grassberger (1983)} 
For a multifractal measure $\mu$ with support in the $d$-dimensional Euclidean space,  the \emph{generalized Rényi dimension} (or \emph{generalized fractal dimension}, or simply \emph{generalized  dimension}, etc.)  of order $q$ (with $q \in \mathds{R}$) is defined by the rate of divergence of Rényi entropy of order $q$ (or Shannon entropy in the case $q=1$) of the `partition' probability distributions $P_{\varepsilon}$ associated with lattice coverings of decreasing box width $\varepsilon$, as follows:     
 \begin{align}
  D_q & := \lim_{\varepsilon \rightarrow 0} \frac{1}{q-1} \frac{\ln\left(   \displaystyle\sum_{k\in K_{\varepsilon}} \mu^q[B_{\varepsilon}(k)]\right) }{\ln\left( \varepsilon\right) } \qquad (q\neq1), \nonumber  \\ & \nonumber \\ \protect{\vspace*{20mm}}
   D_1 & := \lim_{\varepsilon \rightarrow 0} \frac{ \displaystyle\sum_{k\in K_{\varepsilon}} \mu[B_{\varepsilon}(k)] \ln \left( \mu[B_{\varepsilon}(k)]\right) }{\ln \left( \varepsilon\right) }, \nonumber
  \end{align} \bigskip
where $K_{\varepsilon}$ denotes the set of non-$\mu$-null boxes of width $\varepsilon$ (see, for example, Harte 2001). 

Motivated by preliminary developments, and due to their specific interpretation, three particular orders have been of common use in applications: $D_0$, the `capacity'; $D_1$, the `information dimension', and $D_2$, the `correlation exponent' (among various other denominations; see, for example, Grassberger 1983, Hentschel and Procaccia 1983). Further, considering that the generalized dimension curve is monotonically non-increasing, and constant for a monofractal measure, the range $D_{-\infty}-D_{\infty}$ is sometimes refereed to as the `multifractal step', reflecting in a certain sense the richness of the multifractal measure.

\end{paragraph}

\begin{paragraph}{Limiting relation between the generalized complexity measure $C_{\alpha,\beta}(\cdot)$ and the generalized Rényi dimensions $D_q(\cdot)$ -- Angulo and Esquivel (2014), Esquivel, Alonso and Angulo (2017)}
From the approximation (as $\varepsilon \rightarrow 0$)
 \begin{align}
  e^{-H_q(P_{\varepsilon})} \sim \varepsilon^{D_q} \qquad (\forall q),\nonumber
  \end{align}
we have 
  $$
  C_{\alpha,\beta}(P_{\varepsilon})=e^{H_{\alpha}(P_{\varepsilon})-H_{\beta}(P_{\varepsilon})} \sim \varepsilon^{D_{\beta}-D_{\alpha}}.
  $$

Accordingly, the incremental function of the generalized dimension curve, $D_{\alpha}-D_{\beta}$, can be interpreted as a generalized complexity measure in the multifractal domain. In practice, the ($\alpha,\beta$)-maps based on 
$$
D_{\alpha}-D_{\beta},  \qquad \quad  \frac{D_{\alpha}-D_{\beta}}{\alpha-\beta},
$$
and, in particular, the curve of derivatives  
$$
D^{\prime}_{\alpha}
$$
(as the limit main diagonal of the latter), are meaningful for complexity assessment.
\end{paragraph}

\subsection{Divergence and Relative Complexity -- Relative Multifractality}
\label{subsection:relative-entropy-complexity-multifractality}

\begin{paragraph}{Generalized relative R\'{e}nyi dimension}
Under the same fundamental justification underlying the concept of divergence for a `local' assessment in the structural comparison of two probability distributions, a definition of generalized relative Rényi dimensions can be introduced in reference to the limiting behaviour of Rényi divergences of different orders for lattice partitionings of decreasing box width. The formulation adopted here (in contrast to related proposals in the literature) is consistent with the fact that, for two given distribution,  Rényi divergence is a non-decreasing function of the deformation parameter $q$:  For two (multifractal) measures $\mu_1$ and $\mu_2$ (with $\mu_1$ assumed to be absolutely continuous with respect to $\mu_2$), the generalized relative Rényi dimension of order $q$ of $\mu_1$ with respect to $\mu_2$ is defined by 
  \begin{align}
  D_q(\mu_1\|\mu_2) & := \lim_{\varepsilon \rightarrow 0} \frac{1}{1-q} \frac{\ln \left(  \displaystyle\sum_{k\in K_{\varepsilon}} \mu_1^q[B_{\varepsilon}(k)]]\mu_2^{1-q}[B_{\varepsilon}(k)]\right) }{\ln\left( \varepsilon\right) } \qquad (q\neq1), \nonumber \\ 
   D_1(\mu_1\|\mu_2) & := \lim_{\varepsilon \rightarrow 0} \frac{ \displaystyle\sum_{k\in K_{\varepsilon}} \mu_1[B_{\varepsilon}(k)] \ln \left( \frac{\mu_1[B_{\varepsilon}(k)]}{\mu_2[B_{\varepsilon}(k)]} \right)}{\ln \left( \varepsilon\right) }, \nonumber
  \end{align}
where, as before, $K_{\varepsilon}$ denotes the set of non-$\mu_2$-null boxes of width $\varepsilon$. 

\end{paragraph}
\begin{paragraph}{Limiting relation between the generalized relative complexity measure $C_{\alpha,\beta}(\cdot \| \cdot)$ and the generalized relative Rényi dimensions $D_q(\cdot \| \cdot)$}
Similarly to the argument given in the previous subsection \ref{subsection:entropy-complexity-multifractality}, the following approximation holds as $\varepsilon \rightarrow 0$:
 \begin{align}
C_{\alpha,\beta}(P_{1,\varepsilon}\|P_{2,\varepsilon}) \sim \varepsilon^{D_\beta(\mu_1\|\mu_2)-D_\alpha(\mu_1\|\mu_2)}.\nonumber
  \end{align}
Hence, the incremental function of the generalized relative dimension curve can be interpreted as a generalized relative complexity measure in the multifractal domain. Analogous considerations proceed regarding the usefulness of related $(\alpha,\beta)$-maps and, in particular, the curve of derivatives for relative complexity assessment. 
\end{paragraph}

\section{Real Data Illustration:
Seismic Series -- El Hierro (Canary Islands, Spain)}
\label{section:example}

In this section, we analyse a series of 11.142 seismic events occurred in the area of the volcanic island of El Hierro (Canary Islands, Spain), from July 19, 2011, until January 7, 2012, related with the   well-known high activity episode involving the submarine eruption of October 10, 2011. This phenomenon has been the object of investigation in a number of studies focusing on different aspects and from various perspectives (see, for instance,  Angulo and Esquivel 2014, Esquivel and Angulo 2015, and references therein). 

\begin{figure}
\includegraphics[width=12cm]{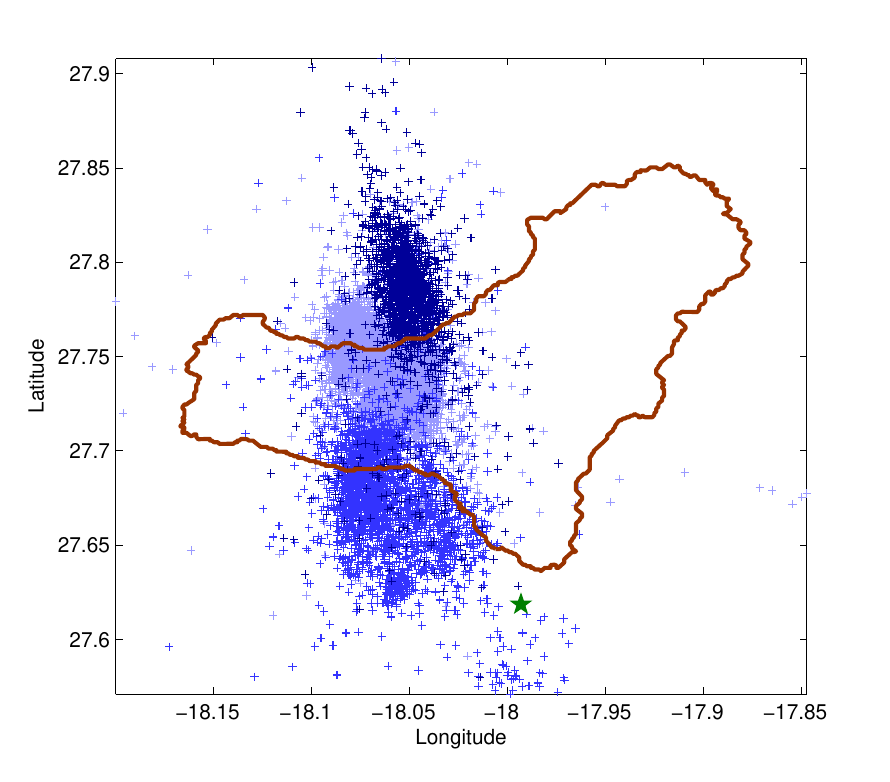}
\includegraphics[width=12cm]{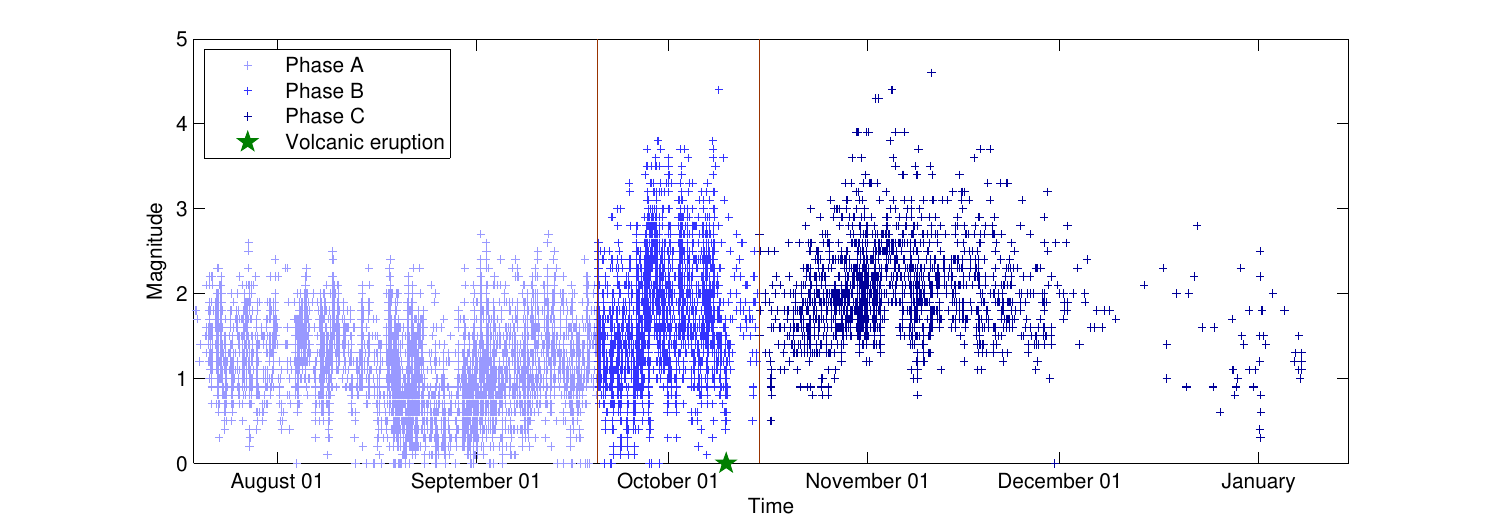}
\caption{El Hierro data: epicenters on contoured island (top), and temporal sequence of magnitudes (bottom). Events corresponding to phases A, B and C are distinguished by different blue colours; the green star corresponds to the main volcanic eruption.}
\label{figure:Hierro-epicenters-magnitudes}
\end{figure}

For comparison purposes, in reference to the evolutionary dynamics, the series is divided into three subperiods (named phases A, B and C), the central of which (phase B) contains the mentioned eruption. Figure \ref{figure:Hierro-epicenters-magnitudes}  (taken from Esquivel and Angulo, 2015) displays the spatial projection of epicenters (top plot) and the temporal sequence of magnitudes (bottom plot) corresponding to the registered events.

In this study, we are interested in assessing possible structural changes, corresponding to the three specified phases, based on the different complexity patterns derived from evaluation of generalized dimensions and generalized relative dimensions, as explained in subsections \ref{subsection:entropy-complexity-multifractality} and \ref{subsection:relative-entropy-complexity-multifractality}.        

\begin{figure}
\includegraphics[width=11cm]{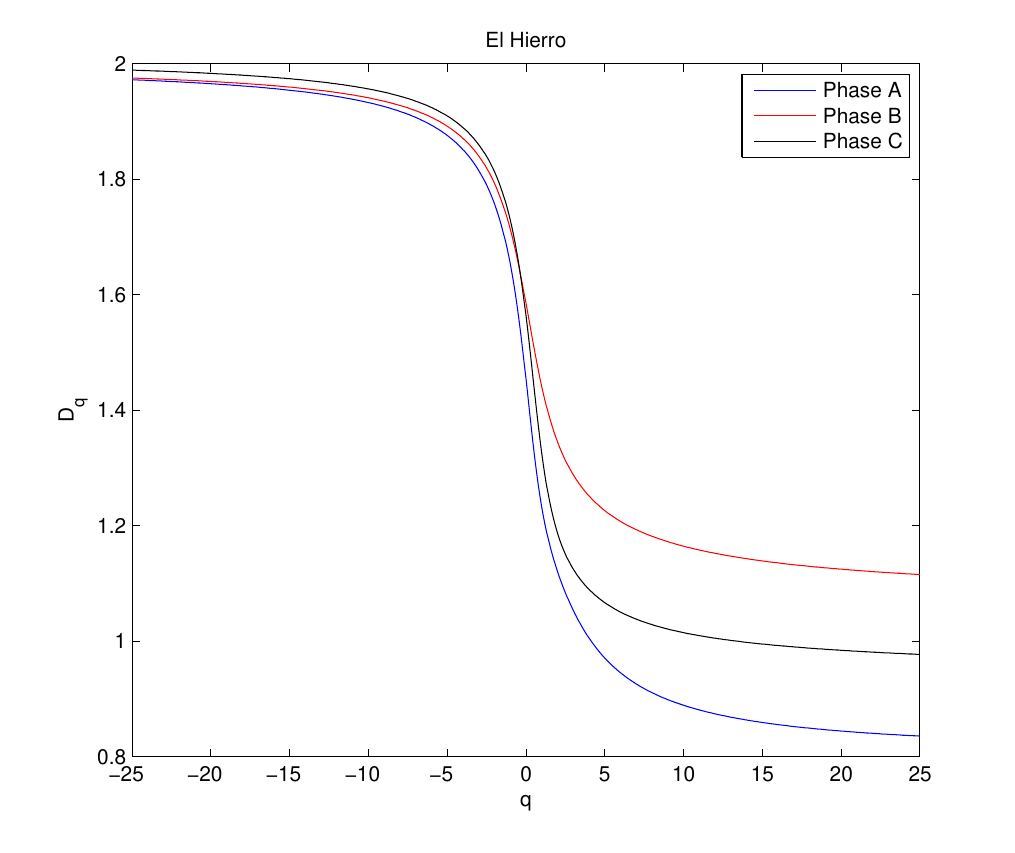}
\centering
\caption{Generalized Rényi dimension curves, showing different multifractal patterns for the three subperiods (phases A, B, and C). }
\label{figure:Hierro-generalized-dimensions}
\end{figure}

\begin{figure}
\includegraphics[width=3.9cm]{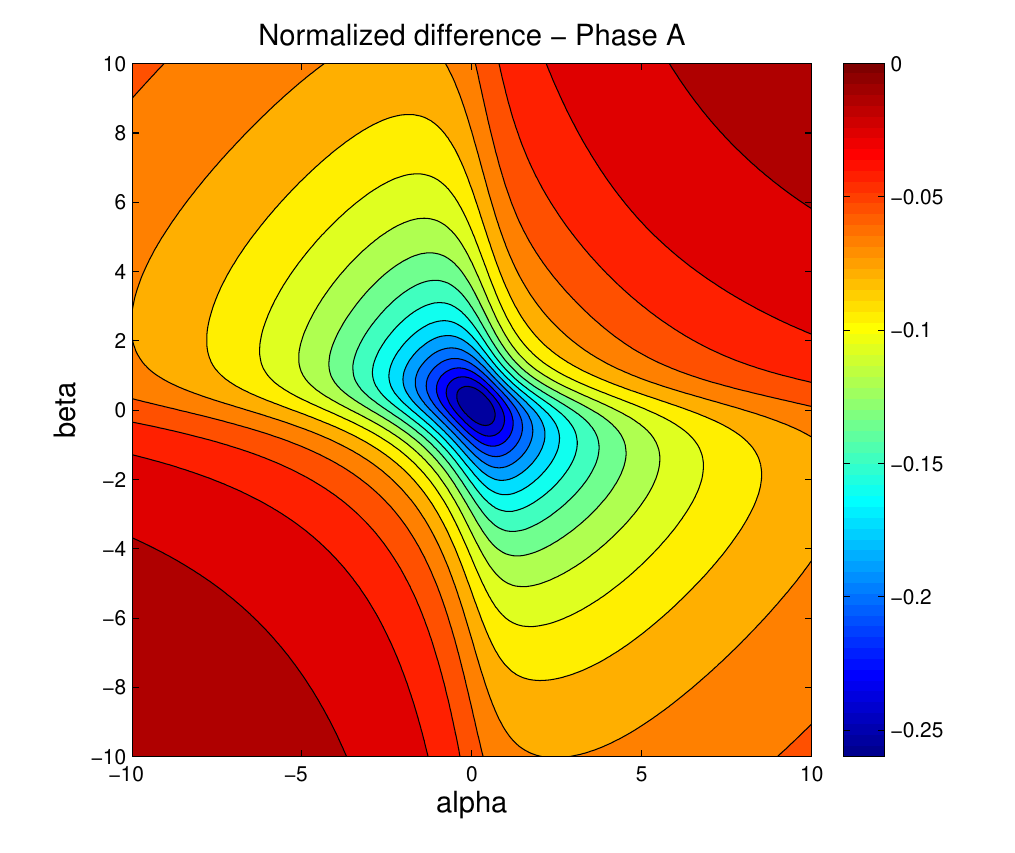}
\includegraphics[width=3.9cm]{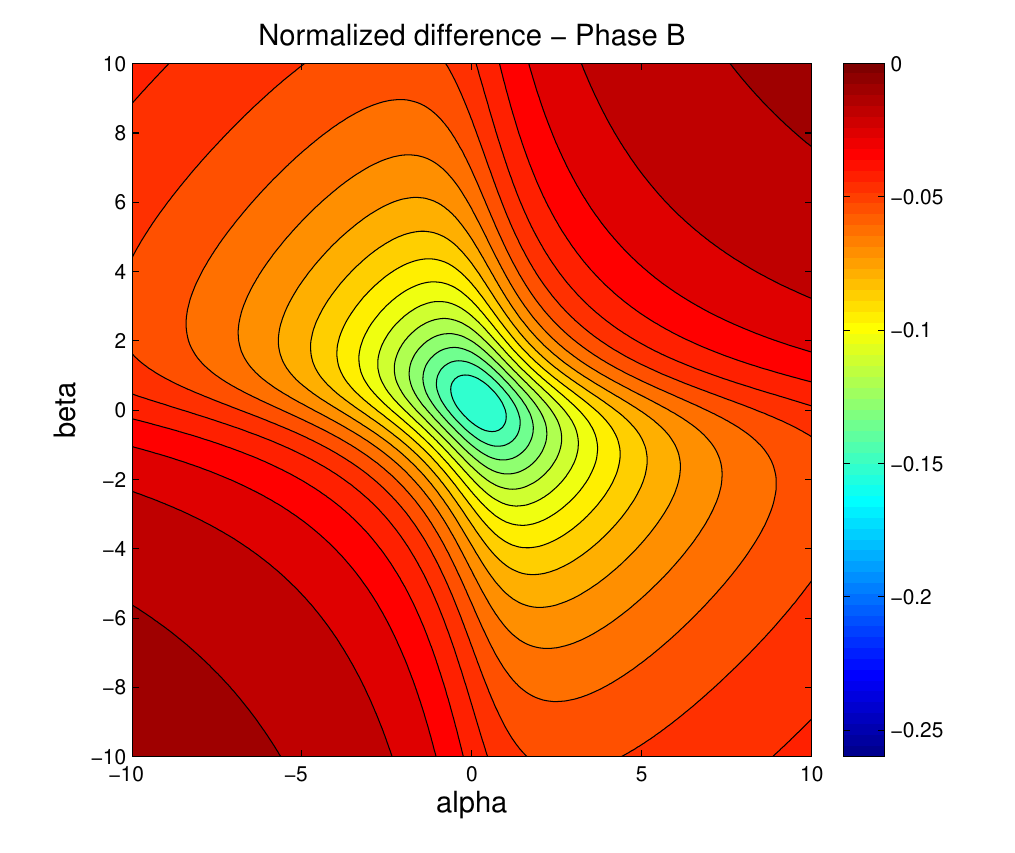}
\includegraphics[width=3.9cm]{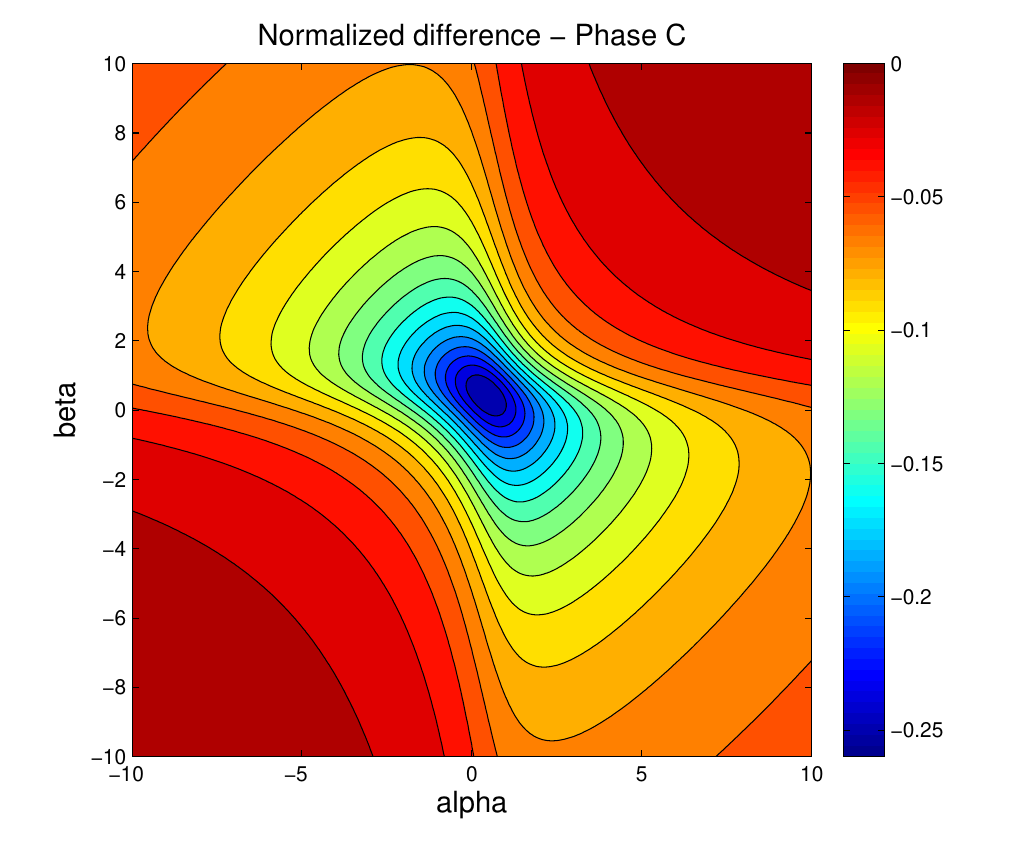}
\caption{Maps of relative increments $\frac{D_{\alpha}-D_{\beta}}{\alpha-\beta}$ for the three subperiods (from left to right: phases A, B and C). The multifractal complexity structural pattern for phase B notably differs from the corresponding patterns for phases A and C, which are relatively similar.}
\label{figure:Hierro-generalized-dimensions-relative-increments}
\end{figure}

\begin{figure}
\includegraphics[width=11cm]{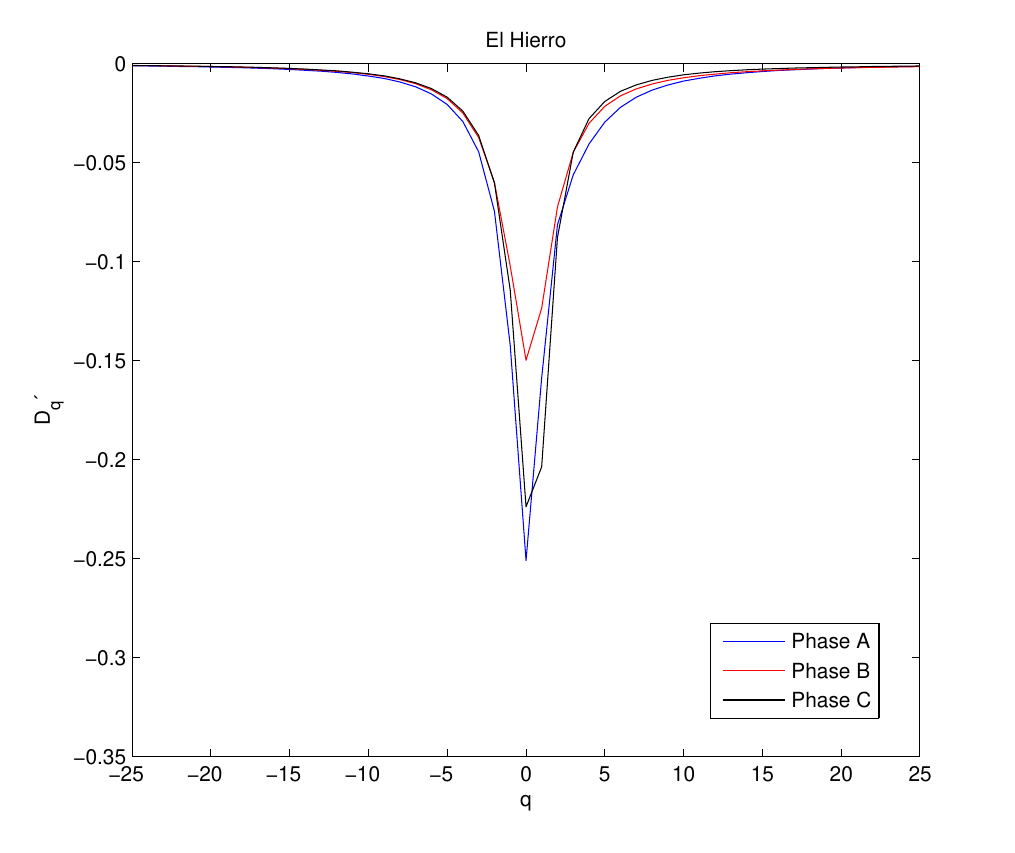}
\centering
\caption{Derivatives of the generalized Rényi dimension curves for the three subperiods (phases A, B and C). The difference between phase B and the two phases A and C is enhanced regarding the corresponding minimum derivative values.}
\label{figure:Hierro-generalized-dimensions-derivatives}
\end{figure}

In the first stage, the analysis only takes into account the temporal distribution of events for each subperiod. Its multifractal nature is reflected in the respective generalized dimension curves, as shown in Figure \ref{figure:Hierro-generalized-dimensions}. There is evidence of an abrupt change from phase A to phase B, with a significant shortening of the range (multifractal step), which tends to be reversed from phase B to phase C. The map of relative increments and, in particular, its limit diagonal representing the derivatives of the generalized dimension curve, respectively displayed in Figures \ref{figure:Hierro-generalized-dimensions-relative-increments} and \ref{figure:Hierro-generalized-dimensions-derivatives}, clearly discriminate central phase B from phases A and C, without much relative  difference between the two latter.

Now, observing that the mean levels, the ranges,  and likely the general patterns of the magnitudes represented in Figure \ref{figure:Hierro-epicenters-magnitudes} (bottom plot) appear to be different, particularly from phase A to phases B and C, we include in the analysis jointly the temporal location and magnitude, paired for each event. As a starting point, it is expected that the physical association between these two aspects, under a relatively regular regime, becomes altered as the dynamics becomes more unstable. For assessment, we then consider, for each partition (with interval length $\varepsilon$) of the temporal domain, two distributions: $\bar{P}_{\varepsilon}$ (the frequency distribution), obtained from box-counting of events, and $\bar{E}_{\varepsilon}$ (the accumulated energy distribution), obtained from `weighted box counting', where each event is weighted by the  amount of released energy according to its magnitude (see Angulo and Esquivel 2014, Esquivel and Angulo 2015, for details). The degree of (multifractal) local coherence between the underlying multifractal measures, respectively denoted  $\bar{p}$ and  $\bar{e}$, within each subperiod, is quantified in terms of the corresponding generalized relative dimensions. Since these are non-symmetric, we calculate $D_{q}(\bar{p}\|\bar{e})$, $D_{q}(\bar{e}\|\bar{p})$, and their arithmetic mean (among other available symmetrization options) for a visual reference.

\begin{figure}
\includegraphics[width=3.9cm]{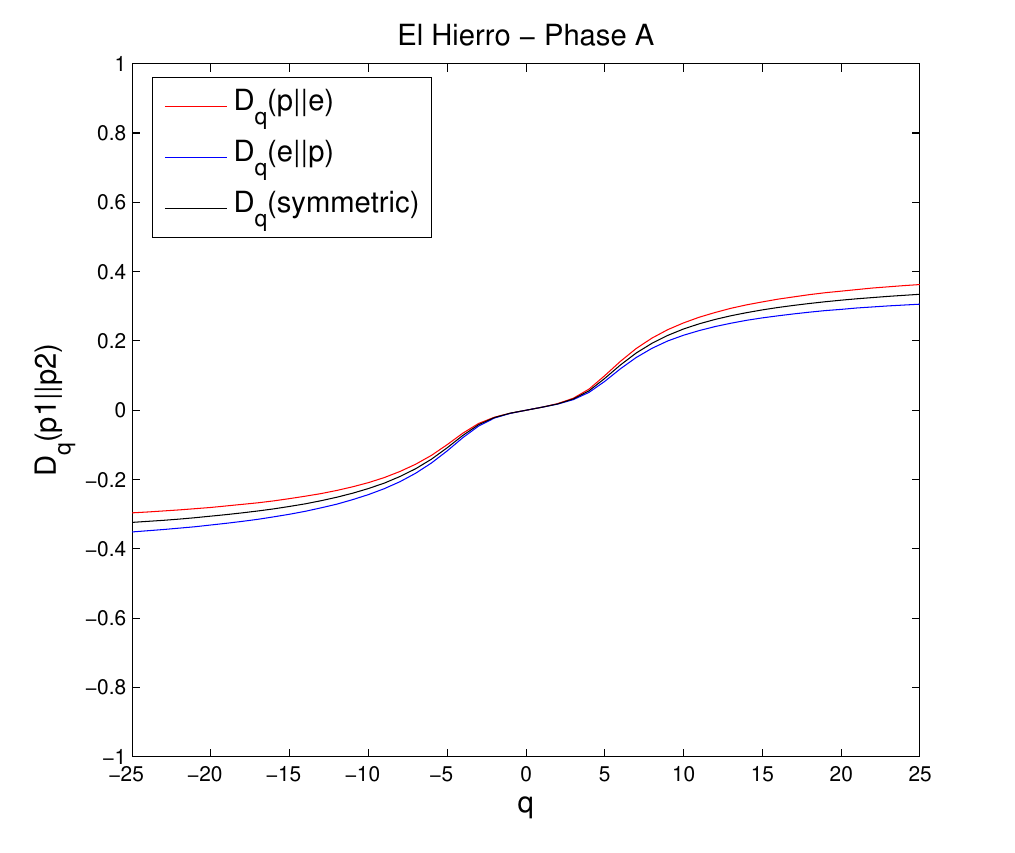}
\includegraphics[width=3.9cm]{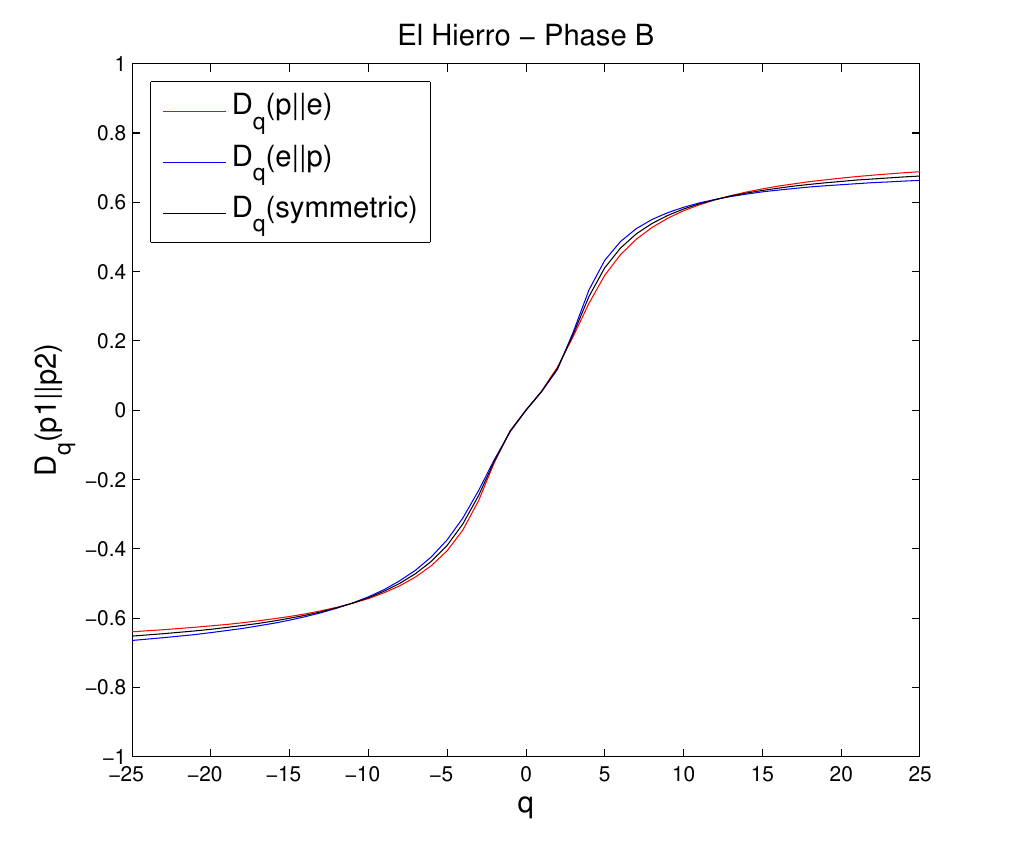}
\includegraphics[width=3.9cm]{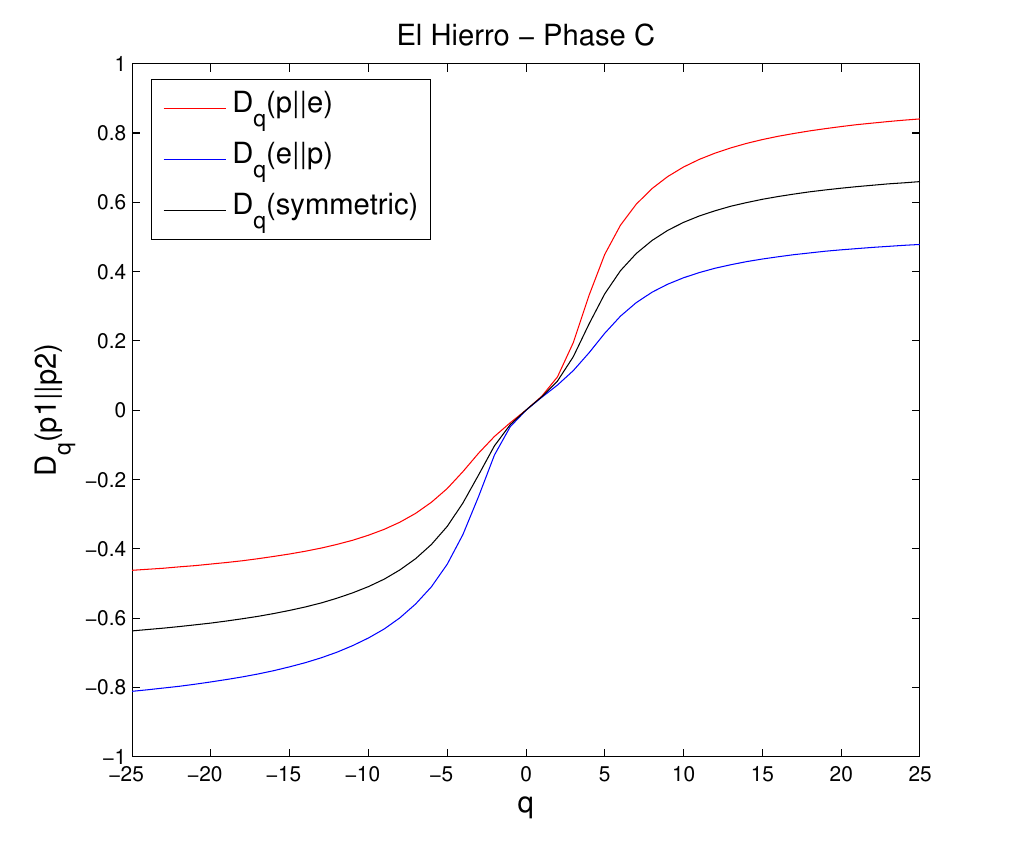}
\caption{Generalized relative Rényi dimension curves for the three subperiods (from left to right: phases A, B and C), based on  $\bar{p}$ as the frequency distribution and $\bar{e}$ as the accumulated energy distribution. A higher dissociation between both distributions is observed in phases B and C with respect to phase A, with a more significant asymmetry in phase C.  }
\label{figure:Hierro-generalized-relative-dimensions}
\end{figure}

\begin{figure}
\includegraphics[width=3.9cm]{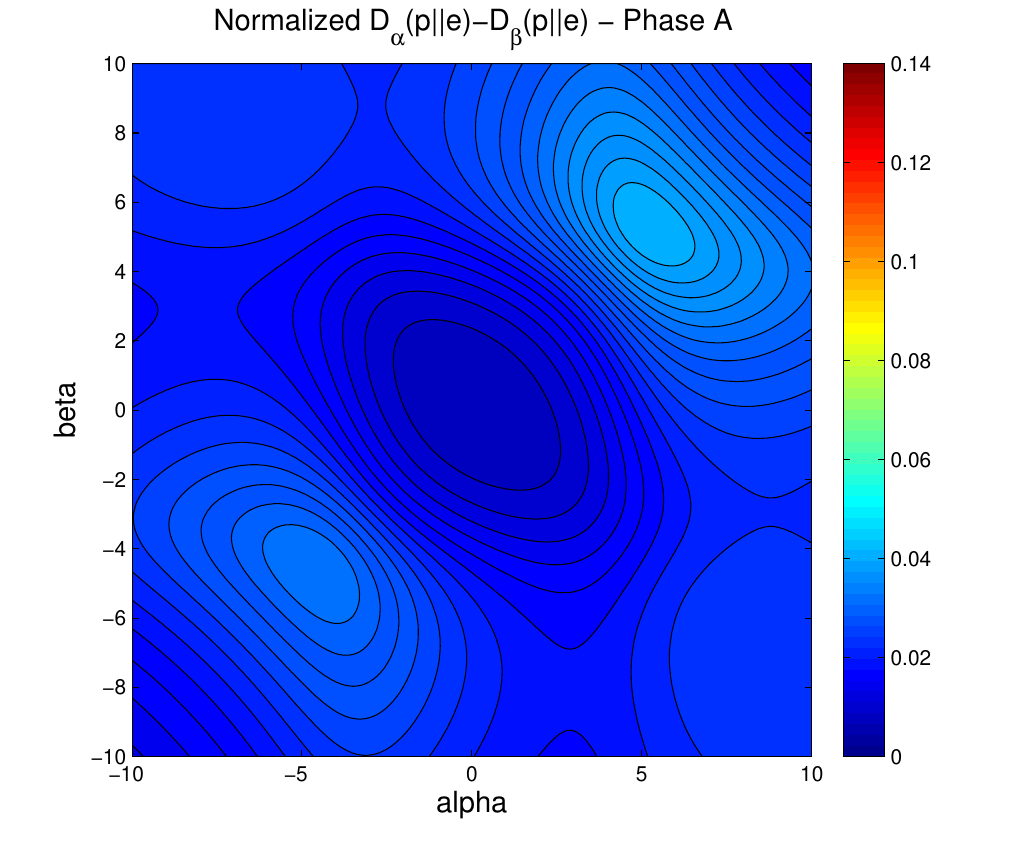}
\includegraphics[width=3.9cm]{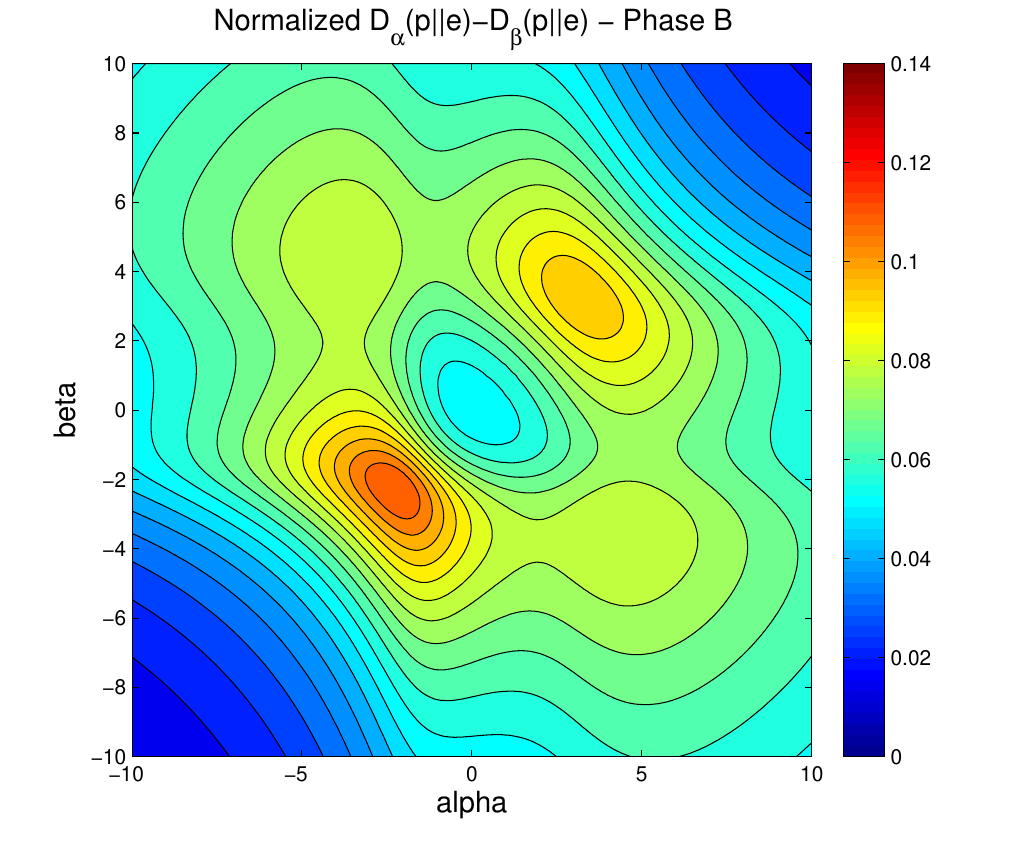}
\includegraphics[width=3.9cm]{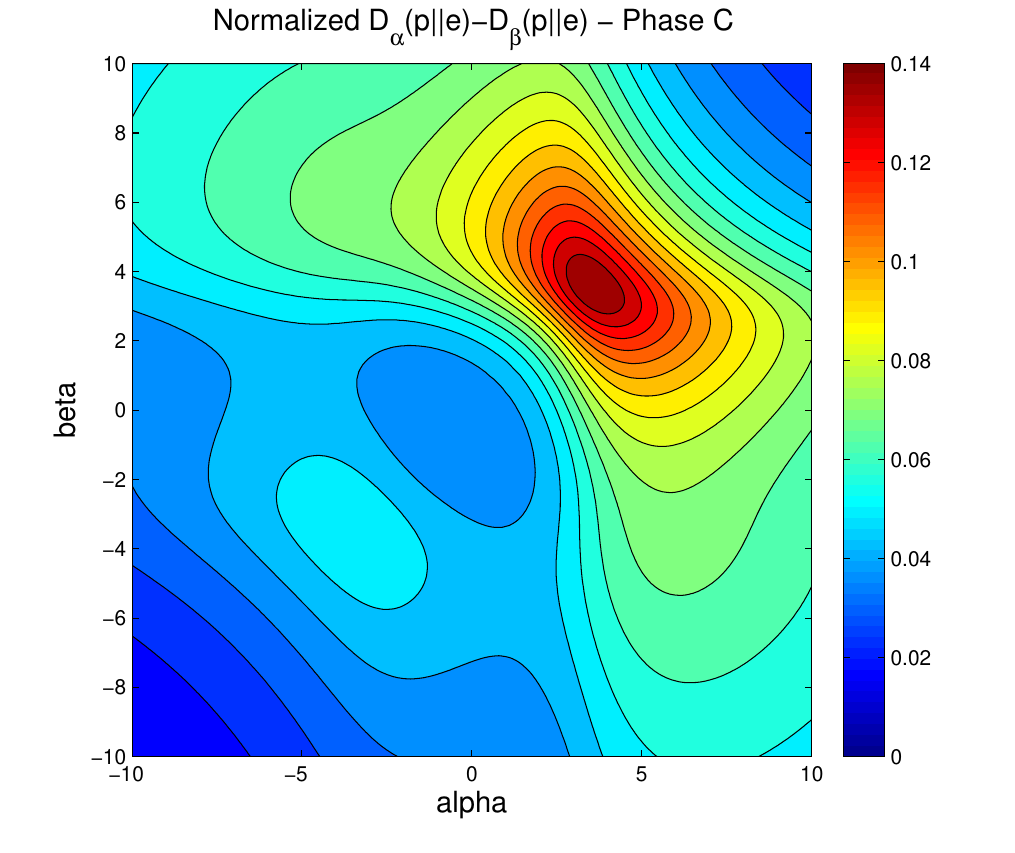}\vspace{3mm}
\\
\includegraphics[width=3.9cm]{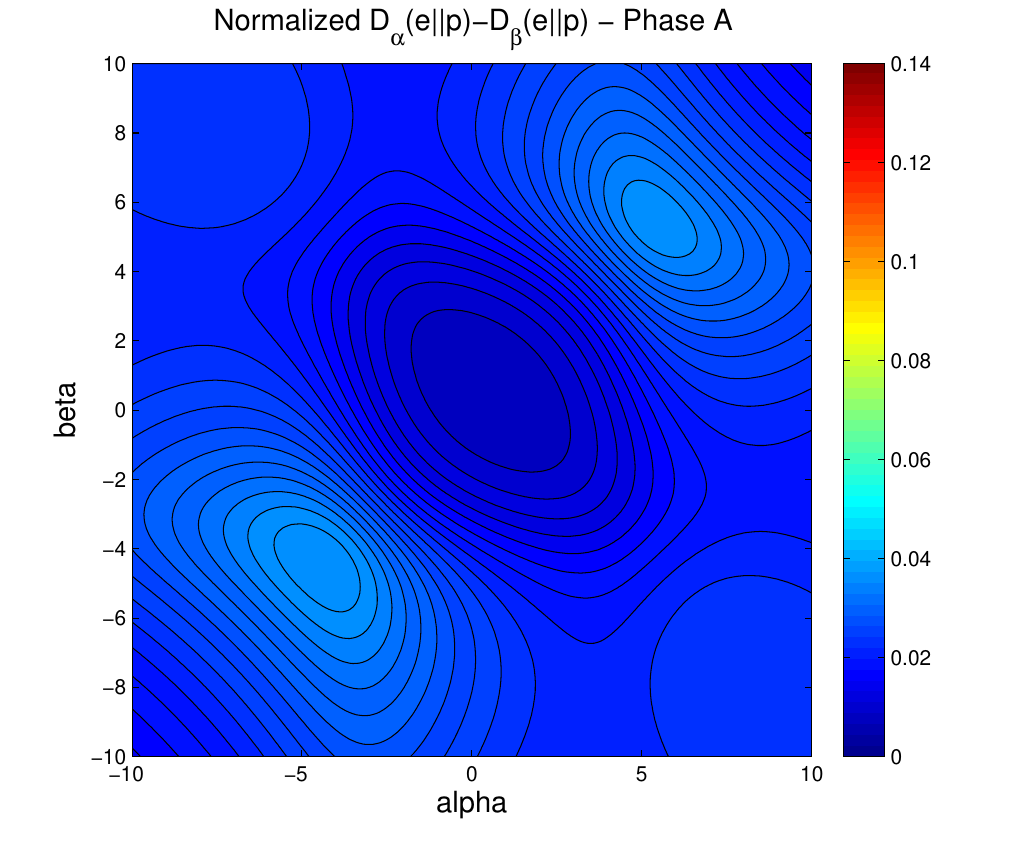}
\includegraphics[width=3.9cm]{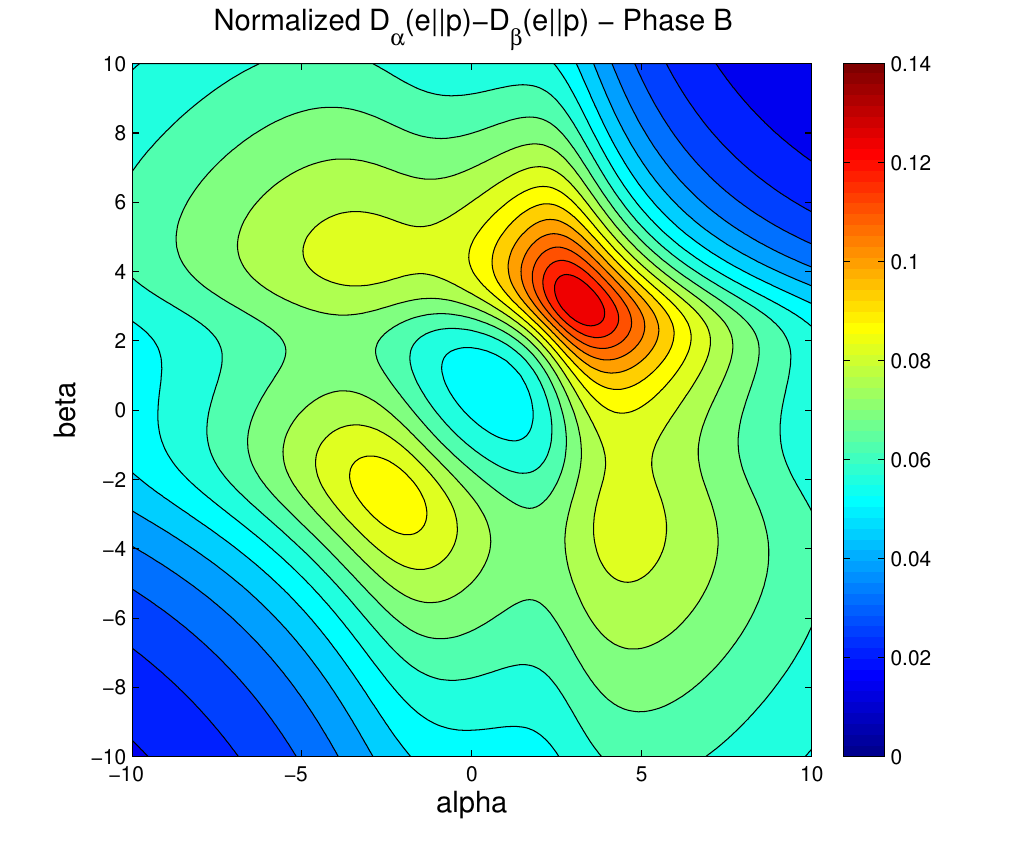}
\includegraphics[width=3.9cm]{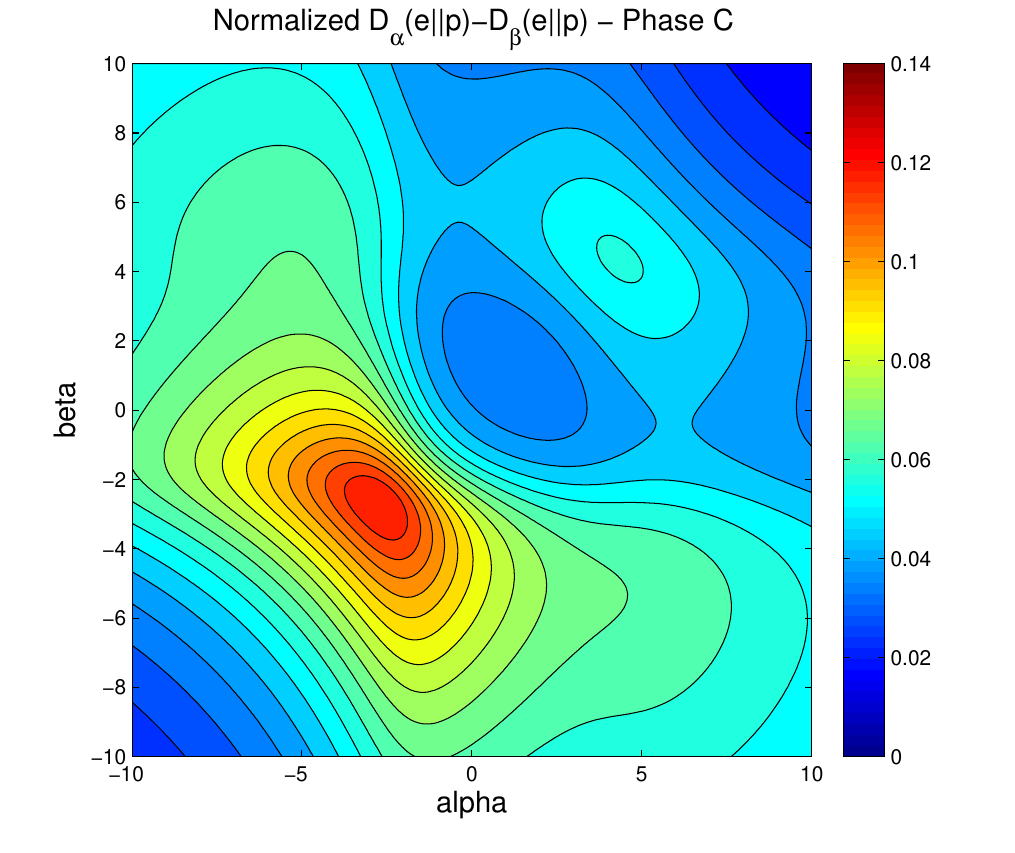}
\caption{Values of $\frac{D_{\alpha}(\overline{p}||\overline{e})-D_{\beta}(\overline{p}||\overline{e})}{\alpha-\beta}$ (top) and $\frac{D_{\alpha}(\overline{e}||\overline{p})-D_{\beta}(\overline{e}||\overline{p})}{\alpha-\beta}$ (bottom), for the three subperiods (from left to right: phases A, B and C). Here, clear dissimilarities between the three patterns reveal the different multifractal complexity structures for phases A, B and C in relation to the joint assessment of the temporal distribution and magnitude of events.}
\label{figure:Hierro-generalized-relative-dimensions-relative-increments}
\end{figure}

\begin{figure}
\includegraphics[height=3.3cm]{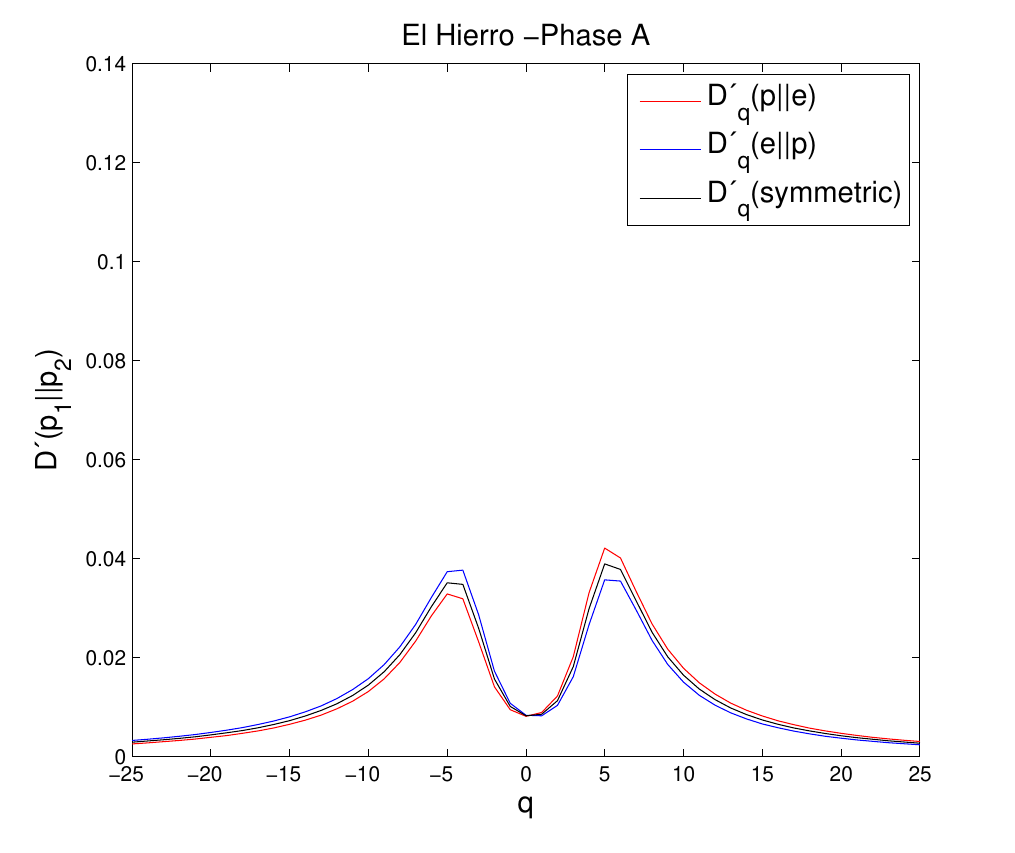}
\includegraphics[height=3.3cm]{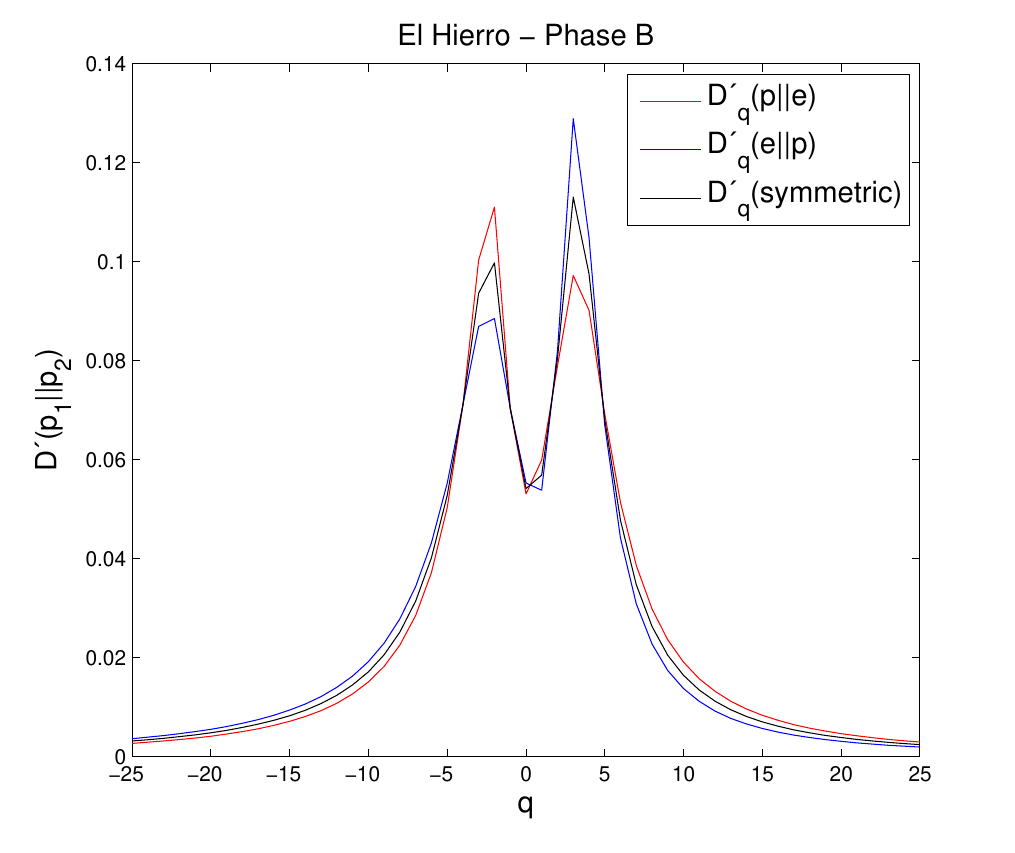}
\includegraphics[height=3.3cm]{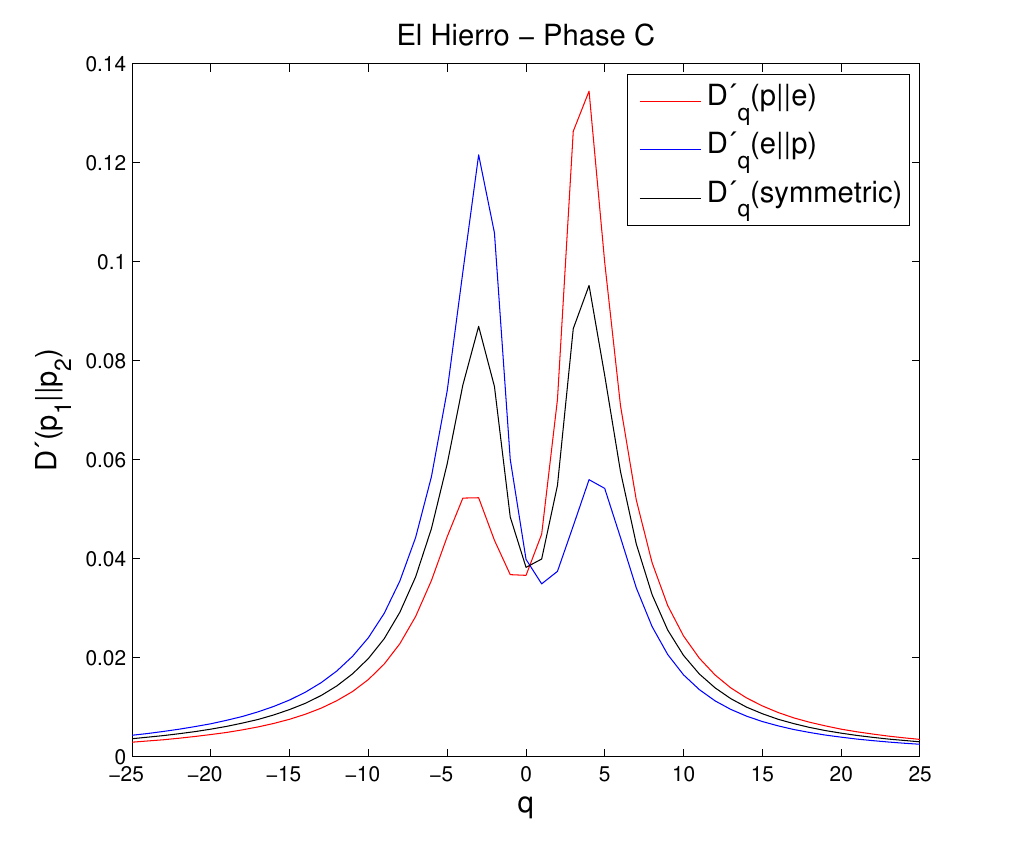}
\caption{Derivatives of the generalized relative Rényi dimension curves for the three subperiods (from left to right: phases A, B and C). Multifractal complexity structural differences, in relation to the joint assessment of the temporal distribution and magnitude of events, are reflected in various features. }
\label{figure:Hierro-generalized-relative-dimensions-derivatives}
\end{figure}

In Figure \ref{figure:Hierro-generalized-relative-dimensions}, we can observe that the generalized relative dimension curves significantly differ between the three subperiods. First, and noting that the curves would be constantly equal to level 0 only in the (unrealistic)  case where $\bar{p} \equiv \bar{e}$, it is clear that there is a higher dissociation between the distributions of events and their magnitudes in phases B and C with respect to phase A. Further comparing phases B and C, we can see that there is a much higher degree of asymmetry, with $D_{q}(\bar{p}\|\bar{e})$  above $D_{q}(\bar{e}\|\bar{p})$, in the latter, which, according to local assessment intrinsic to the definition of Rényi divergence, can be interpreted in the sense that, in this phase, temporal concentrations of events of relatively lower magnitude appear to be predominant in contrast with temporal concentrations of released energy from a lower number of events. These effects are enhanced in the corresponding maps of relative increments shown in Figure \ref{figure:Hierro-generalized-relative-dimensions-relative-increments} and, in particular, in the curves  of derivatives of the generalized relative dimensions displayed in Figure \ref{figure:Hierro-generalized-relative-dimensions-derivatives}, where a distinguishable inversion regarding the levels of the curves $D_{q}(\bar{p}\|\bar{e})$ and $D_{q}(\bar{e}\|\bar{p})$ for the approximate interval $q \in (2,5)$ (and similarly in the negative range) is revealed in the central phase B.

\section{Conclusion and Further Remarks}
\label{section:further-remarks-conclusion}

A synthetic perspective focused on the fundamental role of Shannon and Rényi entropies, as well as Kullback-Leibler and Rényi divergence, as a basis for information and complexity analysis of spatial data, possibly of a multifractal nature, is presented. In each step, the complementary interpretation of entropy-related and divergence-related measures for `global' and `local' structural assessment of the distributions involved is emphasized. 

Well-known product-type generalized complexity measures (based on Rényi entropy) and, being formalized under a parallel conception, product-type generalized relative complexity measures (based on Rényi divergence), originally arisen  in the context of Physics under the notion of complexity as a balance of departure from equilibrium and degeneracy, are reinterpreted in terms of diversity and relative diversity indices. A former proposal of complexity uniquely as divergence from equiprobability is formally identified as a particular case under this approach.  

Parallel limiting connections of the mentioned two-parameter complexity and relative complexity measures with corresponding increments of generalized Rényi dimensions and generalized relative Rényi dimensions with respect to the deformation parameter are highlighted, thus justifying the practical usefulness of related tools for multifractal complexity assessment. These include, in particular, incremental maps and derivative curves, as illustrated in the analysis of a real seismic series, where, among other aspects, the degree of association between the temporal distribution of events and their magnitudes is evaluated on different subperiods for assessment of structural dynamics changes.   

Related complementary aspects and alternative developments in this context are mentioned below. 

Under the original notion and basic formalization of product-type complexity measures for a joint assessment of information and disequilibrium, Martin, Plastino and Rosso (2006) studied different options for the choice of appropriate factors; in particular, they considered the non-extensive Tsallis  entropy and divergence measures (Tsallis 1988), which have also attracted singular attention in diverse fields of application. Angulo and Esquivel (2014) proposed a Tsallis-entropy-based version of multifractal generalized dimensions, used in Angulo and Esquivel (2015) for the formulation of mutual-information-related generalized dependence coefficients in the multifractal domain. In another direction, Alonso, Bueso and Angulo (2016) defined a concept of `generalized mutual complexity' based on generalized Rényi relative complexity measures for structural dependence assessment in a random vector, applied, in particular, to the formulation of a complexity-related optimality criterion for sampling network design. 

Some approaches to `spatial entropy' based on the decomposition of information, explicitly analysing  heterogeneity and dependence in relation to the underlying spatial configuration, have been also introduced in the literature; among related references, see  O'Neill et al. (1988), Riiters, O'Neill, Wickman and Jones (1996), Li and Reynolds (1993), Karlstr\"om and Ceccato (2002), Leibovici (2009), Leibovici, Bastin and Jackson (2011), Leibovici, Claramunt, LeGuyader and Brosset (2014), Altieri, Cocchi and Roli (2017, 2019), etc.

\section*{Acknowledgements}

The authors are grateful to the editors and reviewers for their valuable comments and suggestions to improve the original manuscript. This work was supported by MCIU/AEI/FEDER, UE grant PGC2018-098860-B-I00, and by grant A-FQM-345-UGR18 cofinanced by FEDER Operational Programme 2014-2020 and the Economy and Knowledge Council of the Regional Government of Andalusia, Spain.


\begin{thebibliography}{}


\bibitem{ABA2016} 
Alonso FJ, Bueso MC and Angulo JM (2016). Dependence assessment based on generalized complexity: Application to sampling network design. {\em Methodology and Computing in Applied Probability} 18:921--933

\bibitem{ARC2017} 
Altieri L, Cocchi D and Roli G (2017). A new approach to spatial entropy measures. {\em Environmental and Ecological Statistics} 25: 95--110

\bibitem{ARC2019} 
Altieri L, Cocchi D and Roli G (2019). Advances in spatial entropy measures. {\em Stochastic Environmental Research and Risk Assessment} 33: 1223--1240

\bibitem{AE2014} 
Angulo JM and Esquivel FJ (2014). Structural  complexity in space-time seismic event data. {\em Stochastic Environmental Research and Risk Assessment} 28:1187--1206

\bibitem{AE2015} 
Angulo JM and Esquivel FJ (2015). Multifractal dimensional dependence assessment based on Tsallis mutual information. {\em Entropy} 17:5382--5401

\bibitem{B1974} 
Batty M (1974). Spatial entropy. {\em Geographical Analysis} 6:1--31

\bibitem{BMMS2014} 
Batty M, Morphet R, Masucci P and Stanilov K (2014). Entropy, complexity, and spatial information. {\em Journal of Geographical Systems} 16:363--385

\bibitem{C1966} 
Campbell LL (1966). Exponential entropy as a measure of extent of a distribution. {\em Zeitschrift f\"ur Wahrscheinlichkeitstheorie und Verwandte Gebiete} 5:217--225

\bibitem{CGL-R2002} 
Catal\'an RG, Garay J and L\'opez-Ruiz R (2002). Features of the extension of a statistical measure of complexity to continuous systems. {\em Physical Review E} 66:011102
  
\bibitem{EAA2017} 
Esquivel FJ, Alonso FJ and Angulo JM (2017). Multifractal complexity analysis in space–time based on the generalized dimensions derivatives. {\em Spatial Statistics} 22:469--480
 
\bibitem{EA2015} 
Esquivel FJ and Angulo JM (2015). Non-extensive analysis of the seismic activity involving the 2011 volcanic eruption in El Hierro. {\em Spatial Statistics} 14(B):208--221

\bibitem{G1983} 
Grassberger P (1983). Generalized dimensions of strange attractors. {Physics Letters A} 97:227--230 

\bibitem{H2001} 
Harte D (2001). {\em Multifractals: Theory and Applications}. CRC Press: Boca Raton

\bibitem{H1928} 
Hartley RV (1928). Transmission of information. {\em The Bell System Technical Journal} 7:535--563

\bibitem{HP1983} 
Hentschel HGE and Procaccia I (1983). The infinite number of generalized dimensions of fractals and strange attractors. {\em Physica D} 1983:435--444

\bibitem{HH1986} 
Huberman BA and Hogg T (1986). Complexity and adaptation. {\em Physica D} 22:376--384

\bibitem{KC2002} 
Karlstr\"om A and Cecatto V (2002). A new information theoretical measure of global and local spatial association. {\em Jaharb Regionalwissensc} 22:13--40

\bibitem{KL1951} 
Kullback S and Leibler RA (1951). On information and sufficiency. {\em Annals of Mahematical Statistics} 22:79--86

\bibitem{L2009} 
Leibovici DG (2009). Defining spatial entropy from multivariate distributions of co-occurrences. In {\em  Spatial information theory. COSIT 2009. Lecture Notes in Computer Science}, Hornsby KS, Claramunt C, Denis M, Ligozat G (eds), vol. 5756:392--404. Springer: Berlin, Heidelberg

\bibitem{LBJ2013} 
Leibovici DG, Bastin L and Jackson M (2011). Higher-order co-occurrences for exploratory point pattern analysis and decision tree clustering on spatial data. {\em Computers \& Geosciences}
37:382--389

\bibitem{LCLB2014} 
Leibovici DG, Claramunt C, Le Guyader C and Brosset D (2014). 
 Local and global spatio-temporal entropy indices based on distance ratios and co-occurrences distributions. {\em  International Journal of Geographical Information Science} 28:1061--1084
 
\bibitem{LR1993} 
Li H and Reynolds JF (1993). A new contagion index to quantify spatial
patterns of landscapes. {\em Landscape Ecology} 8:155--162 

\bibitem{L-RMC1995} 
L\'opez-Ruiz R, Mancini HL and Calbet X (1995). A statistical measure of complexity. {\em Physics Letters A} 209:321--326

\bibitem{L-RNRS2009} 
L\'opez-Ruiz R, Nagy A, Romera E and Sa\~nudo J (2009). A generalized statistical complexity measure: Applications to quantum systems.
{\em Journal of Mathematical Physics} 50:123528 

\bibitem{MPR2006} 
Martin MT, Plastino O and  Rosso A (2006). Generalized statistical complexity measures: geometrical and analytical properties.
{\em Physica A} 369:439--462

\bibitem{OEtAl1988} 
O'Neill RV, Krummel JR, Gardner RH, Sugihara G, Jackson B, DeAngelis DL, Milne BT, Turner MG, Zygmunt B, Christensen SW, Dale VH and Graham RL (1988). Indices of landscape pattern. {\em Landscape Ecology} 1:153--162

\bibitem{R1961} 
R\'{e}nyi A (1961). On measures of entropy and information. In {\em Proceedings of the Fourth Berkeley Symposium on Mathematical Statistics and Probability}, Berkeley, CA, USA, June 20--July 30 1960; Neyman, J (ed) vol. 1:547--561. University of California Press: Berkeley

\bibitem{ROWJ996} 
Riitters KH, O'Neill RV, Wickham JD and Jones KB (1996). A note on
contagion indices for landscape analysis. {\em Landscape Ecology} 11:197--202

\bibitem{RSN2011} 
Romera E, Sen KD and Nagy \'A (2011). A generalized relative complexity measure. {\em Journal of Statistical Mechanics: Theory and Experiment} 09:P09016

\bibitem{S1948} 
Shannon CE (1948). A mathematical theory of communication. {\em The Bell System Technical Journal} 27:379--423

\bibitem{T1967}
Theil H (1967). {\em Economics and Information Theory}. North Nolland: Amsterdam 


\bibitem{T1972} 
Theil H (1972). {\em Statistical Decomposition Analysis}. North Holland:
Amsterdam

\bibitem{T1988} 
Tsallis C (1988). Possible generalization of Boltzmann-Gibbs statistics. {\em Journal of Statistical Physics} 52:479--487


\end{thebibliography}
\end{document}